\documentclass[twoside]{siamltex}
\pdfoutput=1
\usepackage{amsmath}
\usepackage{amssymb}
\usepackage{amsfonts}
\usepackage[]{hyperref}
\usepackage{xcolor}
\hypersetup{
 colorlinks, 
 linkcolor= {red!80!black},  
 citecolor={green!80!black},
 urlcolor={blue!50!black},
}

\usepackage{graphicx}
\usepackage{subcaption}
\newtheorem{rem}[theorem]{Remark}
\newtheorem{assu}[theorem]{Assumption}

\renewcommand{\O}{\Omega}
\newcommand{\Si}{\Sigma}
\newcommand{\pO}{\partial\O}
\newcommand{\GD}{{\Gamma_D}}
\newcommand{\GN}{{\Gamma_N}}
\newcommand{\ad}{{\text{ad}}}

\newcommand{\field}[1]{\mathbb{#1}}
\newcommand{\R}{\field{R}}

\renewcommand{\AA}{\mathcal{A}}
\newcommand{\BB}{\mathcal{B}}
\newcommand{\DD}{\mathcal{D}}
\newcommand{\EE}{\mathcal{E}}
\newcommand{\HH}{\mathcal{H}}
\newcommand{\JJ}{\mathcal{J}}
\newcommand{\KK}{\mathcal{K}}
\newcommand{\LL}{\mathcal{L}}
\newcommand{\fL}{\mathfrak{L}}
\newcommand{\MM}{\mathcal{M}}

\newcommand{\s}{\mathcal{S}}
\newcommand{\UU}{\mathcal{U}}
\newcommand{\U}{\mathbb{U}}
\newcommand{\embeds}{\hookrightarrow}
\newcommand{\olom}{\overline{\Omega}}
\newcommand{\olQ}{\overline{Q}}
\newcommand{\olJ}{\overline{J}}
\newcommand{\dd}{\mathrm{d}}
\newcommand{\refA}{\textup{\textbf{(A)}}}
\newcommand{\refS}{\textup{\textbf{(S)}}}
\newcommand{\X}{\mathbb{X}}

\newcommand{\mfr}{\mathfrak{r}}

\newenvironment{proofof}[1]{\emph{Proof of #1.}}{\endproof}

\DeclareMathOperator{\dom}{dom}    
\DeclareMathOperator{\dive}{div}   

\DeclareMathOperator*{\esssup}{ess\,sup}
\DeclareMathOperator*{\essinf}{ess\,inf}

\makeatletter
\def\romlabel#1{\begingroup
    \def\@currentlabel{{\upshape(\roman{rmnum})}}%
    \label{#1}\endgroup
}
\makeatother 

\begin{document}

\author{ H. Meinlschmidt\footnotemark[2], C.~Meyer\footnotemark[3],
  J.~Rehberg\footnotemark[4]}

\title{Optimal control of the thermistor problem in three spatial dimensions}

\renewcommand{\thefootnote}{\fnsymbol{footnote}}

\footnotetext[2]{Faculty of Mathematics, TU Darmstadt, Dolivostrasse
  15, D-64283 Darmstadt, Germany.}  \footnotetext[3]{Faculty of
  Mathematics, Technical University of Dortmund, Vogelpothsweg 87,
  D-44227 Dortmund, Germany.}  \footnotetext[4]{Weierstrass Institute
  for Applied Analysis and Stochastics Mohrenstr.~39, D-10117 Berlin,
  Germany.}

\renewcommand{\thefootnote}{\arabic{footnote}}
\pagestyle{myheadings} \thispagestyle{plain}
\markboth{H.~MEINLSCHMIDT, C.~MEYER, J.~REHBERG} {OPTIMAL CONTROL OF
  THE THERMISTOR PROBLEM}
 
\maketitle

\begin{abstract}
  This paper is concerned with the state-constrained optimal control
  of the three-dimensional thermistor problem, a fully quasilinear
  coupled system of a parabolic and elliptic PDE with mixed boundary
  conditions. This system models the heating of a conducting material
  by means of direct current. Local existence, uniqueness and
  continuity for the state system are derived by employing maximal
  parabolic regularity in the fundamental theorem of Pr\"uss. Global
  solutions are addressed, which includes analysis of the linearized
  state system via maximal parabolic regularity, and existence of
  optimal controls is shown if the temperature gradient is under
  control. The adjoint system involving measures is investigated using
  a duality argument. These results allow to derive first-order
  necessary conditions for the optimal control problem in form of a 
  qualified optimality system. The theoretical findings are illustrated by 
  numerical results.
\end{abstract}

\begin{keywords}
  Partial differential equations, optimal control problems, state
  constraints
\end{keywords}

\begin{AMS}
  35K59, 35M10, 49J20, 49K20
\end{AMS}


\section[Introduction] {Introduction}
\label{s-Introduction}

In this paper, we consider the state-constrained optimal control of
the three-dimensional thermistor problem. In detail the optimal
control problem under consideration looks as follows:
\begin{equation}\tag{P}\label{P}
  \left.
    \begin{aligned}
      \min \quad & 
      \frac{1}{2}\|\theta(T_1)\!-\!\theta_d\|_{L^2(E)}^2\! +\!
      \frac{\gamma}{s} \|\nabla\theta\|_{L^{s}(T_0,T_1;L^q(\Omega))}^{s}\! +\!
      \frac{\beta}{2} \int_{\Si_N} (\partial_t u)^2\! +\! |u|^p \, \dd
      \omega \,\dd t
      \\[0.1cm]
      \text{s.t.}  \quad
      & \text{\eqref{e-strongparabolic}--\eqref{e-strongellipticdirichlet}}\\[1em]
      \text{and} \quad &
      \begin{aligned}[t]
        \theta(x,t) &\leq \theta_{\max}(x,t)& &\text{a.e.~in }\O
        \times (T_0,T_1),
        \\ 0\leq u(x,t) &\leq u_{\max}(x,t)& &\text{a.e.~on }\Gamma_N \times (T_0, T_1)
      \end{aligned}
    \end{aligned}
    \enspace
  \right\}
\end{equation}
where~\eqref{e-strongparabolic}--\eqref{e-strongellipticdirichlet}
refer to the following coupled PDE system consisting of the
instationary nonlinear heat equation and the quasi-static potential equation,
which is also known as {\em thermistor problem}:
\begin{align}
  \partial_t\theta - \dive(\eta(\theta)\kappa\nabla \theta) &=
  (\sigma(\theta)\varepsilon \nabla \varphi)\cdot \nabla \varphi
  & &\text{in }Q := \O\times (T_0,T_1)\label{e-strongparabolic}\\
  \nu \cdot \kappa \nabla \theta + \alpha\theta &= \alpha\theta_l &
  &\text{on }\Si:=\pO\times (T_0,T_1)\label{e-strongparabolicrobin}\\
  \theta(T_0) &= \theta_0 & &\text{in }\O\label{e-initial}\\[1em]
  -\dive(\sigma(\theta)\varepsilon \nabla \varphi) &= 0& &\text{in
  }Q
  \\
  \nu \cdot \sigma(\theta) \varepsilon \nabla \varphi &= u& &\text{on }\Si_N:=\GN\times
 (T_0,T_1)\label{e-strongellipticneumann}\\
  \varphi&=0& &\text{on }\Si_{D}:=\GD\times (T_0,T_1).
  \label{e-strongellipticdirichlet} 
\end{align}
%
Here $\theta$ is the temperature in a conducting material covered by
the three dimensional domain $\O$, while $\varphi$ refers to the
electric potential. The boundary of $\O$ is denoted by $\pO$ with the
unit normal $\nu$ facing outward of $\O$ in almost every boundary
point (w.r.t.\ the boundary measure $\omega$). In
addition, for the boundary we have $\GD~\dot{\cup}~\GN = \pO$, where
$\GD$ is closed within $\pO$. The functions $\eta(\cdot)\kappa$ and
$\sigma(\cdot)\varepsilon$ represent heat- and electric conductivity.
While $\kappa$ and $\varepsilon$ are given, prescribed functions,
$\eta$ and $\sigma$ are allowed to depend on the temperature
$\theta$. Moreover, $\alpha$ is the heat transfer coefficient and
$\theta_l$ and $\theta_0$ are given boundary-- and initial data,
respectively. Finally, $u$ stands for a current which is induced via
the boundary part $\GN$ and is to be controlled. The bounds in the
optimization problem~\eqref{P} as well as the desired temperature
$\theta_d$ are given functions and $\beta$ is the usual Tikhonov
regularization parameter.  The precise assumptions on the data
in~\eqref{P}
and~\eqref{e-strongparabolic}--\eqref{e-strongellipticdirichlet} will
be specified in \S\ref{s-not}. In all what follows, the
system~\eqref{e-strongparabolic}--\eqref{e-strongellipticdirichlet} is
frequently also called \emph{state system}.

The PDE
system~\eqref{e-strongparabolic}--\eqref{e-strongellipticdirichlet}
models the heating of a conducting material by means of a direct
current, described by $u$, induced on the part $\GN$ of the boundary,
which is done for some time $T_1-T_0$. At the grounding $\GD$, homogeneous
Dirichlet boundary conditions are given, i.e., the potential is zero,
inducing electron flow. Note that, usually, $u$ will
be zero on a subset $\Gamma_{N_0}$ of $\GN$, which corresponds to
having insulation at this part of the boundary. We emphasize that the
different boundary conditions are essential for a realistic modeling
of the process. The objective of~\eqref{P} is to adjust the induced
current $u$ to minimize the $L^2$-distance between the desired and the
resulting temperature at end time $T_1$ on the set $E \subseteq \O$, the
latter representing the area of the material in which one is
interested -- realized in the objective functional by the first
term. The other terms are in present to minimize thermal
stresses (second term) and to ensure a certain smoothness of the
controls (third term), whose influence to the objective functional, however, may be
controlled by the weights $\gamma$ and $\beta$. The actual form of
these terms is motivated by functional-analytic considerations, see \S\ref{subs-existenceoc}.
 Moreover, the optimization is subject to
pointwise control and
state constraints. The  control constraints
 reflect a maximum heating power, while the
state constraints limit the temperature evolution to prevent possible
damage, e.g.\ by melting of the material. Similarly to the mixed
boundary conditions, the inequality constraints in~\eqref{P} are
essential for a realistic model as demonstrated by the numerical
example within this paper. Problem~\eqref{P} is relevant in various
applications, such as for instance the heat treatment of steel by
means of an electric current. The example considered in the numerical
part of this paper deals with an application of this type.

The state
system~\eqref{e-strongparabolic}--\eqref{e-strongellipticdirichlet}
exhibits some non-standard features, in particular due to the
quasilinear coupling of the parabolic and the elliptic PDE, the mixed
boundary conditions
in~\eqref{e-strongellipticneumann}--\eqref{e-strongellipticdirichlet},
and the inhomogeneity in the heat equation~\eqref{e-strongparabolic}
as well as the temperature-dependent heat conduction coefficients.
Besides the quasilinear state system, the pointwise state constraints
on the temperature represent another challenging feature of the
optimal control problem under consideration.  The Lagrange multipliers
associated with constraints of this kind only provide poor regularity
in general, which especially complicates the analysis of the adjoint
equation.

We briefly describe the genuine aspects of our work.  First of all,
the discussion of the quasilinear state system alone requires
sophisticated up-to-date tools from maximal elliptic and parabolic
regularity theory.  This concerns already local-in-time existence for
solutions
of~\eqref{e-strongellipticneumann}--\eqref{e-strongellipticdirichlet},
let alone the characterization of global-in-time solutions.  The
corresponding maximal regularity results were established only
recently, see e.g.~\cite{ABHDR, hal/reh, hal/rehPara} for the
parabolic case and~\cite[Appendix]{krum/reh},~\cite{DKR} for the
elliptic one.  Our key ingredient for the proof of local-in-time
existence is a general result of Pr\"uss on quasilinear parabolic
equations~\cite{pruess}.  To verify the assumptions required for the
application of Pr\"uss' result, we heavily rely on an isomorphism
property of the elliptic differential operators in both equations of
the state
system. 
Assuming this isomorphism property only for the case of pure diffusion
coefficients $\kappa$ and $\varepsilon$ in the differential operators,
see Assumption~\ref{a-reg} below, we show that the nonlinear
differential operators involving $\eta(\theta)$ and $\sigma(\theta)$
then also enjoy it, by a technique developed in~\cite{krum/reh}.
However, this analysis only guarantees the local-in-time existence,
and the counterexample in~\cite{AC94} involving a blow-up criterion for a
similar model of the thermistor system demonstrates that one can, in
general, not expect global-in-time solutions. Nevertheless, based on
recent results on non-autonomous parabolic equations~\cite{hajo15}, we
prove that there are control functions that admit global-in-time
solutions.  Moreover, using the implicit function theorem, we show
that these control functions form an open set, which is essential for
the derivation of optimality conditions in qualified form that are
useful for numerical computations.  Concerning the existence of global
minimizers for~\eqref{P}, we benefit from the pointwise state
constraints and the second addend in the objective functional
involving the gradient of the temperature.  Both terms prevent a
blow-up of the temperature and its gradient and allow to restrict the
discussion of the optimization problem to control functions that admit
a global-in-time solution of the state system. This approach is
inspired by~\cite{AQ05}, where a similar technique was used to
establish the existence of optimal controls.

Let us put our work into perspective.  Up to the authors' best
knowledge, there are only few contributions dealing with the optimal
control of the thermistor problem. We refer to~\cite{LS05, cimatti07,
  HLP08}, where two-dimensional problems are discussed.
In~\cite{LS05}, a completely parabolic problem is discussed,
while~\cite{HLP08} considers the purely elliptic counterpart
to~\eqref{e-strongparabolic}--\eqref{e-strongellipticdirichlet}.
In~\cite{cimatti07, SA07}, the authors investigate a
parabolic-elliptic system similar
to~\eqref{e-strongparabolic}--\eqref{e-strongellipticdirichlet},
assuming a particular structure of the controls.  In contrast
to~\cite{LS05, HLP08}, mixed boundary conditions are considered
in~\cite{cimatti07}.  However, all these contributions do not consider
pointwise state constraints and non-smooth data.  Thus,~\eqref{P}
differs significantly from the problems considered in the
aforementioned papers.  In a previous paper~\cite{h/m/r/r}, two of the
authors investigated the two-dimensional counterpart of~\eqref{P}.
This contribution also accounts for mixed boundary conditions,
non-smooth data, and pointwise state constraints.  However, the
analysis in~\cite{h/m/r/r} substantially differs from the three
dimensional case considered here. First of all, in two spatial
dimensions, the isomorphism-property of the elliptic operators
mentioned above directly follows from the classical
paper~\cite{groeger89}. Moreover, the heat conduction coefficient
in~\eqref{e-strongparabolic} is assumed not to depend on the
temperature in~\cite{h/m/r/r}.  Both features allow to derive a global
existence result for a suitable class of control functions.  Hence,
main aspects of the present work do not appear in the two-dimensional
setting.  Let us finally take a broader look on state-constrained
optimal control problems governed by PDEs.  Compared to semilinear
state-constrained optimal control problems, the literature concerning
optimal control problems subject to quasilinear PDEs and pointwise
state constraints is rather scarce. We exemplarily refer
to~\cite{CY95, CT09}, where elliptic problems are studied.  The vast
majority of papers in this field deals with problems that possess a
well defined control-to-state operator. By contrast, as indicated
above, the
state-system~\eqref{e-strongparabolic}--\eqref{e-strongellipticdirichlet}
in general just admits local-in-time solutions, which requires a
sophisticated treatment of the optimal control problem under
consideration.

The paper is organized as follows: We set the stage with notations and
assumptions in \S\ref{s-not} and discuss the state-system in
\S\ref{s-statesys}. More precisely, \S\ref{subs-prereq} collects
preliminary results, also interesting in their own sake, while
\S\ref{subs-proof} is devoted to the actual proof of existence and
uniqueness of local-in-time solutions. We then proceed with the
optimal control problem in \S\ref{s-oc}. Before stating first order
necessary conditions for~\eqref{P} in \S\ref{subs-noc}, we give
sufficient conditions for controls to produce global solutions and
establish continuous differentiability of the control-to-state
operator for global solutions in \S\ref{subs-existenceoc}. The paper
is wrapped-up with an illustrative numerical example in \S\ref{sec:numerics}.


\section{Notations and general assumptions}\label{s-not}

We introduce some notation and the relevant function
spaces. All function spaces under our consideration are \emph{real}
ones. Let, for now, $\Omega$ be a domain in $\R^3$. We give precise
geometric specifications for $\Omega$ in \S\ref{sec:geometric} below.

Let us fix some notations: The underlying time interval is called $J =
(T_0,T_1)$ with $T_0 < T_1$. The boundary measure for the domain
$\O$ is called $\omega$. Generally, given an integrability order $q
\in (1,\infty)$, we denote the conjugated of $q$ by $q'$, i.e., it
always holds $1/q + 1/q' = 1$. 

\begin{definition} \label{d-gebiet} For $q \in (1,\infty)$, let
  $W^{1,q}(\O)$ denote the usual Sobolev space on $\O$. If $\Xi
  \subset \partial \O$ is a closed part of the boundary $\pO$, we set
  $W^{1,q}_\Xi(\O)$ to be the closure of the set $\left\{\psi|_\O:\;
    \psi \in C^\infty_0(\R^3), \, \supp~\psi \cap \Xi = \emptyset
  \right\}$ with respect to the $W^{1,q}$-norm.
\end{definition}

The dual space of
$W^{1,q'}_\Xi(\O)$ is denoted by $W^{-1,q}_\Xi (\O)$; in particular,
we write $W_{\emptyset}^{-1,q}(\O)$ for the dual of $W^{1,q'}(\O)$
(see Remark~\ref{r-konsist} below regarding consistency). The H\"older
spaces of order $\delta$ on $\O$ or order $\varrho$ on $Q$ are denoted
by $C^\delta(\O)$ and $C^\varrho(Q)$, respectively (note here that
H\"older continuous functions on $\O$ or $Q$, respectively, possess an
unique uniformly continuous extension to the closure of the domain,
such that we will mostly use $C^\delta(\olom)$ and
$C^\varrho(\overline Q)$ to emphasize on this).

We will usually abbreviate the function spaces on $\O$ by leaving out
the $\O$, e.g.\ we write $W^{1,q}_\Xi$ instead of $W^{1,q}_\Xi(\O)$ or
$L^p$ instead of $L^p(\O)$. 
Lebesgue spaces on subsets of $\partial \O$
are always to be considered with respect to the boundary measure
$\omega$, but we abbreviate $L^p(\pO,\omega)$ by $L^p(\pO)$ and do so
analogously for any $\omega$-measurable subset of the boundary.  The
norm in a Banach space $X$ will be always indicated by $\|\cdot \|_X$.
For two Banach spaces $X$ and $Y$, we denote the space of linear,
bounded operators from $X$ into $Y$ by $\LL(X;Y)$. The symbol
$\LL\HH(X;Y)$ stands for the set of linear homeomorphisms between $X$
and $Y$.  If $X, Y$ are Banach spaces which form an interpolation
couple, then we denote
by $(X,Y)_{\tau,r}$ the real interpolation space, see~\cite{T78}.  We
use $M_3$ for the set of real, symmetric $3\times 3$-matrices. In the
sequel, a linear, continuous injection from $X$ to $Y$ is called an
\emph{embedding}, abbreviated by $X \embeds Y$. For
Lipschitz continuous functions $f$, we denote the Lipschitz constants
by $L_f$, while for bounded functions $g$ we denote their bound by
$M_g$ (both over appropriate sets, if necessary). Finally, $c$ denotes
a generic positive constant.

\subsection{Geometric setting for $\Omega$ and $\GD$}
\label{sec:geometric}

In all what follows, the symbol $\O$ stands for a bounded Lipschitz
domain in $\R^3$ in the sense of~\cite[Ch.\ 1.1.9]{mazya};
cf.~\cite{hal/rehcoer} for the boundary measure $\omega$ on such a domain.

\begin{rem} \label{r-Lipschitz} 
The thus defined notion is different
  from \emph{strong Lipschitz domain}, which is more restrictive and
  in fact identical with \emph{uniform cone domain}, see
  again~\cite[Ch.\ 1.1.9]{mazya}).
\end{rem}

A Lipschitz domain is formed e.g.\ by the topologically regularized
union of two crossing beams (see~\cite[Ch.~7]{hal/reh}), which is
\emph{not} a strong Lipschitz domain. Moreover, the interior of any
three-dimensional connected polyhedron is a Lipschitz domain, if the
polyhedron is, simultaneously, a $3$-manifold with boundary,
cf.~\cite[Thm.~3.10]{hhkrz}. However, a ball minus half of the
equatorial plate is \emph{not} a Lipschitz domain, and a chisel, where
the blade edge is bent onto the disc, is also not.

\begin{rem} \label{r-konsist} The Lipschitz property of $\O$ implies
  the existence of a linear, continuous extension operator
  $\mathfrak{E}:W^{1,q}(\O) \to W^{1,q}(\R^3)$
  (see~\cite[p.165]{gil}).  This has the following consequences:
  \begin{itemize}
  \item Since any element from $W^{1,q}(\R^3)$ may be approximated by
    smooth functions in the $W^{1,q}$-norm, any element from
    $W^{1,q}(\O)$ may be approximated by restrictions of smooth
    functions in the $W^{1,q}(\O)$-norm. This tells us that the
    definitions of $W^{1,q}(\O)$ and $W^{1,q}_\Xi(\O)$ are consistent
    in case of $\Xi=\emptyset$, i.e., one has $W^{1,q}(\O)=W_\emptyset
    ^{1,q}(\O)$. See also the detailed discussion
    in~\cite[Ch.~1.3.2]{G85}.
  \item It is not hard to see that $\mathfrak{E}$ also provides a
    continuous extension operator $\mathfrak{E}:C^\delta(\O) \to
    C^\delta(\R^3)$ and $\mathfrak{E}:L^p(\O) \to L^p(\R^3)$, where
    $\delta \in (0,1), p \in [1,\infty]$.
  \item Finally, the existence of the extension operator
    $\mathfrak{E}$ provides the usual Sobolev embeddings $W^{1,q}(\O)
    \hookrightarrow L^p(\O)$. In particular, this yields, by duality,
    the embedding $L^{q/2}(\O) \hookrightarrow W^{-1,q}_\emptyset(\O)$
    if $q$ exceeds the space dimension three.
  \end{itemize}
\end{rem}

Next we define the geometric setting for the domains $\O$ and the
Dirichlet boundary part. For this, we denote by $K$ the open unit cube
in $\R^n$, centered at $0 \in \R^n$, by $K_-$ the lower half cube $K
\cap \{\mathrm x \colon x_n < 0 \}$, by $\Si_K = {K} \cap \{\mathrm x \colon x_n =
0 \}$ the upper plate of $K_-$ and by $\Si_K^0$ the left half of
$\Si$, i.e. $\Si_K^0 = \Si_K \cap \{\mathrm x \colon x_{n-1} \le 0 \}$.

\begin{definition}\label{def:groeger}
 Let $\Xi\subset \partial\Omega$ be closed within $\partial\Omega$. 
  \begin{romannum}
  \item \romlabel{a-groeger-groeger} 
    We say that $\Omega \cup \Xi$ is \emph{regular} (in the sense of Gr\"oger), 
    if for any point $\mathrm x
    \in \partial \O$ there is an open neighborhood $U_{\mathrm{x}}$ of
    $\mathrm x$, a number $a_\mathrm{x} > 0$ and a bi-Lipschitz mapping $\phi_{\mathrm{x}}$ from
    $U_\mathrm x$ onto $a_\mathrm{x} K$ such that $\phi_\mathrm x(\mathrm
    x)=0 \in \R^3$, and we have either $\phi_{\mathrm{x}} \bigl( (\O
    \cup \Xi) \cap U_{\mathrm{x}} \bigr) = a_\mathrm{x}K_-$ or $a_\mathrm{x}(K_-
    \cup \Si_K)$ or $a_\mathrm{x}(K_- \cup \Si_K^0)$.
  \item \romlabel{a-groeger-volume} 
    The regular set $\Omega \cup \Xi$ is said to satisfy the \emph{volume-conservation condition}, if
    each mapping $\phi_{\mathrm{x}}$ in
    Condition~\ref{a-groeger-groeger} is volume-preserving.
  \end{romannum}  
\end{definition}

Generally, $\Xi$ is allowed to be empty in
Definition~\ref{def:groeger}. Then
Definition~\ref{def:groeger}~\ref{a-groeger-groeger} merely describes
a Lipschitz domain. Some further comments are in order:
\begin{rem} \label{r-groegerreg} \begin{romannum}
  \item Condition~\ref{a-groeger-groeger} exactly characterizes
    Gr\"oger's \emph{regular sets}, introduced in his pioneering
    paper~\cite{groeger89}. Note that the volume-conservation
    condition also has been required in several contexts,
    cf.~\cite{ggkr} and~\cite{groeger92}.

    Clearly, the properties $\phi_{\mathrm{x}} (U_\mathrm
    x)=a_\mathrm{x} K$ and $\phi_{\mathrm{x}} \bigl( \O \cap
    U_{\mathrm{x}} \bigr) = a_\mathrm{x}K_-$ are already ensured by
    the Lipschitz property of $\O$; the crucial point is the behavior
    of $\phi_\mathrm{x}(\Xi \cap U_\mathrm x)$.
  \item A simplifying topological characterization of Gr\"oger's
    regular sets in the case of three space dimensions reads as
    follows (cf.~\cite[Ch.~5]{HDMR08}):
    \begin{remunerate}
    \item $\Xi$ is the closure of its interior within $\pO$,
    \item the boundary $\partial \Xi$ within $\pO$ is locally
      bi-Lipschitz diffeomorphic to the open unit interval $(0,1)$.
    \end{remunerate}
  \item In particular, all domains with Lipschitz boundary
    (synonymous: \emph{strong} Lipschitz domains) satisfy
    Definition~\ref{def:groeger}: if, after a shift and an orthogonal
    transformation, the domain lies locally beyond a graph of a
    Lipschitz function $\psi$, then one can define $\phi(x_1, \ldots,
    x_d) = (x_1 - \psi(x_2, \ldots, x_d) ,x_2, \ldots,
    x_d)$. Obviously, the mapping $\phi$ is then bi-Lipschitz and the
    determinant of its Jacobian is identically $1$.
  \item It turns out that regularity together
    with the volume-conservation condition is not a too
    restrictive assumption on the mapping $\phi_x$. In particular, there are such
    mappings---although not easy to construct---which map the ball
    onto the cylinder, the ball onto the cube and the ball onto the
    half ball, see~\cite{ghkr, fonse}.  The general message is that
    this class has enough flexibility to map ``non-smooth'' objects
    onto smooth ones.
  \item \romlabel{r-groegerreg-interpol} The spaces $W^{1,q}_{\Xi}$
    and $W^{-1,q}_{\Xi}$ still exhibit the usual interpolation
    properties, see~\cite{ggkr} for details.
  \item If $\Xi$ is nonempty and $\Omega \cup \Xi$ is regular, then $\Xi$ has interior points (with respect to the boundary topology
    in $\partial \O$), and, consequently, never has boundary measure
    $0$.
  \end{romannum}
\end{rem}

The following assumption is supposed to be valid for all the remaining
considerations in the paper.

\begin{assu} \label{a-groeger} 
 The set $\Omega \cup \Gamma_D$ is regular with $\Gamma_D  \neq \emptyset$.
\end{assu}

 For the moment, it is sufficient to impose only the regularity condition from 
 Assumption~\ref{a-groeger}~\ref{a-groeger-groeger} on $\Omega\cup\Gamma_D$. 
 The volume-conservation condition is not needed until Section~\ref{s-oc}, 
 cf.\ Assumption~\ref{a-sigmanemytskii} below. As explained in
 Remark~\ref{r-groegerreg}, Assumption~\ref{a-groeger} in particular implies that
 $\omega(\Gamma_D) > 0$.

\subsection{General assumptions on~\eqref{P}}
\label{sec:generalassu}

Now we are in the position to state the main assumptions for the
quantities in~\eqref{P}. Please note that in order to obtain sharp
results we just give the assumptions on the quantities in
\eqref{e-strongparabolic}--\eqref{e-strongellipticdirichlet} which are
needed to obtain existence, uniqueness, and continuity of solutions to
the state system. For further considerations in \S\ref{s-oc}, in
particular those which include Fr\'echet-differentiability of the
associated solution operator, one has to require more restrictive
conditions on the nonlinearities, which are formulated in
Assumption~\ref{a-sigmanemytskii}, see \S\ref{s-oc}.

We first address the assumptions regarding (local) existence and uniqueness for the 
state equation~\eqref{e-strongparabolic}--\eqref{e-strongellipticdirichlet}. 
This means in particular that we treat $u$ as a fixed, given inhomogeneity in this context, 
whereas it is an unknown control function when considering the optimal control problem~\eqref{P}.

\begin{assu} \label{a-assu1} On the quantities in the state
  system~\eqref{e-strongparabolic}--\eqref{e-strongellipticdirichlet} we generally impose:
  \begin{romannum}
  \item \romlabel{a-assu1-sigmaeta} The functions $\sigma:\R \to
    (0,\infty)$ and $\eta \colon \R \to (0,\infty)$ are bounded and
    Lip\-schitzian on any bounded interval,
  \item \romlabel{a-assu1-eps} the function $\varepsilon \in
    L^\infty(\O;M_3)$ takes \emph{symmetric} matrices as values, and
    satisfies the usual ellipticity condition, i.e., \begin{equation*}
      \essinf_{x\in\O} \sum_{i,j=1}^3
      \varepsilon_{ij}(x)_{ij}\,\xi_i\,\xi_j\geq \underline
      \varepsilon\,\|\xi\|^2_{\R^3} \quad\forall\,\xi\in\R^3
    \end{equation*}
    with a constant $\underline \varepsilon >0$,
  \item \romlabel{a-assu1-kappa} the function $\kappa \in
    L^\infty(\O;M_3)$ also takes \emph{symmetric} matrices as values,
    and, additionally, satisfies an ellipticity condition, that is,
    \begin{equation*} \essinf_{x\in\O} \sum_{i,j=1}^3
      \kappa_{ij}(x)\,\xi_i\,\xi_j\geq \underline
      \kappa\,\|\xi\|^2_{\R^3} \quad\forall\,\xi\in\R^3
    \end{equation*}
    holds with a constant $\underline{\kappa} > 0$,
  \item \romlabel{a-assu1-thetal} $\theta_l \in
    L^\infty(J;L^\infty(\pO))$,
  \item \romlabel{a-assu1-alpha} $\alpha \in L^\infty(\pO)$ with
    $\alpha(x) \ge 0$ a.e.~on $\pO$\ and $\int_{\pO} \alpha \,d\omega >0$,
  \item \romlabel{a-assu1-u} $u \in L^{2r}(J; W^{-1,q}_{\Gamma_D})$ for some $q>3$ to be specified 
  in Assumption~\ref{a-reg} below and $r > \frac{2q}{q-3}$, cf.\ Definition~\ref{d-criticalexp} and 
  Theorem~\ref{t-locexist} below.
  \end{romannum}
\end{assu}

\begin{rem} \label{r-einbett/ident} In assumption~\ref{a-assu1-u},
  we implicitly made use of the embedding $L^{\mathfrak p}(\Gamma_N) \embeds
  W_{\Gamma_D}^{-1,q}$ for $\mathfrak p > \frac{2}{3}q$, realized by the adjoint
  operator of the continuous trace operator $\tau_{\Gamma_N} \colon W^{1,q'}_{\Gamma_D} \to
  L^{\mathfrak p'}(\Gamma_N)$. In this sense, a function $u\in L^{2r}(J; L^\mathfrak{p}(\Gamma_N))$
  is considered as an element of $L^{2r}(J;W_{\Gamma_D}^{-1,q})$. In the same manner, we will treat 
  the function $\alpha \theta_l \in L^\infty(J;L^\infty(\pO))$ as an element of
  $L^\infty(J;W_\emptyset^{-1,q})$, see~\cite[Lemma~2.7]{hoemkrum} for the required embeddings/trace
  operators. 
\end{rem}

Next we turn to the assumptions concerning the optimal control problem~\eqref{P}. 
Now, $u$ plays the role of the searched-for variable or function, whose regularity 
is implicitly determined by the objective functional in~\eqref{P}. 
As we will see in the sequel of \S~\ref{s-oc}, our hypotheses on the objective functional 
stated below imply that it suffices to restrict to control functions in a function space $\U$, 
see~\eqref{eq:ctrlspace}, which continuously embeds in $L^{2r}(J;
W^{-1,q}_{\Gamma_D})$ as required in
Assumption~\ref{a-assu1}~\ref{a-assu1-u}, see
Proposition~\ref{p-compactcontrolspace} below. 

\begin{assu}\label{assu-p}
  The remaining quantities in~\eqref{P} fulfill:
  \begin{romannum}
  \item The integrability exponents in the objective functional satisfy 
  $p > \frac{4}{3}q-2$ and $s > \frac{2q}{q-3}(1 - \frac{3}{q} + \frac{3}{\varsigma})$, 
  where $q$ and $\varsigma$ are specified in Assumption~\ref{a-reg} and Definition~\ref{def:barq} below.
  \item $E$ is an open (not necessarily proper) subset of $\O$.
  \item $\theta_d \in L^2(E)$.
  \item \romlabel{assu-p-thetamax} $\theta_{\max}\in C(\overline Q)$
    with $\max(\max_{\olom} \theta_0,\esssup_{\Sigma}\theta_l) \leq
    \theta_{\max}(x,t)$ for all $(x,t) \in \overline Q$ and
    $\theta_0(x) < \theta_{\max}(T_0,x)$ for all $x \in \olom$.
    \item $u_{\max}$ is a given function with $u_{\max}(x,t)\geq 0$ a.e.~on $\Si_N$.
   \item $\beta>0$.
  \end{romannum}
\end{assu}

Note that we do not impose any regularity assumptions on the function $u_{\max}$. 
In particular, it is allowed that $u_{\max} \equiv \infty$ so that no upper bound is present.


\section{Rigorous formulation, existence and uniqueness of solutions for the thermistor problem}
\label{s-statesys}

In this chapter we will present a precise analytical formulation for
the thermistor-problem, see Definition~\ref{d-loes} below. In order to do so, we first recall some
background material. One of the most crucial points is the requirement of
suitable mapping property for Poisson's operator, cf.\ Assumption~\ref{a-reg}. 
The reader should note that a similar
condition was also posed in~\cite[Ch.~3]{AC94} in order to get
smoothness of the solution; compare also~\cite{gajewski}, where
exactly this regularity for the solution of Poisson's equation is
needed in order to show uniqueness for the semiconductor equations. We
prove, in particular, some preliminary results which are needed later
on and which may be also of independent interest.  After having
properly defined a solution of the thermistor problem, we establish
some more preparatory results and afterwards show existence (locally in
time) and uniqueness of the solution of the thermistor problem in Section~\ref{subs-proof}.
Finally, we show that our concept to treat the problem is not
accidental, but---more or less---inevitable.

\subsection{Prerequisites: Elliptic and parabolic regularity}\label{subs-prereq}
We begin this subsection with the definition of the
divergence operators. First of all, let us introduce the brackets
$\langle \cdot, \cdot \rangle$ as the symbol for the dual pairing
between $W^{-1,2}_{\Xi}$ and $W^{1,2}_{\Xi}$, extending the scalar product in $L^2$.

\begin{definition} \label{d-opera0} Let $\Xi \subset \pO$ be
  closed. Assume that $\mu $ is any bounded, measurable, $M_3$-valued
  function on $\O$ and that $\gamma \in L^\infty(\pO \setminus \Xi)$
  is nonnegative.  We define the operators $-\nabla \cdot \mu \nabla$
  and $-\nabla \cdot \mu + \tilde{\gamma}$, each mapping
  $W^{1,2}_{\Xi}$ into $W^{-1,2}_{\Xi}$, by
  \begin{equation} \label{e-ellipoper0} \left\langle -\nabla \cdot \mu
      \nabla \psi, \xi \right\rangle := \int_{\O} \mu \nabla \psi
    \cdot \nabla \xi \, \dd x \quad \text{for} \quad \psi, \xi \in
    W^{1,2}_{\Xi}
  \end{equation}
  and
  \begin{equation} \label{e-ellipoper}\left\langle(-\nabla \cdot \mu
      \nabla + \tilde{\gamma})\psi,\xi\right\rangle =
    \left\langle-\nabla \cdot \mu \nabla\psi,\xi\right\rangle
    +\int_{\pO \setminus \Xi} \gamma \, \psi \, \xi \, \dd\omega \quad
    \text{for} \quad \psi, \xi \in W^{1,2}_{\Xi}.
  \end{equation}
  In all what follows, we maintain the same notation for the
  corresponding maximal restrictions to $W^{-1,q}_\Xi$, where $q >2$.

\end{definition}
\begin{rem} \label{r-einschr} Let us denote the domain for the
  operator $-\nabla \cdot \mu \nabla$, when restricted to
  $W^{-1,q}_\Xi$ ($q >2$), by $\mathcal D_{q}$, equipped with the
  graph norm.  Then the estimate
  \begin{equation} \label{e-contrac} \|-\nabla \cdot \mu \nabla
    \psi\|_{W^{-1,q}_\Xi}=\sup_{\|\varphi\|_{W^{1,q'}_\Xi}=1} \left
      |\int_\O \mu \nabla \psi \cdot \nabla \varphi \,d \mathrm x
    \right| \le \|\mu\|_{L^\infty } \|\psi\|_{W^{1,q}_\Xi}
  \end{equation}
  shows that $W^{1,q}_\Xi$ is embedded in $\mathcal D_{q}$ for every
  bounded coefficient function $\mu$. It is also known that $\DD_q
  \embeds C^\alpha(\olom)$ for some $\alpha > 0$ whenever $q > 3$, see~\cite[Thm.~3.3]{HDMR08}. Additionally,~\eqref{e-contrac}
  implies that the mapping
  \[
  L^\infty(\O; M_3) \ni \mu \mapsto \nabla \cdot \mu \nabla \in \LL
  (W^{1,q}_{\Xi};W^{-1,q}_{\Xi})
  \]
  is a linear and continuous contraction for every $q \in (1,\infty)$.
\end{rem}

In the following, we consider the operators defined in
Definition~\ref{d-opera0} mostly in two incarnations: firstly, the
case $\Xi = \emptyset$ and $\mu = \kappa$; and secondly $\Xi = \GD$
with $\mu = \varepsilon$. We write $-\nabla \cdot \kappa \nabla$ and
$-\nabla \cdot \kappa\nabla + \tilde{\alpha}$ in the first, and
$-\nabla \cdot \varepsilon \nabla$ in the second case. We recall
various properties of operators of the form $-\nabla \cdot \mu\nabla$.

\begin{proposition} \label{p-divproperty} Let $\O\cup\Xi$ be regular in the sense 
of Definition~\ref{def:groeger} and suppose
  that the coefficient function $\mu$ in~\eqref{e-ellipoper} is real,
  bounded and elliptic.
  \begin{romannum}
  \item Suppose that either $\omega (\Xi)>0$ or $\Xi = \emptyset$ and
    $\int_{\pO} \gamma \, d \omega >0$.
    \begin{remunerate}
    \item \cite{hal/rehcoer} The quadratic form corresponding
      to~\eqref{e-ellipoper} is coercive.
    \item \romlabel{p-divproperty-groeger}\cite{groeger89} There is a
      number $q_0 >2$ such that
      \[
      -\nabla \cdot \mu \nabla + \tilde{\gamma} \colon W^{1,q}_\Xi \to
      W^{-1,q}_\Xi
      \]
      is a topological isomorphism for all $q \in [2,q_0]$. The number
      $q_0$ may be chosen uniformly for all coefficient functions
      $\mu$ with the same ellipticity constant and the same
      $L^\infty$-bound. Moreover, for each  $q \in [2,q_0]$, the norm of the inverse of 
      $\nabla \cdot \mu \nabla + \tilde{\gamma}$ as a mapping from
      $W^{-1,q}_\Xi$ to $W^{1,q}_\Xi$ may be estimated again
      uniformly for all 
      coefficient functions with the same ellipticity constant and the same
      $L^\infty$-bound.
    \end{remunerate}
  \item \romlabel{p-divproperty-sqrroot} Assume that $\gamma$ is a
    nonnegative function from $L^\infty(\partial \O \setminus \Xi)$
    and that the coefficient function $\mu$ takes symmetric matrices
    as values.
    \begin{remunerate}
    \item \cite[Cor.~5.21]{hal/rehPara} The operator $-\nabla \cdot
      \mu \nabla +\tilde \gamma +1$ is a \emph{positive} one on any
      space $W^{-1,q}_\Xi$, if $q \in [2,6]$, i.e., one has the
      resolvent estimate
      \[
      \sup_{\lambda \in [0,\infty )}(\lambda+1) \|(-\nabla \cdot \mu \nabla
      +\tilde \gamma +1+\lambda)^{-1} \|_{\mathcal L (W^{-1,q}_\Xi)} <
      \infty.
      \]
      In particular, all fractional powers of $-\nabla \cdot \mu
      \nabla +\tilde \gamma +1$ are well-defined and possess the usual
      properties, cf.~\cite[Ch.~1.14]{T78}.
    \item \cite[Thm.~4.2]{hal/rehPara} The square root satisfies
      $(-\nabla \cdot \mu \nabla +\tilde \gamma +1)^{-1/2} \in
      \mathcal L(W^{-1,q}_\Xi;L^q)$, or in other words, $\dom\bigl
      ((-\nabla \cdot \mu \nabla +\tilde \gamma +1)^{1/2}\bigr )$
      embeds into $L^q$, if $q \in [2,\infty)$.
    \end{remunerate}
  \end{romannum}
\end{proposition}


See also~\cite{ABHDR} for recent results as in
Proposition~\ref{p-divproperty}~\ref{p-divproperty-sqrroot} in a
broader context. Our next aim is to introduce the solution concept for the thermistor
problem. To this end, we make the following assumption (cf.\ also
Remark~\ref{r-explain} below):

\begin{assu} \label{a-reg} There is a $q \in (3,4)$ such that the
  mappings
  \begin{equation} \label{e-topisoq} -\nabla \cdot \varepsilon \nabla:
    W^{1,q}_{\GD} \to W^{-1,q}_{\GD}
  \end{equation}
  and
  \begin{equation} \label{e-topisoq2} -\nabla \cdot \kappa \nabla + 1:
    W^{1,q} \to W^{-1,q}_{\emptyset}
  \end{equation}
  each provide a topological isomorphism.
\end{assu}

The papers~\cite[Appendix]{krum/reh} and~\cite{DKR} provide a zoo of
arrangements such that Assumption~\ref{a-reg} is satisfied. Note that
it is not presumptous to assume that both differential operators
provide topological isomorphisms at the same time, since the latter
property mainly depends on the behaviour of the discontinuous
coefficient functions (versus the geometry of $\GD$), and these
correspond to the material properties in the workpiece described by
the domain $\O$, i.e., the coefficient functions should exhibit
similar properties with regard to jumps or discontinuities in general,
the main obstacles to overcome for the isomorphism property. Since
$\kappa$ is not assumed to be continuous,
Assumption~\eqref{e-topisoq2} is not satisfied \emph{a priori}, even
though no mixed boundary conditions are present, see
\cite[Ch.~4]{e/r/s} for a striking example. In this sense, mixed
boundary conditions are not a stronger obstruction against higher
regularity in the range
$q \in (3,4)$ than discontinuous coefficient functions are.

\begin{rem} \label{r-erlaeutern} In case of mixed boundary conditions
  it does not make sense to demand Assumption~\ref{a-reg}---even if
  all data are smooth---for a $q \ge 4$, due to Shamir's famous
  counterexample~\cite{shamir}. Note further that the isomorphism
  properties in~\eqref{e-topisoq} and~\eqref{e-topisoq2} remain valid for
  all other $\tilde q \in [2,q)$ due to interpolation, cf.\
  Remark~\ref{r-groegerreg}~\ref{r-groegerreg-interpol}.
\end{rem}

In order to treat the quasilinearity in~\eqref{e-strongparabolic}, we
need to ensure a certain uniformity of domains of the differential
operator $-\nabla \cdot \eta(\theta)\kappa\nabla$ during the
evolution. To this end, we first note that the isomorphism-property
for $-\nabla\cdot\kappa\nabla+1$ from Assumption~\ref{a-reg} extends
to a broader class of coefficient functions.

\begin{definition} \label{d-stetiglowbound} Let $\underline C(\olom)$
  denote the set of positive functions on $\O$ which are uniformly
  continuous and admit a positive lower bound.
\end{definition}

\begin{lemma} \label{l-multiplierkappa} Assume that
  Assumption~\ref{a-reg} holds for some number $q \in [2,4)$. If $\xi
  \in \underline C(\olom)$, then~\eqref{e-topisoq} and
  ~\eqref{e-topisoq2} remain topological isomorphisms, if
  $\varepsilon$ and $\kappa$ are replaced by $\xi \varepsilon$ and
  $\xi \kappa$, respectively.
\end{lemma}

A proof can be found in~\cite[Ch.~6]{DKR}.

\begin{corollary}\label{c-quasidomains} Assume that~\eqref{e-topisoq2}
  is a topological isomorphism for some $q \in [2,4)$.  Then, for
  every $\xi \in \underline C(\olom)$, the domain of the operator
  $-\nabla \cdot \xi\kappa\nabla + \tilde{\alpha}$, considered in
  $W_\emptyset^{-1,q}$, is still $W^{1,q}$. In particular, for every
  function $\zeta \in C(\olom)$, the operator $-\nabla \cdot
  \eta(\zeta)\kappa\nabla + \tilde{\alpha}$ has domain $W^{1,q}$.
\end{corollary}

\begin{proof}
  The first assertion follows from Lemma~\ref{l-multiplierkappa} and
  relative compactness of the boundary integral in $\tilde{\alpha}$
  with respect to $-\nabla \cdot \xi \kappa \nabla$, compare
  ~\cite[Ch. IV.1.3]{kato}. For the second assertion, note that $\eta$
  is assumed to be Lipschitzian on bounded intervals and bounded from
  below by $0$ as in Assumption~\ref{a-assu1}. Thus, $\eta(\zeta)$ is
  uniformly continuous and has a strictly positive lower bound.
  \hfill\end{proof}

We are now in the position to define what is to be understood as a
\emph{solution} to the
system~\eqref{e-strongparabolic}--\eqref{e-strongellipticdirichlet}.

\begin{definition}\label{d-quasilinear} We define
  \[\AA(\zeta) := -\nabla \cdot \eta(\zeta)\kappa\nabla +
  \tilde{\alpha}\] as a mapping $\AA \colon C(\olom) \to
  \LL(W^{1,q};W^{-1,q}_\emptyset)$.
\end{definition}

\begin{definition}
  \label{d-criticalexp} The number $r^*(q) = \frac{2q}{q-3}$
  is called the \emph{critical exponent}.
\end{definition}

\begin{definition} \label{d-loes} Let $q > 3$ and let $r$ be from
  $(r^*(q),\infty)$.  For given $J = (T_0,T_1)$, we call the
  pair $(\theta,\varphi)$ a \emph{solution} of the thermistor-problem,
  if it satisfies the equations
  \begin{align}
    \theta'(t) + \AA(\theta(t))\theta(t) & =
    \left(\sigma(\theta(t))\varepsilon \nabla \varphi(t) \right )
    \cdot \nabla \varphi(t) + \alpha\theta_l(t)
    && \text{in } W^{-1,q}_\emptyset, \label{e-paraex}\\
    {-\nabla \cdot \sigma(\theta(t) ) \varepsilon \nabla} \varphi(t) &
    = u(t) && \text{in } W^{-1,q}_{\GD} \label{e-elliptex}
  \end{align}
  with $\theta(T_0)=\theta_0$ for almost all $t \in (T_0,T_1)$, where
  \begin{equation}
    \varphi \in L^{2r} (J;W^{1,q}_{\GD}) \quad \text{and} \quad \theta \in W^{1,r}(J;W^{-1,q}_\emptyset) \cap
    L^{r}(J;W^{1,q}).
    \label{e-solreg}
  \end{equation}
  We call $(\theta,\varphi)$ a \emph{local solution}, if it
  satisfies~\eqref{e-paraex} and~\eqref{e-elliptex} in the above
  sense, but only on $(T_0,T^\bullet) \subseteq (T_0,T_1)$.
\end{definition}

\begin{rem} \label{r-Loesung}
  \begin{romannum}
  \item In the context of Definition~\ref{d-loes}, $\theta'$ always means
    the time derivative of $\theta$ in the sense of vector-valued
    distributions, see~\cite[Ch.~III.1]{A95} or~\cite[Ch.~IV]{ggz}.
  \item Via~\eqref{e-maxreghoelderembed} and
    Corollary~\ref{c-intpolcompact} below, we will see that a solution
    $\theta$ in the above sense is in fact H\"older-continuous on
    $\overline{\O \times J}$. In particular, $\theta(t)$ is uniformly
    continuous on $\O$ for every $t \in J$, such that $\AA(\theta(t))$
    is well-defined according to Definition~\ref{d-quasilinear}.
  \item \romlabel{r-boundcond} The reader will verify that the
    boundary conditions imposed on $\varphi$
    in~\eqref{e-strongellipticneumann}
    and~\eqref{e-strongellipticdirichlet} are incorporated in this
    definition in the spirit of~\cite[Ch.~II.2]{ggz}
    or~\cite[Ch.~1.2]{cia}.  For an adequate interpretation of the
    boundary conditions for $\theta$ as
    in~\eqref{e-strongparabolicrobin}, see~\cite[Ch.~3.3.2]{lions} and
    the in-book references there.
  \end{romannum}
\end{rem}

We are now going to formulate the main result of this part.

\begin{theorem} \label{t-locexist} Let $q \in (3,4)$ be a number for
  which Assumption~\ref{a-reg} is satisfied, $r > r^*(q)$ and
  $u \in L^{2r}(J;W^{-1,q}_\GD)$, where $r^*(q)$ 
  is the critical exponent from Definition~\ref{d-criticalexp}.
  If $\theta_0$ is from $(W^{1,q},W^{-1,q}_\emptyset)_{\frac{1}{r},r}$, then there is a
  unique \emph{local} solution of~\eqref{e-paraex}
  and~\eqref{e-elliptex} in the sense of Definition~\ref{d-loes}.
\end{theorem}

The proof of this theorem is given in the next subsection.

\subsection{Local existence and uniqueness for the state system: the proof} \label{subs-proof}

Let us first briefly sketch the proof of Theorem~\ref{t-locexist} by giving an overview 
over the steps:
\begin{itemize}
 \item The overall proof is based on a local existence result of Pr\"uss 
 for abstract quasilinear parabolic equations, whose principal part satisfies a certain 
 maximal parabolic regularity property, see~\cite{pruess} and Proposition~\ref{p-pruess}.
 \item For the application of this abstract result to our problem, we 
 reduce the thermistor system to an equation in the temperature
 $\theta$ only by solving the elliptic equation for $\varphi$ in
 dependence of $\theta$.
 This gives rise to a nonlinear operator $S$ appearing in the reduced equation for $\theta$, see Definition~\ref{d-ellright} and 
 Proposition~\ref{p-assuS}. 
 \item The key tool to verify the assumptions on $S$ for the application of Pr\"uss' result 
 is Lemma~\ref{l-multiplierkappa}, which is the basis for the proof of 
 Lemma~\ref{t-righthands}. The application of Lemma~\ref{l-multiplierkappa} requires 
 to treat the temperature in a space which (compactly) embeds into $C(\overline{\Omega})$. 
 This issue is addressed by Corollary~\ref{c-intpolcompact}.
\end{itemize}


Before we start with the proof itself, let us first recall the concept of maximal parabolic regularity, a
crucial tool in the following considerations, and point out some basic
facts on this:

\begin{definition} \label{d-maxpar} Let $X$ be a Banach space and $A$
  be a closed operator with dense domain $\dom(A) \subset X$.  Suppose
  $\mfr \in (1,\infty)$. Then we say that $A$ has \emph{maximal
    parabolic $L^\mfr(J;X)$-regularity}, iff for every $f \in L^\mfr(J;X)$
  there is a unique function $w \in W^{1,\mfr}(J;X) \cap L^{\mfr}(J;\dom(A))$
  which satisfies \begin{equation}\label{e-para} w'(t) + Aw(t)=f(t),
    \quad w(T_0)=0\end{equation} in $X$ for almost every $t \in J=
  (T_0,T_1)$.
\end{definition}

\begin{rem} \label{r-adequat}
  \begin{romannum}
  \item As in Remark~\ref{r-Loesung}, $w'$ in Definition~\ref{d-maxpar} also always means
    the time derivative of $w$ in the sense of vector-valued
    distributions.
  \item We consider the concept of maximal parabolic regularity as
    adequate for the solution since it allows for discontinuous (in
    time) right hand sides---as are required in our context and in
    many other applications.
  \end{romannum}
\end{rem}

\begin{rem} \label{r-maxparregfacts} The following results on maximal
  parabolic $L^\mfr(J;X)$-regularity are well-known:
  \begin{romannum}
  \item \romlabel{r-maxparregfacts-int} If $A$ satisfies maximal
    parabolic $L^\mfr(J;X)$-regularity, then it does so for any other
    (bounded) time interval, see~\cite{dore}.
  \item \romlabel{r-maxparregfacts-exp} If $A$ satisfies maximal
    parabolic $L^\mfr(J;X)$-regularity for \emph{some} $\mfr \in
    (1,\infty)$, then it satisfies maximal parabolic
    $L^\mfr(J;X)$-regularity for \emph{all} $\mfr \in (1,\infty)$,
    see~\cite{sob} or~\cite{dore}.
  \item Let $Y$ be another Banach space, being dense in $X$ with $Y
    \embeds X$. Then there is an embedding
    \begin{equation} \label{e-maxreghoelderembed} W^{1,\mfr}(J;X) \cap
      L^{\mfr}(J;Y) \embeds C^\rho(\olJ;(Y,X)_{\zeta,1})
    \end{equation}
    where $0 < \rho \le \zeta - \frac {1}{\mfr}$,
    see~\cite[Ch.~3,~Thm.~3]{A01}. In the immediate context of maximal
    parabolic regularity, $Y$ is taken as $\dom(A)$ equipped with the
    graph norm, of course.
  \end{romannum}
  According to~\ref{r-maxparregfacts-int}
  and~\ref{r-maxparregfacts-exp}, we only say that $A$ satisfies
  \emph{maximal parabolic regularity on $X$}.
\end{rem}

In the following, we establish some preliminary results for the proof
of Theorem~\ref{t-locexist}, which will heavily rest on the following
fundamental theorem of Pr\"uss, see~\cite{pruess}:

\begin{proposition} \label{p-pruess} Let $Y,X$ be Banach spaces, $Y$
  dense in $X$, such that $Y \embeds X$ and set $J = (T_0,T_1)$ and $\mfr
  \in (1,\infty)$. Suppose that $A$ maps $(Y,X)_{\frac{1}{\mfr},\mfr}$ into
  $\LL(Y;X)$ such that $A(w_0)$ satisfies maximal parabolic regularity
  on $X$ with $\dom (A(w_0)) = Y$ for some $w_0 \in
  (Y,X)_{\frac{1}{\mfr},\mfr}$.  Let, in addition, $S \colon J \times
  (Y,X)_{\frac {1}{\mfr},\mfr} \to X$ be a Carath\'eodory map and $S(\cdot,
  0)$ be from $L^\mfr(J;X)$. Moreover, let the following two assumptions
  be satisfied:
  \begin{remunerate}
  \item[\refA] \label{a-pruessA} For every $M>0$, there is a constant
    $L(M)$ such that for all $w,\bar{w} \in (Y,X)_{\frac{1}{\mfr},\mfr}$,
    where
    $\max(\|w\|_{(Y,X)_{\frac{1}{\mfr},\mfr}},\|\bar{w}\|_{(Y,X)_{\frac{1}{\mfr},\mfr}})
    \leq M$, we have
    \begin{equation} \label{e-lipA} \| A(w) - A(\bar{w})\|_{\LL(Y;X)}
      \le L(M) \| w - \bar{w} \|_{(Y,X)_{\frac{1}{\mfr},\mfr}}.
    \end{equation}
        %
  \item[\refS] \label{a-pruessS} For every $M>0$, assume that there is
    a function $h_M \in L^\mfr(J)$ such that for all $w,\bar{w} \in
    (Y,X)_{\frac{1}{\mfr},\mfr}$, where
    $\max(\|w\|_{(Y,X)_{\frac{1}{\mfr},\mfr}},\|\bar{w}\|_{(Y,X)_{\frac{1}{\mfr},\mfr}})
    \leq M$, it is true that
    \begin{equation} \label{e-lipf} \| S(t,w) - S(t,\bar{w})\|_X \le
      h_M(t) \| w - \bar{w} \|_{(Y,X)_{\frac{1}{\mfr},\mfr}}
    \end{equation}
    for almost every $t \in J$.
  \end{remunerate}

  Then, for each $w_0 \in (Y,X)_{\frac{1}{\mfr},\mfr}$, there exists $T_{\max}
  \in J$ such that the problem
  \begin{equation} \label{e-pruessloes} \left\{ \begin{aligned} w'(t)
        + A(w(t)) w(t) &= S(t,w(t)) \quad \text{in }
        X 
        , \\
        w(T_0) &= w_0
      \end{aligned} \right.
  \end{equation}
  admits a \emph{unique} solution $w \in W^{1,\mfr}(T_0,T_\bullet;X) \cap
  L^\mfr(T_0,T^*_\bullet;Y)$ on $(T_0,T_\bullet)$ for every $T_\bullet \in (T_0,T_{\max})$.
\end{proposition}

\begin{rem} \label{r-minimal}
  It is known that the solution of the thermistor problem possibly
  ceases to exist after finite time in general, cf.~\cite[Ch.~5]{AC94}
  and the references therein. Thus, one has to expect here, in
  contrast to the two-dimensional case treated in~\cite{h/m/r/r}, only
  a local-in-time solution. In this scope, Pr\"uss' theorem will prove
  to be the adequate instrument.
\end{rem}

As indicated above, we will prove Theorem~\ref{t-locexist} by reducing 
the thermistor system to an equation in the temperature only and apply Proposition~\ref{p-pruess}
to this equation. To be more precise, we first establish the assumptions~\refA\ for $\mfr =r > r^*(q)$ and $\AA$ as
defined in Definition~\ref{d-quasilinear}. We then
solve the elliptic equation~\eqref{e-elliptex} for $\varphi$
(uniquely) for every time point $t$ in dependence of a function
$\zeta$ and $u(t)$, where $\zeta$ enters the equation inside the
coefficient function $\sigma(\zeta)\varepsilon$. Then the right-hand
side of the parabolic equation~\eqref{e-paraex} may be written also as
a function $S$ solely of $t$ and $\zeta$. We then show that this
function satisfies the suppositions~\refS\ in Pr\"uss' theorem.

To carry out this concept, we need several prerequisites: here our
first central aim is to show that indeed the mapping
$(W^{1,q},W^{-1,q}_\emptyset)_{\frac{1}{r},r} \ni \zeta \mapsto
\AA(\zeta)$ from Definition~\ref{d-quasilinear} satisfies the
assumptions from Proposition~\ref{p-pruess} for $r > r^*(q)$,
cf.\ Lemma~\ref{p-assuA} below. For doing so, we first investigate the
spaces $(W^{1,q},W^{-1,q}_\emptyset)_{\zeta,1}$ in view of their
embedding into H\"older spaces. For later use, the subsequent result
is formulated slightly broader as presently needed.

\begin{theorem} \label{p-maxparregembed} Let $q \in (3,4)$ and
  $\varsigma \in [2,q]$.  For every $\tau \in
  (0,\frac{q-3}{2q}(1-\frac{3}{q}+\frac{3}{\varsigma})^{-1})$, the
  interpolation space $(W^{1,q},W^{-1,\varsigma}_\emptyset)_{\tau,1}$
  embeds into some H\"older space $C^\delta(\olom)$ with $\delta > 0$.
\end{theorem}
\begin{proof}
 We apply the reiteration theorem~\cite[Ch.~1.10.2]{T78} to obtain
  \begin{multline}\label{e-reitera}
    (W^{1,q},W_\emptyset^{-1,\varsigma})_{\tau,1} =
    (W^{1,q},(W^{1,q},W_\emptyset^{-1,\varsigma})_{\frac
      {1}{2},1})_{2\tau,1} \\ \hookrightarrow
    (W^{1,q},(W^{1,\varsigma},W_\emptyset^{-1,\varsigma})_{\frac
      {1}{2},1})_{2\tau,1} \hookrightarrow
    (W^{1,q},(W_\emptyset^{-1,\varsigma},\mathcal
    D_\varsigma)_{\frac{1}{2},1})_{2\tau,1},
  \end{multline}
  where $\mathcal D_\varsigma$ denotes the domain of the Laplacian
  $-\Delta+1$ acting on the Banach space $W^{-1,\varsigma}_\emptyset$,
  cf.\ Remark~\ref{r-einschr}.  Denoting the domain of $(-\Delta
  +1)^{1/2}$, considered on the same space, by $\mathcal
  D_{\varsigma}^{\frac{1}{2}}$, one has
  $(W_\emptyset^{-1,\varsigma},\mathcal D_\varsigma)_{\frac{1}{2},1}
  \hookrightarrow \mathcal D_{\varsigma}^{\frac {1}{2}}$,
  cf.~\cite[Ch.~1.15.2]{T78}. Due to
  Proposition~\ref{p-divproperty}~\ref{p-divproperty-sqrroot}, we
  already know the embedding $\mathcal D_{\varsigma}^{\frac {1}{2}}
  \hookrightarrow L^\varsigma$.  Inserting in~\eqref{e-reitera}, this
  altogether yields $(W^{1,q},W_\emptyset^{-1,\varsigma})_{\tau,1}
  \hookrightarrow (W^{1,q},L^\varsigma)_{2\tau,1}$.

  We define $p:=\bigl (\frac {1-2\tau}{q}+\frac
  {2\tau}{\varsigma}\bigr )^{-1}$ and observe that $\delta:=1 - 2\tau
  - \frac{3}{p} \in (0,1)$, due to our condition on $\tau$.  Denoting
  by $H^{t,p}$ the corresponding space of Bessel potentials (cf.~\cite
  [Ch.~4.2.1]{T78}) one has the embedding $H^{1-2\tau,p}
  \hookrightarrow C^\delta(\olom) $, see~\cite[Thm.~4.6.1]{T78}.
  This, combined with the interpolation inequality for $H^{1-2\tau,p}$
  (\cite[Thm.~3.1]{ggkr}) gives for any $\psi \in W^{1,q}$ the
  estimate
  \begin{equation} \label{e-intineq} \|\psi\|_{C^\delta(\olom)} \le
    \|\psi\|_{H^{1-2\tau,p}} \le \| \psi\|_{W^{1,q}}^{1-2\tau} \|
    \psi\|^{2\tau}_{L^\varsigma}.
  \end{equation}
  But it is well-known (cf.~\cite[Ch.~1.10.1]{T78} or~\cite[Ch.~5,
  Prop.~2.10]{bennet}) that an inequality of type~\eqref{e-intineq} is
  constitutive for the embedding $(W^{1,q},L^\varsigma)_{2\tau,1}
  \embeds C^\delta(\olom)$. \hfill
\end{proof}

\begin{corollary} \label{c-intpolcompact}
  \begin{romannum}
  \item \romlabel{c-intpolcompact-real} Let $q > 3$ and $\varsigma \in
    [2,q]$.  Then, for every $s >
    \frac{2q}{q-3}(1-\frac{3}{q}+\frac{3}{\varsigma})$, the
    interpolation space $(W^{1,q},W^{-1,\varsigma}_\emptyset)_{\frac
      {1}{s},s}$ embeds into some H\"older space $C^\delta(\olom)$,
    and thus even compactly into $C(\olom)$.
  \item \romlabel{c-intpolcompact-mprspace} Under the same
    supposition, there exists a $\varrho > 0$ such that
    \[W^{1,s}(J;W^{-1,\varsigma}_\emptyset) \cap L^{s}(J;W^{1,q})
    \embeds C^\varrho(\olJ;C^\varrho(\olom)).\]
  \item \romlabel{c-intpolcompact-mpr} Let Assumption~\ref{a-reg} hold
    true for some $q \in (3,4)$. Then the operator $\AA(\zeta)$
    satisfies maximal parabolic regularity on $W^{-1,q}_\emptyset$
    with domain $W^{1,q}$ for every $\zeta \in
    (W^{1,q},W^{-1,q}_\emptyset)_{\frac{1}{r},r}$ with $r > r^*(q)$, where $r^*(q)$ 
    is the critical exponent from Definition~\ref{d-criticalexp}.
  \end{romannum}
\end{corollary}

\begin{proof}
  \ref{c-intpolcompact-real}~We have
  $(W^{1,q},W^{-1,\varsigma}_\emptyset)_{\frac{1}{s},s} \embeds
  (W^{1,q},W^{-1,\varsigma}_\emptyset)_{\iota,1}$ for every $\iota \in
  (\frac{1}{s},1)$. The condition on $s$ implies that the interval
  $\mathcal I:=(\frac{1}{s}
  ,\frac{q-3}{2q}(1-\frac{3}{q}+\frac{3}{\varsigma})^{-1})$ is
  non-empty. Taking $\iota$ from $\mathcal I$, the assertion follows
  from Theorem~\ref{p-maxparregembed}.
  \ref{c-intpolcompact-mprspace}~follows from
  Theorem~\ref{p-maxparregembed} and
  Remark~\ref{r-maxparregfacts}. \ref{c-intpolcompact-mpr}~The claim
  follows from uniform continuity of functions from
  $(W^{1,q},W^{-1,q}_\emptyset)_{\frac{1}{r},r}$
  by~\ref{c-intpolcompact-real}, Lemma~\ref{l-multiplierkappa} for
  $\xi := \eta(\zeta)$
  and~\cite[Thm.~5.4/5.19~(ii)]{hal/rehPara}. \hfill\end{proof}

Setting $\varsigma = q$ in
Corollary~\ref{c-intpolcompact}~\ref{c-intpolcompact-real}
and~\ref{c-intpolcompact-mprspace} gives the condition $r > r^*(q) =
\frac{2q}{q-3}$ for the assertions to hold with $s = r$. We will use this special
case frequently in the course of the remaining part of this
section. Let us now turn to the operator $\AA$.

\begin{proposition}
  \label{p-assuA} Suppose that Assumption~\ref{a-reg} holds true for
  some $q \in (3,4)$ and that $\theta_0 \in
  (W^{1,q},W^{-1,q}_\emptyset)_{\frac{1}{r},r}$ where $r > r^*(q)$. 
  With $\mathcal A$ as in Definition~\ref{d-quasilinear}, the
  function $(W^{1,q},W^{-1,q}_\emptyset)_{\frac{1}{r},r} \ni \zeta
  \mapsto \AA(\zeta)$ then satisfies the assumptions from
  Proposition~\ref{p-pruess} for the spaces $X = W^{-1,q}_\emptyset$
  and $Y = W^{1,q}$.
\end{proposition}

\begin{proof}
  With $\varsigma = q$, Corollary~\ref{c-intpolcompact} shows that
  $(W^{1,q},W^{-1,q}_\emptyset)_{\frac{1}{r},r} \embeds C(\olom)$,
  such that the operator $\AA$ indeed maps
  $(W^{1,q},W^{-1,q}_\emptyset)_{\frac{1}{r},r}$ into
  $\LL(W^{1,q};W^{-1,q}_\emptyset)$ by
  Corollary~\ref{c-quasidomains}. Using Lipschitz continuity of $\eta$
  on bounded sets and Remark~\ref{r-einschr}, we also obtain~\refA:
  Let $w,\bar{w} \in (W^{1,q},W^{-1,q}_\emptyset)_{\frac{1}{r},r}$
  with norms bounded by $M > 0$. Then we have
  \begin{align*}
    \left\|\AA(w) -
      \AA(\bar{w})\right\|_{\LL(W^{1,q};W^{-1,q}_\emptyset)} & =
    \left\|\nabla\cdot \left(\eta(w)-\eta(\bar{w})\right)\kappa\nabla
    \right\|_{\LL(W^{1,q},W^{-1,q}_\emptyset)} \\
    & \leq L_\eta \|\kappa\|_{L^\infty}\|w-\bar{w}\|_{C(\olom)} \\
    & \leq C L_\eta\|\kappa\|_{L^\infty}
    \|w-\bar{w}\|_{(W^{1,q},W^{-1,q}_\emptyset)_{\frac{1}{r},r}}.
  \end{align*}
  Finally, the property of maximal parabolic regularity for
  $\AA(\theta_0)$ follows immediately from
  Corollary~\ref{c-intpolcompact}.  \hfill\end{proof}

%
%
%
%
Next we will establish and investigate the right hand hand side of
\eqref{e-pruessloes}.  For doing so, we now turn our attention to the
elliptic equation~\eqref{e-elliptex}.

\begin{lemma}
  \label{l-rhsbilinear} For $q \geq 2$ and $\zeta \in C(\olom)$,
  $\mathfrak a_\zeta(\varphi_1,\varphi_2) :=
  (\sigma(\zeta)\varepsilon\nabla\varphi_1)\cdot\nabla\varphi_2$
  defines a continuous bilinear form $\mathfrak a_\zeta \colon W^{1,q}_\GD \times
  W^{1,q}_\GD \to L^{q/2}$.  Moreover, $(\zeta,\varphi) \mapsto
  \mathfrak a_\zeta(\varphi,\varphi)$ is Lipschitzian over bounded sets in
  $C(\olom) \times W^{1,q}_\GD$.
\end{lemma}

\begin{proof}
  Bilinearity and continuity of each $\mathfrak a_\zeta$ are clear. The second
  assertion follows from a straightforward calculation with the
  resulting estimate
  \begin{align*}\|\mathfrak a_{\zeta_1}(\varphi_1,\varphi_1)-\mathfrak a_{\zeta_2}(\varphi_2,\varphi_2)\|_{L^{q/2}}
    & \leq
    \|\sigma(\zeta_1)-\sigma(\zeta_2)\|_{L^\infty}\|\varepsilon\|_{L^\infty}
    \|\varphi_1\|_{W^{1,q}_\GD}^2 \\& \quad +
    2\|\sigma(\zeta_2)\|_{L^\infty}\|\varepsilon\|_{L^\infty}\|\varphi_1\|_{W^{1,q}_\GD}\|\varphi_1-\varphi_2\|_{W^{1,q}_\GD},
  \end{align*}
  Lipschitz continuity of $\sigma$ and boundedness of the underlying
  sets. \hfill
\end{proof}

Let us draw some further conclusions from
Lemma~\ref{l-multiplierkappa}. For this, we assume
Assumption~\ref{a-reg} for the rest of this chapter.

\begin{theorem} \label{t-contoncont} The mapping
  \begin{equation} \label{e-muli} \underline C(\olom) \ni \phi \mapsto
    (-\nabla \cdot \phi \varepsilon \nabla)^{-1} \in \LL\HH
    (W^{-1,q}_{\GD};W^{1,q}_{\GD} )
  \end{equation}
  is well-defined and even continuous.
\end{theorem}
\begin{proof}
  The well-definedness assertion results from
  Lemma~\ref{l-multiplierkappa}.  The second assertion is implied by
  the first, Remark~\ref {r-einschr} and the continuity of the mapping
  $\LL \HH (X;Y) \ni B \mapsto B^{-1} \in \LL \HH (Y;X)$,
  see~\cite[Ch.~III.8]{schwartz}.  \hfill\end{proof}

\begin{corollary} \label{c-Calpha} Let $\underline{\mathfrak C}
  \subset \underline{C}(\olom)$ be a compact set in $C(\olom)$ which
  admits a common lower positive bound. Then the function
  \[
  \underline{\mathfrak C} \ni \phi \mapsto \JJ(\phi) := \left (-\nabla
    \cdot \phi \varepsilon \nabla \right )^{-1} \in \LL\HH
  (W^{-1,q}_{\GD};W^{1,q}_{\GD} )
  \]
  is bounded and even Lipschitzian. The same holds for
  $\underline{\mathfrak C} \times \mathfrak B \ni (\phi,v) \mapsto
  \JJ(\phi)v \in W^{1,q}_{\GD}$ for every bounded set $\mathfrak B
  \subset W^{-1,q}_{\GD}$.
\end{corollary}
\begin{proof}
  Theorem~\ref{t-contoncont} and the compactness of
  $\underline{\mathfrak C}$ in $C(\olom)$ immediately imply
  boundedness of $\JJ$ on $\underline{\mathfrak C}$. In turn,
  Lipschitz continuity of $\JJ$ is obtained from boundedness and the
  resolvent-type equation
  \begin{multline} \label{e-02} (-\nabla \cdot \phi_1 \varepsilon
    \nabla)^{-1}-(-\nabla \cdot \phi_2 \varepsilon \nabla)^{-1} \\ =
    (-\nabla \cdot \phi_1 \varepsilon \nabla)^{-1}(-\nabla \cdot
    (\phi_2-\phi_1) \varepsilon \nabla) (-\nabla \cdot \phi_2
    \varepsilon \nabla)^{-1}
  \end{multline}
  (read: $A^{-1} - B^{-1} = A^{-1}(B-A)B^{-1}$) and
  Remark~\ref{r-einschr}. 
  Considering the assertion on the combined mapping, boundedness is
  obvious and further we have for $\phi_1,\phi_2 \in
  \underline{\mathfrak C}$ and $v_1,v_2 \in \mathfrak B$:
  \begin{align*}\left\|\JJ(\phi_1)v_1 -
      \JJ(\phi_2)v_2\right\|_{W^{1,q}_{\GD}} 
    & \leq
    \left\|\JJ(\phi_1)-\JJ(\phi_2)\right\|_{\LL(W^{-1,q}_{\GD},W^{1,q}_{\GD})}
    \|v_1\|_{W^{-1,q}_{\GD}} \\ & \qquad +
    \|\JJ(\phi_2)\|_{\LL(W^{-1,q}_{\GD},W^{1,q}_{\GD})}
    \|v_1-v_2\|_{W^{-1,q}_{\GD}}. 
  \end{align*}
  With Lipschitz continuity and boundedness of $\JJ$ over
  $\underline{\mathfrak C}$ and boundedness of $\mathfrak B$, this
  implies the claim. \hfill\end{proof}

\begin{rem} \label{r-explain} At this point we are in the position to
  discuss the meaning of Assumption~\ref{a-reg} in some detail. Under
  Assumption~\ref{a-groeger}~\ref{a-groeger-groeger} for a closed
  subset $\Xi$ of $\pO$, it is known that, even for arbitrary
  measurable, bounded, elliptic coefficient functions $\mu$,
  $(\DD_q,W^{-1,q}_\Xi)_{\tau,1}$ embeds into a H\"older space for
  suitable $\tau$, cf.~\cite[Cor.~3.7]{HDMR08} (for $\DD_q$, see
  Remark~\ref{r-einschr}). In particular, one does \emph{not} need an
  assumption for the ismorphism property between $W^{1,q}_\Xi$ and
  $W^{-1,q}_\Xi$ for this result. The crucial point behind
  Assumption~\ref{a-reg} is to achieve both independence of the
  domains for the operators $-\nabla \phi \mu \nabla$ within a
  suitable class of functions $\phi$, as well as a well-behaved
  dependence on $\phi$ in the space $\mathcal L(\mathcal
  D_q;W^{-1,q}_\Xi)$, cf.\ Lemma~\ref{l-multiplierkappa} and
  Corollaries~\ref{c-quasidomains} and~\ref{c-Calpha}.
\end{rem}

The next lemmata establish the right-hand side in~\eqref{e-pruessloes}
with the correct regularity and properties. Moreover,
Lipschitz continuity with respect to the control $u$ in the elliptic
equation is shown along the way, which will become useful in later
considerations. Recall that $\sigma \colon \R \to \R_+$ is Lipschitzian on
any finite interval by Assumption~\ref{a-assu1}.

\begin{definition} \label{d-ellright} We assign to $\zeta \in
  C(\olom)$ and $v \in W^{-1,q}_{\GD}$ the solution $\varphi_v$ of
  $-\nabla \cdot \sigma (\zeta) \varepsilon \nabla \varphi_v=v$ via
  $\varphi_v = \JJ(\sigma(\zeta))v$ with $\JJ$ as in
  Corollary~\ref{c-Calpha}. Moreover, set \[\Psi_v(\zeta) :=
\mathfrak   a_\zeta(\JJ(\sigma(\zeta))v,\JJ(\sigma(\zeta))v) 
  \] for $\zeta \in C(\olom)$ with $\mathfrak a_\zeta$ as in
  Lemma~\ref{l-rhsbilinear}.
\end{definition}

\begin{lemma} \label{t-righthands}
  Let $\mathfrak C$ be a compact subset of $C(\olom)$ and $\mathfrak
  B$ a bounded set in $W^{-1,q}_{\GD}$.  Then $(v,\zeta) \mapsto
  \Psi_v(\zeta)$ is Lipschitzian from $\mathfrak C \times \mathfrak B$
  into $L^{q/2}$ and the Lipschitz constant of $\zeta \mapsto
  \Psi_v(\zeta)$ is bounded over $v \in \mathfrak B$.
\end{lemma}
\begin{proof}
  For every $\zeta \in \mathfrak C$, the function $\sigma(\zeta)$
  belongs to $\underline C(\olom)$, thus $\JJ(\sigma(\zeta))v$ is
  indeed from $W^{1,q}_{\GD}$ thanks to
  Lemma~\ref{l-multiplierkappa}. Hence, $\Psi_v(\zeta) \in L^{q/2}$ is
  clear by H\"older's inequality. Let us show the Lipschitz property
  of $\Psi$: First, note that Nemytskii operators induced by
  Lipschitz functions preserve compactness in the space of continuous
  functions, and note further that the set of all $\sigma(\zeta)$ for
  $\zeta \in \mathfrak C$ admits a common positive lower bound by the
  Lipschitz property of $\sigma$. Hence, the set $\{ \sigma(\zeta)
  \colon \zeta \in \mathfrak C\}$ satisfies the assumptions in
  Lemma~\ref{l-rhsbilinear} and Corollary~\ref{c-Calpha}. For
  $\zeta_1,\zeta_2 \in \mathfrak C$ and $v_1,v_2 \in W^{-1,q}_{\GD}$,
  we first obtain via Lemma~\ref{l-rhsbilinear}
  \[\|\Psi_{v_1}(\zeta_1)-\Psi_{v_2}(\zeta_2)\|_{L^{q/2}} \leq
  L_\mathfrak a\left(\|\zeta_1-\zeta_2\|_{C(\olom)} +
    \|\JJ(\sigma(\zeta_1))v_1-\JJ(\sigma(\zeta_2))v_2\|_{W^{1,q}_\GD}\right)\]
  and further with Corollary~\ref{c-Calpha}
  \begin{equation*}\|\JJ(\sigma(\zeta_1))v_1-\JJ(\sigma(\zeta_2))v_2\|_{W^{1,q}_\GD}
    \leq L_\JJ\left(\|\sigma(\zeta_1)
      -\sigma(\zeta_2)\|_{C(\olom)}+\|v_1-v_2\|_{W^{-1,q}_\GD}\right).
  \end{equation*}
  The assertion follows since $\sigma$ was
  Lipschitz continuous. Uniformity of the Lipschitz constant of $\zeta
  \mapsto \Psi_v(\zeta)$ is immediate from the previous
  considerations.
  \hfill\end{proof}

Following the strategy outlined above, we will specify the mapping $S$ from Proposition~\ref{p-pruess} 
for our case and show that it satisfies the required conditions. 

\begin{proposition} \label{p-assuS} Let $q \in (3,4)$ be such that
  Assumption~\ref{a-reg} is satisfied, $r > r^*(q)$, and $u \in
  L^{2r}(J;W^{-1,q}_\GD)$. We set \[S(t,\zeta):= \Psi_{u(t)}(\zeta) +
  \alpha\theta_l(t).\] Then $S$ satisfies the conditions from
  Proposition~\ref{p-pruess} for the spaces $X = W^{-1,q}_\emptyset$
  and $Y = W^{1,q}$.
\end{proposition}
\begin{proof} We show that $S(\cdot,0) \in
  L^r(J;W^{-1,q}_\emptyset)$. The function $\alpha\theta_l$ is
  essentially bounded in time with values in $W^{-1,q}_\GD$ by virtue
  of Remark~\ref{r-einbett/ident} and thus poses no problem here. For
  almost all $t \in J$, we further have
  \[
  \left\|\Psi_{u(t)}(0)\right\|_{L^{q/2}} \leq |\sigma(0)|
  \|\varepsilon\|_{L^\infty}
  \|\JJ(\sigma(0)\|_{\LL(W^{-1,q}_{\GD};W^{1,q}_{\GD})}^2
  \|u(t)\|^2_{W^{-1,q}_{\GD}}.
  \]
  Since $u$ is $2r$-integrable in time, this means that
  $\Psi_{u(t)}(0) \in L^{r}(J;L^{q/2})$.  Due to $q >3$ and thus
  $L^{q/2} \embeds W^{-1,q}_\emptyset$ (cf.\ Remark~\ref{r-konsist}),
  we hence have $S(\cdot,0) \in L^r(J;W^{-1,q}_\emptyset)$. \\
  Let us now show the Lipschitz condition~\eqref{e-lipf}. If
  $\mathfrak C \subset (W^{1,q},W^{-1,q}_\emptyset)_{\frac {1}{r},r}$
  is bounded, its closure $\overline{\mathfrak C}$ with respect to the
  $\sup$-norm on $\olom$ forms a compact set in $C(\olom)$ by
  Corollary~\ref{c-intpolcompact}. The desired Lipschitz estimate for
  $S(t,\cdot)$ now follows immediately from Lemma~\ref{t-righthands}.
  \hfill\end{proof}

Note that this is the point where the supposition on the time-integrability of $u$ from Assumption~\ref{a-assu1}~\ref{a-assu1-u}
comes into play. Essentially, $\Psi_{u(t)}(\zeta)$ only admits half
the time-integrability of $u$, but Propositions~\ref{p-assuA}
and~\ref{p-assuS} both require $r > r^*(q)$ to make use of the
(compact) embedding $(W^{1,q},W^{-1,q}_\emptyset)_{\frac1r,r} \embeds
C(\olom)$. Hence, we need more than $2r^*(q)$-integrability for $u$ in time.

Now we have established all ingredients to prove Theorem~\ref{t-locexist}.
For this purpose, let the assumptions of Theorem~\ref{t-locexist} hold. Combining
Propositions~\ref{p-assuA} and~\ref{p-assuS} with
Proposition~\ref{p-pruess}, we obtain a local-in-time solution
$\theta$ of the equation
\[
\theta'(t)+ \AA(\theta(t)) \theta(t) =S(t,\theta(t)), \quad
\theta(T_0)=\theta_0
\]
on $(T_0,T_*)$ with $T_* \in (T_0,T_1]$, such that
\[\theta \in
W^{1,r}(T_0,T_*;W^{-1,q}_\emptyset) \cap L^r(T_0,T_*;W^{1,q}) \embeds
C([T_0,T_*];(W^{1,q},W^{-1,q}_\emptyset)_{\frac{1}{r},r}).\] If $T_* <
T_1$, we may apply Proposition~\ref{p-pruess} again on the interval
$(T_*,T_1)$ with initial value $\theta(T_*) \in
(W^{1,q},W^{-1,q}_\emptyset)_{\frac{1}{r},r}$, thus obtaining another
local solution on a subinterval of $(T_*,T_1)$, ``glue'' the solutions
together and start again (note that $\AA(\theta(t))$ again satisfies
maximal parabolic regularity for every $t \in [T_*,T_1)$ by
Corollary~\ref{c-intpolcompact}). As we may let these intervals of local existence overlap,  
the uniqueness of local solutions by Proposition~\ref{p-pruess} implies that the ``glued'' solution satisfies 
the claimed regularity for the solutions as in~\eqref{e-solreg}. In this way, we either obtain a
solution on the whole prescribed interval $(T_0,T_1)$ or end up with a
maximal interval of existence, denoted by $J_{\max} = (T_0,T_{\max})$,
such that there exists a solution $\theta$ in the above sense on every
interval $(T_0,T_\bullet)$ where $T_\bullet \in J_{\max}$ (or
equivalently $(T_0,T_\bullet] \subset (T_0,T_{\max})$). The maximal
time of existence $T_{\max}$ is characterized by the property that
$\lim_{t \nearrow T_{\max}} \theta(t)$ does not exist in
$(W^{1,q},W^{-1,q}_\emptyset)_{\frac{1}{r},r}$,
see~\cite[Cor.~3.2]{pruess}.

Consider such $T_\bullet \in J_{\max}$. We now {\em define} the
function $\varphi(t)$ for each $t \in (T_0,T_\bullet)$ as the solution
of $-\nabla \cdot \sigma (\theta(t))\varepsilon \nabla \varphi=u(t)$,
that is, \begin{equation}\varphi(t) :=
  \JJ(\sigma(\theta(t)))u(t).\label{e-varphidef}\end{equation} Then
$\varphi$ indeed belongs to $L^{2r}(T_0,T_\bullet;W^{1,q}_{\GD})$,
since $\JJ(\sigma(\theta(t))$ is uniformly bounded in
$\LL(W^{-1,q}_{\GD};W^{1,q}_{\GD})$ over $[T_0,T_\bullet]$ due to the
compactness of the set $\{\theta(t) \colon t \in [T_0,T_\bullet]\}$ in
$C(\olom)$ (cf.\ Corollary~\ref{c-intpolcompact} and
Corollary~\ref{c-Calpha}), and $u$ was from
$L^{2r}(J;W^{-1,q}_{\GD})$.

Obviously, $(\theta,\varphi)$ is then a solution of the
thermistor-problem on $(T_0,T_\bullet)$ in the spirit of
Definition~\ref{d-loes} as claimed in Theorem~\ref{t-locexist}.

  We end this chapter with
  some explanations why the chosen setting in spaces of the kind
  $W^{-1,q}_\emptyset$ and $W^{-1,q}_\GD$ with $q > 3$ is adequate for the problem
  under consideration. 

  Let us inspect the requirements on the spaces in which the
  equations are formulated. Clearly, they need to contain Lebesgue
  spaces on $\Omega$ as well as on the boundary $\Gamma$ (or on a
  subset of the boundary like $\GN$), in order to
  incorporate the nonhomogenenous Neumann boundary data present in
  both equations. The boundary conditions should be reflected
  by the formulation of the equations in an adequate way,
  cf.~Remark~\ref{r-Loesung}~\ref{r-boundcond}. These demands already strongly prejudice spaces of type
  $W^{-1, q_p}_\emptyset$ for the parabolic equation and
  $W^{-1,q_e}_\GD$ for the elliptic equation with probably different
  integrability orders $q_p$ and $q_e$ for each equation. Finally, in order to
  treat the nonlinear parabolic equation, we need maximal parabolic
  regularity for the second order divergence operators $\AA(\zeta)$ over
  $W^{-1,q_p}_\emptyset$, which is generally available by
  Corollary~\ref{c-intpolcompact}~\ref{c-intpolcompact-mpr}
  or~\cite[Thm.~5.16/Rem.~5.14]{hal/reh} in a general context.

  Further, aiming at continuous solutions $\theta$, which are
  needed for having fulfillable Constraint Qualifications
  for~\eqref{P} in the
  presence of state constraints, it is necessary that the domain
  $\DD_{q_p}(\zeta)$ of the
  differential operators $\AA(\zeta)$, cf.\ Remark~\ref{r-einschr}, embeds into the space of
  continuous functions on $\overline Q$. But it is known that solutions $y$ to equations
  $-\nabla \cdot \mu\nabla y = f$ for $\mu \in L^\infty(\Omega,M_n)$ elliptic with $f \in W^{-1,n}_\emptyset$,
  where $n$ denotes the space dimension, may in general even be
  unbounded, see~\cite[Ch.~1.2]{lady}. On the other hand,
  $\DD_{q_p}(\zeta)$ embeds into a H\"older space if $q_p > 3$, see Remark~\ref{r-einschr}. These two
  facts make the requirement $q_p  > n = 3$ expedient. Let us now
  assume that the elliptic
  equation admits solutions whose gradient is integrable up to some order $q_g$. Then the right hand side in the parabolic equation
  prescribes $q_g \geq \frac{6q_p}{q_p + 3}$ in order to have the
  embedding $L^{
q_g/2} \embeds W^{-1,q_p}_\emptyset$. From the
  requirement $q_p > 3$ then follows $q_g > 3$ as well, i.e., the
  elliptic equation must admit $W^{1,q_g}_\GD$-solutions with $q_g >
  3$. With right-hand sides in $W^{-1,q_e}_\GD$, the best possible constellation
  is thus $q_e = q_g > 3$ again. Having $q_e$ and $q_p$ both in the
  same range, we simply choose $q = q_e = q_p > 3$. 

  Moreover, in
  order to actually have $W^{1,q}_\GD$-solutions to the elliptic
  equations for all right-hand sides from $W^{-1,q}_\GD$, the operator $-\nabla \cdot \sigma(\zeta)\varepsilon\nabla$ must be a
  topological isomorphism between $W^{1,q}_\GD$ and $W^{-1,q}_\GD$. It is
  also a well-established fact that solutions to elliptic equations
  with bounded and coercive, but discontinuous coefficient functions
  may admit almost arbitrarily poor integrability properties for
  gradients of their solutions, see~\cite{mercier} and~\cite[Ch.~4]{e/r/s}. Under
  Assumption~\ref{a-reg}, we know that this is not the case for $-\nabla
  \cdot \varepsilon \nabla$ over $W^{-1,q}_\GD$, but it is clear that
  it is practically impossible to guarantee this also for the
  operators $-\nabla \cdot \sigma(\zeta)\varepsilon \nabla$ for all
  $\zeta$, if $\sigma(\zeta)$ is discontinuous in general. However, from
  Lemma~\ref{l-multiplierkappa} we know that if $\sigma(\zeta)$ if
  uniformly continuous on $\Omega$, then the
  isomorphism property carries over. This shows that
  continuous solutions for the parabolic equation are also needed
  purely from an analytical point of view, without the considerations
  coming from the optimal control problem, and also explains why
  Assumption~\ref{a-reg} is, in a sense, a ``minimal'' assumption.

\section{Global solutions and optimal control}\label{s-oc}

The setting and results of \S~\ref{s-statesys} are assumed as given
from now on, i.e., we consider the assumptions of Theorem~\ref{t-locexist}
to be fulfilled and fixed, that means, $q > 3$ and $r > r^*(q)$ are given from now on. 
In particular, for every $u \in
L^{2r}(J;W^{-1,q}_{\GD})$, there exists a local solution $\theta_u$
such that $\theta_u \in W^{1,r}(T_0,T_\bullet;W^{1,q}) \cap
L^r(T_0,T_\bullet;W^{-1,q}_\emptyset)$ for every $T_\bullet \in
J_{\max}(u)$, the maximal interval of existence for a given control
$u$. We consider $\varphi \in L^{2r}(T_0,T_\bullet;W^{1,q}_\GD)$ to
be given in dependence of $u$ and $\theta_u$ as
in~\eqref{e-varphidef}. Due to $q >3$ and $r > r^*(q)$, each
solution $\theta_u$ is H\"older-continuous on
$[T_0,T_{\bullet}]\times\olom$, cf.\
Corollary~\ref{c-intpolcompact}~\ref{c-intpolcompact-mprspace}.


\begin{rem}
  As noted above, if the solution $\theta_u$ for a given control $u$
  does \emph{not} exist on the whole time interval $J$, there exists
  $T_{\max}(u) \leq T_1$, the maximal time of existence, such that
  $\lim_{t \nearrow T_{\max}(u)} \theta_u(t)$ does not exist in
  $(W^{1,q},W^{-1,q}_\emptyset)_{\frac{1}{r},r}$. For a proof and an
  equivalent formulation in the maximal regularity-norm,
  see~\cite[Cor.~3.2]{pruess}.
\end{rem}

Our aim in the following sections is to characterize the set of control functions
which admit a solution on the \emph{whole} time interval. 
These control functions will be called ``global controls'', see Definition~\ref{d-controltostate}. 
In view of the state constraints and the end time observation in the objective of~\eqref{P}, 
it is natural to restrict the optimal control problem to the set of global controls. 
Our characterization of this set will then allow to establish the existence of 
(globally) optimal controls. Let us give a brief roadmap of the upcoming analysis:
\begin{itemize}
 \item We first show that the set of global controls is not empty, see
   Proposition~\ref{p-hj1}, and that it is an open set, cf.\ Theorem~\ref{t-open}. This property will be of major importance for the derivation of 
 meaningful optimality conditions in Section~\ref{subs-noc}.
 \item For global controls one can define a control-to-state operator
   in function spaces on the whole time interval, see Definition~\ref{d-controltostate}. 
 The proof of Theorem~\ref{t-open} features an application of the
 implicit function theorem and thereby shows that the control-to-state operator is Fr\'{e}chet-differentiable, which is also essential for the derivation 
 of necessary optimality conditions.
 \item We then turn to the existence of optimal controls. The arguments follow the 
 classical direct method of the calculus of variations, see Theorem~\ref{t-existenceoptcont}. 
 To this end, we need to establish a closedness result for the 
 set of global controls in Theorem~\ref{t-globalcompact}. This result requires a certain 
 boundedness of the gradient of the temperatures which is ensured by the 
 second addend in the objective in~\eqref{P}. To pass to the limit in 
 the thermistor system, we additionally need that the control space,
 induced by the third term in the objective functional, 
 compactly embeds into $L^{2r}(J;W^{1,q}_{\Gamma_D})$. This issue is addressed in 
 Proposition~\ref{p-compactcontrolspace}.
 \item Finally, Section~\ref{subs-noc} is devoted to the derivation of necessary optimality conditions. 
 As the set of global controls is open, the standard generalized Karush-Kuhn-Tucker theory applies, 
 see Theorem~\ref{t-lagrangeexistence}. We then introduce an adjoint system in Definition~\ref{d-adjointsys} and~\ref{d-abstractadjointsys}
 and show by means of a classical duality argument that this system admits a unique solution, 
 cf.\ Theorem~\ref{t-weakadjointexists}. 
 This allows to reformulate the necessary conditions in terms of a qualified 
 optimality system, see Theorem~\ref{t-neccondfinal}.
\end{itemize}

\subsection{Global solutions and existence of optimal controls}\label{subs-existenceoc}

In~\cite{AC94}, Antontsev and Chipot show that it is possible to
give concrete conditions under which the solution to a
thermistor-like problem does not exist globally. While the authors
of~\cite{AC94} consider a slightly different setting (in particular no
Robin boundary conditions for the parabolic equation), we devote a
subsection to the question whether there is any relevant
characterization of \emph{global} controls $u$, i.e., controls such
that the corresponding solution $(\theta_u,\varphi_u)$ does exist on
the whole (prescribed) interval $J = (T_0,T_1)$.

We make the following assumption for the rest of this paper:

\begin{assu}\label{a-sigmanemytskii}
  \begin{enumerate}
  \item The functions $\eta$ and $\sigma$, each mapping $ \R \to
    \R_+$, are continuously differentiable. The derivatives $\eta'$
    and $\sigma'$ are each bounded and Lipschitz continuous on bounded
    sets.
  \item In addition to Assumption~\ref{a-groeger}, we from now on require that 
   $\Omega \cup \Gamma_D$ satisfies the volume-conservation condition from 
   Definition~\ref{def:groeger}~\ref{a-groeger-volume}.
  \end{enumerate}
\end{assu}

\begin{definition}\label{d-controltostate}
  We call a control $u \in L^{2r}(J;W^{-1,q}_{\GD})$, $r> r^*(q)$, a
  \emph{global control} if the corresponding solution $\theta$ exists
  on the whole prescribed interval $(T_0,T_1)$ and denote the set of
  global controls by $\UU_g$. Moreover, we define the
  \emph{control-to-state operator}
  \[\s \colon \UU_g \ni u \mapsto \s(u) = \theta_u \in
  W^{1,r}(J;W^{-1,q}_\emptyset) \cap L^r(J;W^{1,q})\] on $\UU_g$.
\end{definition}

Let us firstly show that the previous definition is in fact
meaningful in the sense that $\UU_g \neq \emptyset$. 
The natural candidate for a global control is $u \equiv 0$. One readily
observes that the control $u \equiv 0$ leads to the solution $\varphi
\equiv 0$ for the elliptic equation~\eqref{e-elliptex}, hence the
right-hand side in the parabolic equation reduces to $\alpha
\theta_l(t)$ in this case. Indeed, we will show that there exists a
global solution $\theta_{u=0}$ to the equation
\begin{equation}\partial_t \theta + \AA(\theta)\theta =
  \alpha\theta_l, \quad
  \theta(T_0) = \theta_0.\label{e-control0}
\end{equation}

In order to obtain a global solution to~\eqref{e-control0}, we need
the volume-conservation condition. Under this additional assumption,
the following result has been shown in~\cite[Thm.~5.3]{hajo15}. Note
that the case of $\Omega \cup \GD$ regular is only a special case of
the admissible geometries in~\cite{hajo15}.

\begin{proposition} \label{p-hj1} Assume that $\O \cup \Xi$ is regular
  and in addition satisfies the volume-conservation condition.
  Let $\mu$ be a coefficient function on $\O$, measurable, bounded,
  elliptic. Assume that $\phi:\R \to [\underline\phi,\overline\phi]$,
  where $\underline\phi > 0$, is Lipschitz continuous on bounded
  sets. Suppose further that
  \[
  -\nabla \cdot \mu \nabla :W^{1,q}_\Xi \to W^{-1,q}_\Xi
  \]
  is a topological isomorphism for some $q > 3$. Let $w_0$ be from
  $(W^{1,q}_\Xi,W^{-1,q}_\Xi)_{\frac1r,r}$ with $r > r^*(q) = \frac{2q}{q-3}$. 
  Then, for every $f \in L^{r}(J;W^{-1,q}_\Xi)$, there exists a unique \emph{global} solution
  $w$ of the quasilinear equation
  \begin{equation} \label{e-quasil} w' - \nabla \cdot \phi(w) \mu
    \nabla w =f, \quad w(T_0)=w_0,
  \end{equation}
  which belongs to $W^{1,r}(J;W^{-1,q}_\Xi) \cap
  L^{r}(J;W^{1,q}_\Xi)$.
\end{proposition}

With $w_0 = \theta_0$, $\Xi = \emptyset$, $\phi = \eta$, $\mu
= \kappa$ and $f = \alpha\theta_l$, we may use Proposition~\ref{p-hj1}
to ensure the existence of a \emph{global} solution of~\eqref{e-control0} in the
sense of Definition~\ref{d-loes} under Assumption~\ref{a-reg} -- in
particular, $0 \in \UU_g$ follows. In~\cite{hajo15}, Proposition~\ref{p-hj1} is proven
for the case where the differential operator consists of the
divergence-gradient operator only. However, it is clear that the result extends to 
the operators of the form $\AA$ including the boundary form since the latter 
is relatively compact with respect to the main part, cf.\ Corollary~\ref{c-quasidomains} and
the reference there, see also~\cite[Lem.~5.15]{hal/reh}.

The next theorem establishes continuous differentiability of the
control-to-state operator $\s$. Given a control $u$, we use
$\varphi_u$ for the associated solution of the elliptic equation with
$u$ on the right-hand side, cf.~\eqref{e-varphidef}.

\begin{theorem} \label{t-open} 
  Let $r > r^*(q)$ be given. Then the set of global controls $\UU_g$ forms an open set in $L^{2r}(J;W^{-1,q}_{\GD})$. Moreover,
  the control-to-state operator $\s$ is continuously
  differentiable. For every $h \in L^{2r}(J;W^{-1,q}_\GD)$, its derivative $\zeta_h = \s'(u)h \in
  W^{1,r}(J;W^{-1,q}_\emptyset) \cap L^r(J;W^{1,q})$ is given by the
  unique solution of the equation
  \begin{multline}\partial_t \zeta + \AA(\theta_u)\zeta =
    \left(\sigma'(\theta_u)\zeta\varepsilon\nabla\varphi_u\right)\cdot\nabla\varphi_u
    + \nabla \cdot \eta'(\theta_u)\zeta\kappa\nabla\theta_u \\ -2
    \left(\sigma(\theta_u)\varepsilon\nabla\varphi_u\right)\cdot\nabla
    \left[\JJ(\sigma(\theta_u))\left(-\nabla \cdot
        \sigma'(\theta_u)\zeta \varepsilon \nabla\varphi_u + h
      \right)\right], \label{e-linearized}
  \end{multline} which has to hold for almost every $t \in J$ in the
  space $W^{-1,q}_\emptyset$, with $\zeta(T_0) = 0$.
\end{theorem}

\begin{proof}
  Let $\bar{u} \in L^{2r}(J;W^{-1,q}_{\GD})$ be global, i.e.,
  the associated solution $\theta_{\bar{u}} =: \bar{\theta}$ exists on
  the whole time horizon $(T_0,T_1)$. We intend to apply the implicit
  function theorem. To this end, we show that the mapping
  \begin{multline*}\BB \colon \left (W^{1,r}(J;W^{-1,q}_\emptyset) \cap
      L^r(J;W^{1,q})\right ) \times L^{2r}(J;W^{-1,q}_{\GD})
    \\ \to L^r(J;W^{-1,q}_\emptyset) \times
    (W^{1,q},W^{-1,q}_\emptyset)_{\frac{1}{r},r},\end{multline*} where
  \begin{equation}\BB(\theta,u) = \left(\partial_t \theta +
      \AA(\theta)
      \theta - \Psi_u(\theta) - \alpha \theta_l, \theta(T_0)-\theta_0\right),\label{e-differentialoperator}\end{equation} is continuously differentiable in
  $(\bar{\theta},\bar{u})$, and that the partial derivative
  $\partial_\theta \BB(\bar{\theta},\bar{u})$ is continuously
  invertible. Note that $\BB(\bar{\theta},\bar{u}) = 0$. The term
  $\alpha\theta_l$ does not depend neither on $u$ nor on $\theta$ and is thus neglected for the rest of
  this proof. Let us first consider the partial derivative with
  respect to $u$: For each
  $\theta \in C(\olQ)$, the mapping
  \[L^{2r}(J;W^{-1,q}_{\GD})^2 \ni (u,v) \mapsto
  \left(\sigma(\theta)\varepsilon\nabla\varphi_u(\theta)\right)\cdot\nabla\varphi_v(\theta)
  \in L^r(J;L^{q/2})\] gives rise to a continuous symmetric bilinear
  form $b_\theta(u,v)$ (cf.\ also Lemma~\ref{l-rhsbilinear}), since
  for fixed $\theta \in C(\olQ)$ we have
  \begin{multline*}\|b_\theta(u,v)\|_{L^r(J;L^{q/2})} \leq
    \|\sigma(\theta)\|_{C(\olQ)}\|\varepsilon\|_{L^\infty}\|\JJ(\sigma(\theta))\|_{C(\olJ;\LL(W^{-1,q}_{\GD},W^{1,q}_{\GD}))}^2
    \\ \cdot \|u\|_{L^{2r}(J;W^{-1,q}_{\GD})}\|v\|_{L^{2r}(J;W^{-1,q}_{\GD})}.
  \end{multline*} 
  Accordingly, $u \mapsto
  \Psi_u(\bar{\theta}) = b_{\bar{\theta}}(u,u)$ is continuously
  differentiable, and its derivative in $\bar{u}$ is given by $h
  \mapsto 2b_{\bar{\theta}}(\bar{u},h)$. The second component of $\BB$
  is independent of $u$. Next, we treat the derivative of $\BB$
  w.r.t.\ $\theta$. First, note that, due to
  Assumption~\ref{a-sigmanemytskii}, the Nemytskii operator $\theta
  \mapsto \eta(\theta)$ is continuously differentiable from $C(\olQ)$
  to $C(\olQ)$ and its derivative in $\bar{\theta}$ is given by $h
  \mapsto \eta'(\bar{\theta})h$. With Remark~\ref{r-einschr}, we thus
  find that the derivative of the function $\theta \mapsto \partial_t
  \theta +\AA(\theta)\theta$ as a mapping from
  $W^{1,r}(J;W^{-1,q}_\emptyset) \cap L^r(J;W^{1,q})$ to
  $L^r(J;W^{-1,q}_\emptyset)$ in the point $\bar{\theta}$ is given
  by \begin{equation}h \mapsto \partial_t h - \nabla \cdot
    \eta(\bar{\theta})\kappa\nabla h + \tilde{\alpha}h - \nabla \cdot
    \eta'(\bar{\theta})h\kappa\nabla\bar{\theta} = \partial_t h +
    \AA(\bar{\theta})h - \nabla \cdot
    \eta'(\bar{\theta})h\kappa\nabla\bar{\theta}. \label{e-differentialderivative}
  \end{equation}
  We turn to $\theta \mapsto \Psi_{\bar{u}}(\theta)$. As above, due to
  Assumption~\ref{a-sigmanemytskii}, $\theta \mapsto \sigma(\theta)$
  is continuously differentiable as a mapping from $C(\olQ)$ to
  $C(\olQ)$ and with derivative $h \mapsto \sigma'(\bar{\theta})h$ (in
  a point $\bar{\theta}$). Further, recall that the derivative of the
  (continuously differentiable) mapping $\LL(X;Y) \ni A \mapsto A^{-1}
  \in \LL(Y;X)$ in $A$ is given by $H \mapsto -A^{-1}HA^{-1}$. The
  chain rule and Remark~\ref{r-einschr} thus yield continuous
  differentiability of $\theta \mapsto \JJ(\sigma(\theta))$ as a
  mapping from $C(\olJ;C(\olom))$ to
  $C(\olJ;\LL(W^{1,q}_\GD;W^{-1,q}_{\GD}))$ with the
  derivative
  \[\left[\left(\JJ \circ \sigma\right)'(\bar{\theta})\right]h =
  -\JJ(\sigma(\bar{\theta}))\left[-\nabla \cdot \sigma'(\bar{\theta})h
    \varepsilon \nabla\right]\JJ(\sigma(\bar{\theta})).\] Hence,
  $\theta \mapsto \varphi_{\bar{u}}(\theta) =
  \JJ(\sigma(\theta))\bar{u}$ is also continuously differentiable,
  considered as a mapping from $C(\olJ;C(\olom))$ to
  $L^{2r}(J;W^{1,q}_{\GD})$. Continuous differentiability of
  the function given by $C(\olJ;C(\olom)) \ni \theta \mapsto
  \Psi_{\bar{u}}(\theta) \in L^r(J;L^{q/2}) \embeds
  L^r(J;W^{-1,q}_\emptyset)$ is now straightforward and its derivative
  in $\bar{\theta}$ is given by
  \begin{multline}\label{e-Psiderivative}\left[\partial_\theta
      \Psi_{\bar{u}}(\bar{\theta})\right]h = -2
    \left(\sigma(\bar{\theta})\varepsilon\nabla\left[\JJ(\sigma(\bar{\theta}))\bar{u}\right]\right)\cdot\nabla
    \left[\left(\left[\left(\JJ \circ
            \sigma\right)'(\bar{\theta})\right]h\right)\bar{u}\right]
    \\ +
    \left(\sigma'(\bar{\theta})h\varepsilon\nabla\left[\JJ(\sigma(\bar{\theta}))
        \bar{u}\right]\right)\cdot\nabla\left[\JJ(\sigma(\bar{\theta}))\bar{u}\right].
  \end{multline}
  The second component of $\BB$, i.e., $\theta \mapsto
  \theta(T_0)-\theta_0$, is affine-linear and continuous from the
  maximal regularity space into
  $(W^{1,q},W^{-1,q}_\emptyset)_{\frac{1}{r},r}$ and as such has the
  derivative $h \mapsto h(T_0)$. It remains to show the continuous
  invertibility of $\partial_\theta \BB(\bar{\theta},\bar{u})$. For
  this, we identify for almost every $t \in (T_0,T_1)$ and $h \in
  C(\olJ;C(\olom))$ as follows:
  \[ B(t)h(t) = \left(\left[\partial_\theta
      \Psi_{\bar{u}}(\bar{\theta})\right]h\right)(t) + \nabla\cdot
  \eta'(\bar{\theta}(t))h(t)\kappa\nabla\bar\theta(t),\] such that
  $B(t)$ is from $\LL(C(\olom);W^{-1,q}_\emptyset)$ and $t \mapsto
  B(t) \in
  L^r(J;\LL(C(\olom);W^{-1,q}_\emptyset))$. Combining~\eqref{e-differentialderivative}
  and~\eqref{e-Psiderivative}, in order to prove that $\mathcal
  B_\theta$ is continuously invertible we need to show that the
  equation
  \begin{equation}\partial_t \xi(t) + \AA(\bar{\theta}(t))\xi(t) =
    B(t)\xi(t) + f(t),
    \qquad \xi(T_0) = \xi_0\label{e-linearizedperturb}\end{equation}
  has a unique solution $\xi
  \in W^{1,r}(J;W^{-1,q}_\emptyset)\cap L^r(J;W^{1,q})$ for every $f
  \in L^r(J;W^{-1,q}_\emptyset)$ and $\xi_0 \in (W^{1,q},W^{-1,q}_\emptyset)_{\frac{1}{r},r}$.
  This, however, is exactly what is obtained
  by~\cite[Cor.~3.4]{pruess}, hence we have 
  \[
  \partial_\theta \BB(\bar{\theta},\bar{u}) \in \LL\HH \bigl
  (W^{1,r}(J;W^{-1,q}_\emptyset) \cap
  L^r(J;W^{1,q});L^r(J;W^{-1,q}_\emptyset)\bigr ) \times
  (W^{1,q},W^{-1,q}_\emptyset)_{\frac{1}{r},r})).
  \]
  Thus, all requirements for the implicit function theorem are
  satisfied, which yields neighbourhoods $\mathfrak{V}_{\bar{u}}$ of
  $\bar{u}$ in $L^{2r}(J;W^{-1,q}_{\GD})$ and
  $\mathfrak{V}_{\bar{\theta}}$ of $\bar{\theta}$ in the maximal
  regularity space, such that there exists a continuously
  differentiable mapping $\Phi \colon \mathfrak{V}_{\bar{u}} \to
  \mathfrak{V}_{\bar{\theta}}$ with $\BB(\Phi(u),u) =
  \BB(\bar{\theta},\bar{u}) = 0$ for all $u \in
  \mathfrak{V}_{\bar{u}}$. This shows that the set of global controls is
  open. Moreover, $\Phi$ locally coincides with the
  control-to-state operator $u \mapsto \s(u)$, which implies
  continuous differentiability for the latter.  \\ The stated
  expression for $\s'(u)h$ is obtained by differentiating the
  (constant) function $u \mapsto \BB(\s(u),u)$. From the second
  component, we then find $(\s'(u)h)(T_0) = 0$ in
  $(W^{1,q},W^{-1,q}_\emptyset)_{\frac{1}{r},r}$ for all $h$, and the
  chain rule yields
  \[\s'(u)h = -[\partial_\theta \BB(\s(u),u)]^{-1} \partial_u
  \BB(\s(u),u)h,\] meaning exactly that $\s'(u)h$ is the unique
  solution to the problem~\eqref{e-linearizedperturb} with right-hand
  side $f = - \partial_u \BB(\s(u),u)h$ and initial value
  $0$. Inserting all formulas, we obtain the equation stated in the
  theorem.  \hfill
\end{proof}

\begin{rem}\label{r-linearizedtwoequations}
  One may split the equation solved by $\zeta_h = \s'(u)h$ in the
  previous Theorem~\ref{t-open} back into two equations: Introducing
  \[\Phi(\zeta) := \JJ(\sigma(\theta_u))\left(-\nabla \cdot
    \sigma'(\theta_u)\zeta \varepsilon \nabla\varphi_u + h \right) \in
  L^{2r}(J;W^{1,q}_\GD),\] we find that, for every $h \in
  L^{2r}(J;W^{-1,q}_\GD)$, the pair $(\zeta,\pi) :=
  (\s'(u)h,\Phi(\s'(u)h))$ is the \emph{unique solution} of the system
  \begin{align*} \partial_t \zeta + \AA(\theta_u)\zeta & =
    \left(\sigma'(\theta_u)\zeta\varepsilon\nabla\varphi_u\right)\cdot\nabla\varphi_u
    + \nabla \cdot \eta'(\theta_u)\zeta\kappa\nabla\theta_u +2
    \left(\sigma(\theta_u)\varepsilon\nabla\varphi_u\right)\cdot\nabla
    \pi \\
    - \nabla \cdot \sigma(\theta_u)\varepsilon\nabla \pi & = -\nabla
    \cdot \sigma'(\theta_u)\zeta \varepsilon\nabla\varphi_u + h
  \end{align*}
  with $\zeta(T_0) = 0$ (the first equation is supposed to hold in
  $W^{-1,q}_\emptyset$, the second one in $W^{-1,q}_\GD$, each for
  almost all $t \in J$). These equations are exactly the
  \emph{linearized state system} for~\eqref{e-paraex}
  and~\eqref{e-elliptex}. This also shows, expectedly, that from a
  functional-analytical point of view, it makes no difference working
  with $\theta$ only and considering $\varphi$ as a function obtained
  by $\theta$, instead of considering both functions at once.
\end{rem}

Combining Theorem~\ref{t-open} with Proposition~\ref{p-hj1} as explained above, we obtain the following

\begin{corollary} \label{c-smallball} There is always a neighbourhood
  $\mathfrak V_0$ of $0$ in $L^{2r}(J;W^{-1,q}_{\GD})$,
  containing only global controls, i.e., $\mathfrak{V}_0 \subseteq
  \UU_g$.
\end{corollary}

Now that we have established a certain richness of global controls, we
turn to the question of existence of an optimal control
of~\eqref{P}. Following the standard direct method of the calculus of variations, 
one soon encounters the
problem of lacking uniform boundedness in a suitable space for solutions $\theta_u$ 
associated to a minimizing sequence of global controls $u$, which is a common obstacle
to overcome when treating quasilinear equations. To circumvent this,
we use Proposition~\ref{p-divproperty}~\ref{p-divproperty-groeger} to
show that the solutions $\theta_u$ are uniformly bounded in
$W^{1,s}(J;W^{-1,\varsigma}_\emptyset)$, where $\varsigma \leq 3 < q$ (in general only $\varsigma \sim \frac{3}{2}$) and $s$ is the exponent from the second 
addend in the objective function in~\eqref{P}.
As this term in the objective gives an additional bound in $L^{s}(J;W^{1,q})$, 
we can employ Corollary~\ref{c-intpolcompact} to ``lift'' this boundedness result 
to a H\"older space, which is suitable for passing to the limit with a minimizing sequence.
However, in order to apply Corollary~\ref{c-intpolcompact}, the exponent $s$ has to be sufficiently large. 
The precise bound for $s$ is characterized by the following

\begin{definition}\label{def:barq}
 Let $\mathfrak{q} \in (2, \min\{q_0,3\}]$ be given, where $q_0$ is the number from 
 Proposition~\ref{p-divproperty}~\ref{p-divproperty-groeger}, and set 
 $\varsigma := \frac{3\mathfrak{q}}{6-\mathfrak{q}}$. Then we define the number $\bar r(q,\varsigma) > 0$ by
 \begin{equation}
  \bar r(q,\varsigma) := \frac{2q}{q-3} \left(1-\frac{3}{q}+\frac{3}{\varsigma}\right).
 \end{equation}
\end{definition}

On account of $\varsigma \leq 3 < q$ it follows that $\bar r(q,\varsigma) > r^*(q) = \frac{2q}{q-3}$. 
Therefore, for a given number $s > \bar r(q,\varsigma)$, the previous results, 
in particular the assertions of Theorem~\ref{t-locexist}, Theorem~\ref{t-open}, 
and Corollary~\ref{c-0feasible}, hold with $r = s$. 
The next theorem precisely elaborates the argument depicted before Definition~\ref{def:barq}:

\begin{theorem} \label{t-globalcompact} Let $\UU \subseteq \UU_g$ be
  bounded in $L^{2s}(J;W^{-1,q}_{\GD})$ with $s > \bar r(q,\varsigma)$ 
  and suppose in addition that the associated set 
  of solutions $\KK = \{\theta_u \colon u \in \UU\}$ is bounded in
  $L^{s}(J;W^{1,q})$. Then $\KK$ is even compact in $C(\olQ)$ and
  the closure of $\UU$ in $L^{2s}(J;W^{-1,q}_{\GD})$ is
  still contained in $\UU_g$.
\end{theorem}

As indicated above, the second addend in the objective functional together with
the state constraints will guarantee 
the bound in $L^{s}(J;W^{1,q})$ for the minimizing sequence, see the proof of Theorem~\ref{t-existenceoptcont} 
below.

\begin{proofof}{Theorem~\ref{t-globalcompact}} We show that $\KK$ is bounded in a suitable
  maximal-regularity-like space. To this end, we first investigate the
  right-hand side in the parabolic equation~\eqref{e-paraex}.  Denote
  by $(\theta_u,\varphi_u)$ the solution for a given $u \in
  \UU$. Thanks to Assumption~\ref{a-assu1}~\ref{a-assu1-sigmaeta},
  Proposition~\ref{p-divproperty}~\ref{p-divproperty-groeger} 
  shows that $-\nabla \cdot \sigma(\theta)\varepsilon\nabla$ is a
  topological isomorphism between $W^{1,\mathfrak{q}}_{\GD}$ and
  $W^{-1,\mathfrak{q}}_{\GD}$ with
  \begin{equation}\label{eq:invnormest}
   \sup_{\theta\in\KK}
   \|\left(-\nabla\cdot\sigma(\theta)
   \varepsilon\nabla\right)^{-1}\|_{\LL(W^{-1,\mathfrak{q}}_{\GD};W^{1,\mathfrak{q}}_{\GD})}
   < \infty.  
  \end{equation}   
  Hence, for every $u \in \UU$ there exists a unique $\psi
  = \psi_u \in L^{2s}(J;W^{1,\mathfrak{q}}_{\GD})$ such
  that
  \[\psi_u(t) =
  \left(-\nabla\cdot\sigma(\theta_u(t))\varepsilon\nabla\right)^{-1}
  u(t) \quad \text{in } W^{1,\mathfrak{q}}_{\GD}\] for almost every $t
  \in (T_0,T_1)$, and 
  \[\sup_{u \in\UU} \|\psi_u\|_{L^{2s}(J;W^{1,\mathfrak{q}}_{\GD})} < \infty.\] 
  Since
  $W^{1,q}_{\GD} \embeds W^{1,\mathfrak{q}}_{\GD}$ and, by uniqueness
  of $\psi_u$, we in particular obtain $\varphi_u = \psi_u$, such that
  the family $\varphi_u$ is bounded in $L^{2s}(J;W^{1,\mathfrak{q}}_{\GD})$ as well. Estimating as in
  Lemma~\ref{t-righthands}, we find that also
  \[\sup_{u \in\UU}
  \left\|\left(\sigma(\theta_u)\varepsilon\nabla\varphi_u\right)
  \cdot\nabla\varphi_u\right\|_{L^{s}(J;L^{\mathfrak{q}/2})}
  < \infty.\] Using the boundedness assumption on $\KK$ in
  $L^{s}(J;W^{1,q})$, both the family of functionals
  $\tilde{\alpha}\theta_u$ and, here also employing boundedness of
  $\eta$, the divergence-operators
  $-\nabla\cdot\eta(\theta_u)\kappa\nabla\theta_u$ are uniformly
  bounded over $\UU$, i.e.,
  \[\sup_{u \in \UU} \|\nabla\cdot
  \eta(\theta_u)\kappa\nabla\theta_u\|_{L^{s}(J;W^{-1,q}_\emptyset)}
  + \|\tilde{\alpha}\theta_u\|_{L^{s}(J;W^{-1,q}_\emptyset)} <
  \infty.\] Sobolev embeddings give the embedding $L^{\mathfrak{q}/2}
  \embeds W^{-1,\varsigma}_\emptyset$ for $\varsigma =
  \frac{3\mathfrak{q}}{6-\mathfrak{q}}$, and certainly
  $W^{-1,q}_\emptyset \embeds W^{-1,\varsigma}_\emptyset$ due to $q >
  \varsigma$. Hence,
  \[\partial_t \theta_u = \nabla\cdot
  \eta(\theta_u)\kappa\nabla\theta_u - \tilde{\alpha}\theta_u +
  \left(\sigma(\theta_u)\varepsilon\nabla\varphi_u\right)\cdot\nabla\varphi_u
  + \alpha\theta_l\] is uniformly bounded over $\UU$ in
  $L^{2s}(J;W^{-1,\varsigma}_\emptyset)$. This shows that $\KK$ is
  bounded in the space $W^{1,s}(J;W^{-1,\varsigma}_\emptyset) \cap
  L^{s}(J;W^{1,q})$. By Corollary~\ref{c-intpolcompact}, $\KK$ is
  then also bounded in a H\"older space and thus a (relatively)
  compact set in $C(\olQ)$.

  Next, let us show that the limit of a convergent sequence in $\UU$
  is still a global control. Denote by $(u_n) \subset \UU$ such a
  sequence, converging in $L^{2s}(J;W^{-1,q}_\GD)$ to the
  limit $\bar{u}$. We call the associated states $(\theta_n) :=
  (\theta_{u_n})$. Compactness of $\KK$ as shown above gives a
  subsequence of $(\theta_n)$, called $(\theta_{n_k})$, which
  converges to some $\bar{\theta}$ in
  $C(\olQ)$. Lemma~\ref{t-righthands} shows that
  $\Psi_{u_{n_k}}(\theta_{n_k}) \to \Psi_{\bar{u}}(\bar{\theta})$ as
  $k \to \infty$. Note that $\bar{\theta} = \theta_{\bar{u}}$ is not
  clear yet, but of course we will show exactly this now. By~\cite[Lem.~5.5]{hajo15}, the equations
  \[\partial_t \zeta + \AA(\theta_{n_k})\zeta =
  \Psi_{u_{n_k}}(\theta_{n_k}) + \alpha\theta_l, \quad
  \theta_{n_k}(T_0) = \theta_0\] have solutions $\zeta_{n_k} \in
  W^{1,s}(J;W^{-1,q}_\emptyset)\cap L^{s}(J;W^{1,q})$, which, due
  to uniqueness of solutions, must coincide with $\theta_{n_k}$. This
  means, on the one hand, that $\zeta_{n_k} = \theta_{n_k} \to
  \bar{\theta}$ in $C(\olQ)$ as $k \to \infty$. On the other hand,~\cite[Lem.~5.5]{hajo15} also shows that the sequence
  $(\zeta_{n_k})$ has a limit $\bar{\zeta}$ in the maximal regularity
  space as $k$ goes to infinity, where $\bar{\zeta}$ is the solution
  of the limiting problem
  \[\partial_t \zeta + \AA(\bar{\theta})\zeta =
  \Psi_{\bar{u}}(\bar{\theta}) + \alpha\theta_l, \quad \zeta(T_0) =
  \theta_0.\] We do, however, already know that $\bar{\zeta} =
  \bar{\theta}$, such that $\bar{\theta}$ is the unique global
  solution to the nonlinear problem for the limiting control
  $\bar{u}$, i.e., $\bar{\zeta} = \bar{\theta} =:
  \theta_{\bar{u}}$. Hence, $\bar{u}$ is still a global control.
  \hfill\end{proofof}

\begin{rem}
 Note that we used Proposition~\ref{p-divproperty}~\ref{p-divproperty-groeger} 
 instead of Lemma~\ref{l-multiplierkappa} at the beginning of the proof of Theorem~\ref{t-globalcompact}.
 This is indeed a crucial point, since Proposition~\ref{p-divproperty}~\ref{p-divproperty-groeger} 
 implies the isomorphism property and a \emph{uniform} bound of the inverse for all coefficient functions 
 that share the same ellipticity constant and the same $L^\infty$-bound. Thus, in our concrete 
 situation, the norm of
 $(-\nabla\cdot\sigma(\theta)\varepsilon\nabla)^{-1}$ is completely
 determined by $\Omega \cup \GD$ and the data from 
 Assumption~\ref{a-assu1}~\ref{a-assu1-sigmaeta} and~\ref{a-assu1}~\ref{a-assu1-eps}, 
 which gives the estimate in~\eqref{eq:invnormest}. 
 By contrast, the application of Lemma~\ref{l-multiplierkappa} would require to control 
 the norm of $\sigma(\theta)$ in $C(\overline{\Omega})$, see also Theorem~\ref{t-contoncont}. 
 This however cannot be guaranteed a priori so that Proposition~\ref{p-divproperty}~\ref{p-divproperty-groeger} 
 is indeed essential for the proof of Theorem~\ref{t-globalcompact}. 
 Since the integrability exponent from Proposition~\ref{p-divproperty}~\ref{p-divproperty-groeger} 
 is in general less than 3 and therefore less than $q$, one needs an improved 
 regularity in time to have the continuous embedding in the desired H\"older space, 
 cf.\ Corollary~\ref{c-intpolcompact}. Therefore it is \emph{not} sufficient to require 
 $s > r^*(q)$ and the more restrictive condition $s> \bar r(q,\varsigma)$ is imposed instead.
\end{rem}

%
Next, we incorporate the control- and state constraints in~\eqref{P} into the
control problem. For this purpose, let us introduce the set
\begin{equation}\label{eq:uad}
 \UU^\ad := \{ u \in L^2(J;L^2(\GN)) \colon 0 \leq u \leq u_{\max} \text{ a.e.\ in } \Si_N\}.
\end{equation}

\begin{definition}
  We call a global control $u \in \UU_g$ \emph{feasible}, if $u \in
  \UU^\ad$ and the
  associated state satisfies $\s(u)(x,t)
  \leq \theta_{\max}(x,t)$ for all $(x,t) \in \olQ$.
\end{definition}

While the state constraints give upper bounds on the values of
feasible solutions, lower bounds are natural in the problem and implicitly contained
in~\eqref{e-strongparabolic}--\eqref{e-strongellipticdirichlet} in the
sense that the temperature of the workpiece associated with $\O$ will
not drop below the minima of the surrounding temperature (represented
by $\theta_l$) and the initial temperature distribution $\theta_0$.

\begin{proposition}\label{p-lowerbounds} For every solution
  $(\theta,\varphi)$ in the sense of Theorem~\ref{t-locexist} with
  maximal existence interval $J_{\max}$, we have $\theta(x,t) \geq
  m_{\inf} := \min(\essinf_\Si \theta_l,\min_{\olom} \theta_0)$ for
  all $(x,t) \in \olom \times [T_0,T_\bullet]$, where $T_\bullet
  \in J_{\max}$. 
\end{proposition}

See Proposition~\ref{p-lowerboundsappendix} in the Appendix for a proof. Analogously, we find that $u \equiv 0$ is a feasible control
under Assumption~\ref{assu-p}~\ref{assu-p-thetamax}, the latter demanding that
the surrounding temperature and the initial temperature do not exceed
the state bounds at any point.

\begin{corollary}\label{c-0feasible}
  The control $u \equiv 0$ is a feasible one.
\end{corollary}

\begin{proof}
  By Corollary~\ref{c-smallball}, $u \equiv 0$ is a global control, it
  obviously satifies the control constraints, and using the same
  reasoning as in Proposition~\ref{p-lowerboundsappendix} with
  Assumption~\ref{assu-p}~\ref{assu-p-thetamax}, we obtain
  $\theta_{u\equiv 0} \leq \theta_{\max}$. \hfill
\end{proof}

Let us next introduce a modified control space, fitting the norm in
the objective functional in~\eqref{P}. So far, the controls originated from the space $L^{2s}(J;W^{-1,q}_\GD)$ 
with $s > \bar r(q, \varsigma)$. For the optimization, we now switch to the more
advanced control space
\begin{equation}\label{eq:ctrlspace}
 \U := W^{1,2}(J;L^2(\GN)) \cap L^p(J;L^p(\GN))
\end{equation}
with the standard norm $\|u\|_\U = \|u\|_{W^{1,2}(J;L^2(\GN))} +
\|u\|_{L^p(J;L^p(\GN))}$. Since $p > \frac{4}{3} q - 2$ by Assumption~\ref{assu-p},
this space continuously embeds into $L^{2s}(J;W^{-1,q}_\GD)$, which will give the 
boundedness required for Theorem~\ref{t-globalcompact}. Moreover, this embedding is even \emph{compact}, as the following result shows:

\begin{proposition}
  \label{p-compactcontrolspace} Let $p > 2$. The space $\U$ is embedded into a
  H\"older space $C^\varrho(\olJ;L^{\mathfrak p}(\GN))$ for
  some $\varrho > 0$ and $2<\mathfrak p < \frac{p+2}{2}$. In particular, there exists a \emph{compact}
  embedding $\EE \colon \U \embeds L^s(J;W^{-1,q}_\GD)$ for every $p > \frac{4}{3}q - 2$ and $s \in [1,\infty]$.
\end{proposition}

\begin{proof} From the construction of real interpolation spaces by means of
the trace method it immediately follows that \[\U \embeds
C(\olJ;(L^p(\GN),L^2(\GN)_{\frac{2}{p+2},\frac{p+2}{2}})) =
C(\olJ;L^{\frac{p+2}{2}}(\GN)),\]
see~\cite[Ch.~1.8.1--1.8.3 and Ch.~1.18.4]{T78}. With similar
reasoning as for~\eqref{e-maxreghoelderembed}, see
also~\cite[Lem.~3.17]{h/m/r/r} and its proof, we also may show $\U
\embeds C^\varrho(\olJ;(L^p(\GN),L^2(\GN))_{\tau,1})$ for all $\tau \in
(\frac{2}{2+p},1)$ and some $\varrho = \varrho(\tau) > 0$. Moreover, \[\left(L^p(\GN),L^2(\GN)\right)_{\tau,1} \embeds
\left[L^p(\GN),L^2(\GN)\right]_{\tau} = L^{\mathfrak p}(\GN)\]
with $\mathfrak p = \mathfrak p(\tau)
= (\frac{1-\tau}{p} +\frac{\tau}{2})^{-1} \in (2,\frac{2+p}{2})$ for
$\tau \in (\frac{2}{2+p},1)$, see~\cite[Ch.~1.10.1/3 and
Ch.~1.18.4]{T78}. This means we have $\U \embeds
C^\varrho(\olJ;L^{\mathfrak p}(\GN))$ for all $\mathfrak p \in
(2,\frac{2+p}{2})$, with $\varrho > 0$ depending on $\mathfrak p$. If
$\mathfrak p > \frac{2}{3} q$, then there is an embedding $L^{\mathfrak p}(\GN)) \embeds W^{-1,q}_\GD$, cf.\
Remark~\ref{r-einbett/ident}, and this is even compact in this
case as we will show below. To make $\mathfrak p >
\frac{3}{2}q$ possible, we need $\frac{p+2}{2} >
\frac{2}{3}q$, which is equivalent to $p > \frac{4}{3}q - 2$. Now the
vector-valued Arzel\`{a}-Ascoli Theorem yields the assertion. 

It remains to show that $L^{\mathfrak p}(\GN) \embeds W^{-1,q}_\GD$
compactly for $\mathfrak p > \frac{2}{3}q$, or equivalently
$W^{1,q'}_\GD \embeds L^{\mathfrak p'}(\GN)$
compactly. From~\cite[Ch.~1.4.7, Cor.~2]{mazya} and~\cite[Lem.~3.2]{hal/rehcoer} we
obtain \[\|u\|_{L^{\mathfrak p'}(\partial \Omega)} \leq C
  \|u\|_{W^{1,q'}}^\tau \|u\|_{L^{q'}}^{1-\tau} \quad
  \text{for all } u \in W^{1,q'}\] for $\mathfrak p' \in
  (\frac{2}{3}q',\frac{2q'}{3-q'})$ and $\tau =
  \frac{3}{q'}-\frac{2}{\mathfrak p'}$. Note that $\tau \in (0,1)$ for
  the given range of $\mathfrak p'$.  The preceding inequality implies
  $(L^{q'},W^{1,q'})_{\tau,1} \embeds L^{\mathfrak p'}(\partial
  \Omega)$, cf.~\cite[Lem.~1.10.1]{T78} and hence, due to the compact
  embedding $W^{1,q'} \embeds L^{q'}$ as of~\cite[Ch.~1.4.6, Thm.~2]{mazya}, $W^{1,q'} \embeds
  L^{\mathfrak p'}(\partial \Omega)$ compactly for all $\mathfrak p'
  \in (0,\frac{2q'}{3-q'})$ by~\cite[Ch.~1.16.4]{T78}. With
  $W^{1,q'}_\GD \embeds W^{1,q'}$ and $L^{\mathfrak p'}(\partial
  \Omega) \embeds L^{\mathfrak p'}(\GN)$, this means $W^{1,q'}_\GD \embeds
  L^{\mathfrak p'}(\GN)$ compactly for $\mathfrak p > \frac{2}{3}q$.\hfill
\end{proof}

\begin{definition}\label{def:reduced}
  Consider the embedding $\EE$ from Proposition~\ref{p-compactcontrolspace} with range in 
  $L^{2s}(J;W^{-1,q}_{\Gamma_D})$, where $s > \bar r(q,\varsigma)$ is the integrability exponent 
  from the objective functional. We set 
  $$\U_g := \{ u \in \U\colon \EE(u) \in \UU_g\}$$
  and define the mapping  
  \[ \s_\EE := \s \circ \EE \colon \U_g \to W^{1,s}(J;W^{1,q})\cap L^{s}(J;W^{-1,q}_\emptyset).\]
  Moreover, we define the \emph{reduced objective
    func\-tional} $j$ obtained by reducing the objective functional
  in~\eqref{P} to $u$, i.e.,
  \[j(u) = \frac{1}{2}\int_E \left|\s_\EE(u)(T_1)-\theta_d\right|^2 \,
  \dd x + \frac{\gamma}{s}\|\nabla\s_\EE(u)\|_{L^{s}(J;L^q)}^{s} +
  \frac{\beta}{2} \int_{\Sigma_N} (\partial_tu)^2 + |u|^p \,
  \dd\omega\,\dd t,\] as a function on $\U_g$. Further, let $\U^\ad := \U
  \cap \UU^\ad$ and $\U_g^\ad :=
  \U_g \cap \UU^\ad$, where $\UU^\ad$ is as defined in~\eqref{eq:uad}.
\end{definition}

One readily observes that $\s_\EE$ on $\U_g$ is still continuously
differentiable with the derivative $h \mapsto \s'_\EE(u)h = \s'(\EE u)\EE h$.

\begin{theorem}\label{t-existenceoptcont}
  There exists an optimal solution $\bar{u} \in \U_g^\ad$ to the
  problem
  \begin{equation}\min_{u \in \U_g^\ad} j(u) \quad \text{such that}
    \quad
    \s_\EE(u)(x,t) \leq \theta_{\max}(x,t) \quad \forall (x,t) \in \olQ.\label{e-reducedproblem}\tag{P$_u$}\end{equation}
\end{theorem}

\begin{proof}
  Thanks to the existence of the feasible control $u\equiv 0$, cf.\
  Corollary~\ref{c-0feasible}, the objective functional is bounded from below by $0$. Thus there exists a minimizing sequence
  of feasible controls $(u_n)$ in $\U_g^\ad$ such that $j(u_n) \to
  \inf_{u \in \U_g^\ad} j(u)$ in $\R$. On account of
  \begin{equation}\label{eq:radunbound}
    \int_{\Sigma_N} (\partial_t u)^2 + |u|^p \, \dd\omega\,\dd t
    \longrightarrow \infty \quad \text{when} \quad \|u\|_{\U} \longrightarrow \infty,
  \end{equation}    
  the objective functional is radially unbounded so that
  the minimizing sequence is
  bounded in $\U$ and, due to reflexivity of $\U$, has a weakly convergent subsequence (again
  $(u_n)$), converging weakly to some $\bar{u} \in \U$. As
  $\U^\ad$ is closed and convex, we have $\bar u \in \U^\ad$. By the compact
  embedding from Proposition~\ref{p-compactcontrolspace}, $(u_n)$
  converges strongly in $L^{2s}(J;W^{-1,q}_{\GD})$, also to $\bar{u}
  \in L^{2s}(J;W^{-1,q}_{\GD})$. The fact that state
  constraints are present and Proposition~\ref{p-lowerbounds} imply
  that the family $(\theta_{u_n})$ is uniformly bounded in time and
  space for every feasible control $u$. Together with the gradient
  term in the objective functional, Theorem~\ref{t-globalcompact} now
  shows $\bar{u} \in \U_g$, hence $\bar u \in \U_g^\ad$. Moreover, $\s_\EE(u_n) \to \s_\EE(\bar{u})$
  in $W^{1,s}(J;W^{1,q})\cap L^{s}(J;W^{-1,q}_\emptyset)$, which
  immediately implies convergence of the first two terms in the
  objective functional (each as $n$ goes to infinity). The third term,
  corresponding to $\U$, is clearly continuous and convex on $\U$ and
  as such weakly lower semicontinuous, hence we find
  \[\inf_{u\in\U_g^\ad} j(u) = \lim_{n\to\infty} j(u_n) \geq
  j(\bar{u})\] and thus $j(\bar{u}) = \inf_{u \in \U_g^\ad} j(u)$.
  \hfill\end{proof}


\begin{rem}
  In the proof of Theorem~\ref{t-existenceoptcont}, boundedness of
  minimizing sequence $(u_n)$ in the control space $\U$ was essential
  and followed from the radial unboundedness of the objective
  functional as seen in~\eqref{eq:radunbound}.  Alternatively, one could also
  assume that the upper bound $u_{\max}$ in the control constraints
  satisfies $u_{\max} \in L^p(J;L^p(\Gamma_N))$ with $p > \frac{4}{3}q
  - 2$. In this case, an objective functional of the form
  \begin{equation*}
    \frac{1}{2}\|\theta(T_1)-\theta_d\|_{L^2(E)}^2 +
    \frac{\gamma}{s} \|\nabla\theta\|_{L^{s}(T_0,T_1;L^q(\Omega))}^{s} +
    \frac{\beta}{2} \int_{\Si_N} (\partial_t u)^2 \,\dd\omega \,\dd t
  \end{equation*}
  is sufficient to establish the existence of a globally optimal
  control.
\end{rem}

\subsection{Necessary optimality conditions}\label{subs-noc}

This section is devoted to the derivation of necessary
  optimality conditions for~\eqref{P} in the
  form~\eqref{e-reducedproblem}. To this end, let us start with the
definition of the Lagrangian function. It is well-known that the
Lagrangian multipliers associated to the state constraints may, in
general, only be regular Borel measures, see for
instance~\cite{C93}. Hence, we introduce the space $\MM(\olQ)$ as the
space of regular Borel measures on $\olQ$ and, simultaneously, as the
dual space of $C(\olQ)$.

\begin{definition}\label{d-lagrange}
  The \emph{Lagrangian function} $\fL \colon \U_g \times \MM(\olQ) \to \R$
  associated with~\eqref{e-reducedproblem} is given by
  \[\fL(u,\mu) = j(u) + \langle
  \mu,\s_\EE(u)-\theta_{\max}\rangle_{\MM(\olQ),C(\olQ)},\] where $j$
  is the reduced objective functional.
\end{definition}

\begin{definition}\label{d-qlaplace}
  We denote by $\Delta_q \colon W^{1,q} \to W^{-1,q'}_\emptyset$ the
  \emph{(weak) $q$-Laplacian}, given by
  \[\langle \Delta_q \psi,\xi\rangle :=
  \int_\O |\nabla \psi|^{q-2}\nabla\psi\cdot\nabla\xi \, \dd x\] for
  each $\psi,\xi \in W^{1,q}$.
\end{definition}

The chain rule immediately yields the derivative of $\fL$ with respect
to $u$:

\begin{lemma}\label{l-lagrangianderivative}
  The Lagrangian function $\fL$ is continuously differentiable with
  respect to $u$. Abbreviating the states by $\theta_u := \s_\EE(u)$
  and $\theta'_u = \s_\EE'(u)h$, the partial derivative in direction
  $h \in \U$ is given by
  \begin{equation}
    \begin{split}\partial_u \fL(u,\mu)h & = \int_E (\theta_u(T_1) -
      \theta_d)\theta_u'(T_1) \, \dd x + \gamma\int_{T_0}^{T_1}
      \|\nabla \theta_u(t)\|_{L^q}^{s-q}\langle
      \Delta_q\theta_u(t),\theta'_u(t)\rangle\, \dd t \\ & \qquad +
      \beta\int_{\Si_N} \partial_t u\partial_t h +
      \frac{p}{2}|u|^{p-2}uh\, \dd\omega\,\dd t) + \langle
      \mu,\theta'_u\rangle_{\MM(\olQ),C(\olQ)}
    \end{split}\label{e-lagrangederivative}
  \end{equation}
  with $\Delta_q$ given as in Definition~\ref{d-qlaplace}.
\end{lemma}

Using the Lagrangian function and its derivative, we characterize
local optima of~\eqref{e-reducedproblem}. We say that that a feasible
control $\bar{u}$ is \emph{locally optimal} if there exists an
$\epsilon > 0$ such that $j(\bar{u}) \leq j(u)$ for all feasible $u
\in \U_g^\ad$ with $\|u-\bar{u}\|_{\U} < \epsilon$. As we will see in
the proof of Theorem~\ref{t-lagrangeexistence}, the restriction to
global controls $u \in \U_g$ does not influence the derivation of
optimality conditions, since $\U_g$ is an \emph{open} set by
Theorem~\ref{t-open}.

\begin{definition}\label{d-lagrangemultiplier} A measure $\bar \mu \in \MM(\olQ)$
  is called a \emph{Lagrangian multiplier} associated with the state
  constraint in~\eqref{e-reducedproblem}, if for a locally optimal
  control $\bar{u}$ the \emph{KKT conditions}
  \begin{align}
    \label{e-lagrangemult1}
    \bar \mu & \geq 0, \\
    \label{e-lagrangemult2}
    \langle \bar \mu, \s_\EE(\bar{u}) - \theta_{\max}\rangle_{C(\olQ)} & = 0, \\
    \label{e-lagrangemult3}
    \left\langle\partial_u \fL(\bar{u},\bar \mu),u-\bar
      u\right\rangle_\U & \geq 0 \quad \forall u \in \U^\ad
  \end{align}
  hold true. Here,~\eqref{e-lagrangemult1} means that $\langle \bar
  \mu, f\rangle_{C(\olQ)} \geq 0$ for all $f \in C(\olQ)$ with $f(x,t)
  \geq 0$ for all $(x,t) \in Q$. Note
    that~\eqref{e-lagrangemult3} has to be satisfied for \emph{all} $u\in
    \U^\ad$ instead of only in $\U^\ad_g$, the latter being defined in
    Definition~\ref{def:reduced}.
\end{definition}

It is well-known that, in general, a so-called regularity condition is
needed in order to ensure the existence of a Lagrangian multiplier. In
this case, we rely on the linearized Slater condition, which is a
special form of Robinson's regularity condition.

\begin{theorem}\label{t-lagrangeexistence}
  Let $\bar{u}$ be a locally optimal control and let the following
  so-called \emph{linearized Slater condition} be satisfied: There
  exists $\hat u \in \U_g^\ad$ such that there is a $\delta >
  0$ with the property
  \begin{equation}\s_\EE(\bar{u})(x,t) +
    \s_\EE'(\bar{u})(\hat u-\bar{u})(x,t) \leq
    \theta_{\max}(x,t) - \delta \quad \text{for all } (x,t) \in Q. \label{e-linearizedslater}\end{equation}
  Then there exists a Lagrangian multiplier $\bar\mu \in \MM(\olQ)$
  associated with the state constraint in~\eqref{e-reducedproblem},
  i.e., such that~\eqref{e-lagrangemult1}-\eqref{e-lagrangemult3} is
  satisfied.
\end{theorem}

  \begin{proof}
    Since $\U \embeds L^{2s}(J;W^{-1,q}_{\Gamma_D})$ as seen in
    Proposition~\ref{p-compactcontrolspace}, Theorem~\ref{t-open}
    implies that there is an open ball $B_\delta(\bar u) \subset \U$
    around $\bar u$ with radius $\delta > 0$ such that $B_\delta(\bar
    u) \cap \U^\ad \subset \U^\ad_g$.  We consider the auxiliary
    problem
    \begin{equation}\tag{P$_\textup{aux}$}\label{eq:paux}
      \left.
        \begin{aligned}
          \min & \quad j(u)\\
          \text{s.t.} & \quad u \in B_\delta(\bar u) \cap \U^\ad,
          \quad \s_\EE(u)(x,t) \leq \theta_{\max}(x,t) \quad \forall
          (x,t) \in \olQ.
        \end{aligned}
        \qquad \right\}
    \end{equation}
    Clearly, $\bar u$ is also a local minimizer of this
    problem. Moreover, in contrast to $\U^\ad_g$ appearing
    in~\eqref{e-reducedproblem}, the feasible set $B_\delta(\bar u)
    \cap \U^\ad$ is now \emph{convex}. Therefore the standard
    Karush-Kuhn-Tucker (KKT) theory in function space can be applied
    to~\eqref{eq:paux}, see e.g.~\cite[Thm.~3.1]{ZK79},~\cite[Thm.~5.2]{C93}
    or~\cite[Thm.~3.9]{BS00}.  Hence, on account of the linearized
    Slater condition in~\eqref{e-linearizedslater}, there exists a
    Lagrange multiplier $\bar\mu \in \MM(\overline{Q})$ so
    that~\eqref{e-lagrangemult1},~\eqref{e-lagrangemult2}, and
    \begin{equation}\label{eq:gradeq}
      \left\langle\partial_u \fL(\bar{u},\bar \mu),v-\bar u\right\rangle_\U \geq 0 \quad \forall v \in B_\delta(\bar u) \cap \U^\ad
    \end{equation}
    are fulfilled. Now, let $u\in \U^{ad}$ be arbitrary. Then, due to
    convexity, $\bar u + \tau (u -\bar u) \in B_\delta(\bar u) \cap
    \U^\ad$ for $\tau > 0$ sufficiently small such that this function
    can be chosen as test function in~\eqref{eq:gradeq}, giving in
    turn~\eqref{e-lagrangemult3}. \hfill
  \end{proof}

Let us now transform~\eqref{e-lagrangemult1}-\eqref{e-lagrangemult3}
into an optimality system involving an adjoint state. To this
end, we aim to reformulate the derivative expression for $\partial_u
\fL(\bar{u},\mu)$ from Lemma~\ref{l-lagrangianderivative} in a
designated locally optimal point $\bar{u}$. For brevity, we define
\[\X = W^{1,s}(J;W^{-1,q}_\emptyset) \cap L^{s}(J;W^{1,q}) \quad
\text{and} \quad X_{s} =
(W^{1,q},W^{-1,q}_\emptyset)_{\frac{1}{s},s}.\] The plan is to use the
adjoint of the derivative of the control-to-state operator. We will
show that $\s_\EE'(\bar{u})^*$ is associated to the solution operator
(in an appropriate sense) to the \emph{adjoint system}, which we
formally introduce as follows:

\begin{definition}\label{d-adjointsys} For given, fixed functions $\theta$
  and $\varphi$, given terminal value $\vartheta_T$ and
  inhomogeneities $f_1,f_2,g_1,g_2$, we call the following system the
  \emph{adjoint system}:
  \begin{equation}\label{e-adjoint}
    \left.\begin{aligned}
        -\partial_t\vartheta - \dive(\eta(\theta)\kappa\nabla\vartheta) &
        = (\sigma'(\theta)\vartheta\varepsilon\nabla\varphi)\cdot\nabla\varphi -
        (\sigma'(\theta)\varepsilon\nabla\varphi)\cdot\nabla\psi  \\
        & \qquad - (\eta'(\theta)\kappa\nabla\vartheta)\cdot\nabla\theta
        + f_1 && \emph{in }
        Q,
        \\ 
        \nu \cdot \eta(\theta)\kappa\nabla\vartheta + \alpha\vartheta &
        = f_2 && \emph{on } \Si, \\
        \vartheta(T_1) & = \vartheta_T && \emph{in } \O, \\[1em]
        -\dive(\sigma(\theta)\varepsilon\nabla\psi) & =
        -2\dive(\sigma(\theta)\vartheta\varepsilon\nabla\varphi) + g_1 && \emph{in }
        Q, \\
        \nu \cdot \nabla\sigma(\theta)\varepsilon\nabla\psi & = 2\nu \cdot
        \sigma(\theta)\vartheta\varepsilon\nabla\varphi + g_2 && \emph{on }
        \Si_N, \\
        \psi & = 0 && \emph{on } \Si_D.
      \end{aligned} \right\}
  \end{equation}
\end{definition}

More specified assumptions about the inhomogeneities $f_1,f_2,g_1,g_2$
and the terminal value $\vartheta_T$ will be given in the
following. Note that~\eqref{e-adjoint} is only a \emph{formal}
representation of the adjoint of the linearized system
of~\eqref{e-strongparabolic}-\eqref{e-strongellipticdirichlet}. We
will work with the abstract version, referring to~\eqref{e-paraex}
and~\eqref{e-elliptex} and its linearizations,
cf.~\eqref{e-linearized} or Remark~\ref{r-linearizedtwoequations}.

\begin{definition}\label{d-abstractadjointsys} Let $\theta \in
  \X$ be fixed and set $\varphi = \JJ(\sigma(\theta))u$. Further, let
  $f \in \X'$, $\vartheta_T \in X_r'$, and $g \in L^{
    (2s)'}(J;W^{-1,q'}_\GD)$ be given with $(2s)' = \frac{2s}{2s - 1}$. The
  \emph{abstract adjoint system} is given by
  \begin{multline}\label{e-abstractadjoint}-\partial_t \vartheta
    + \partial_\theta \AA(\theta)\vartheta = -
    (\eta'(\theta)\kappa\nabla\vartheta)\cdot\nabla \theta +
    (\sigma'(\theta)\vartheta\varepsilon\nabla
    \varphi)\cdot\nabla\varphi + \delta_{T_1}\otimes\vartheta_T -
    \delta_{T_0}\otimes\chi +f \\ -
    (\sigma'(\theta)\varepsilon\nabla\varphi)\cdot\nabla\left[\JJ(\sigma(\theta))^*(-2\nabla\cdot\sigma(\theta)\vartheta\varepsilon\nabla\varphi
      + g)\right].
  \end{multline}
  Here, $\delta_{T_0}$ and $\delta_{T_1}$ are Dirac measures in $T_0$
  and $T_1$, obtained as the adjoints of the point evaluations in
  $T_0$ and $T_1$, respectively. The latter are continuous mappings
  from $C(\olJ;X_s)$ to $X_s$, such that
  $\delta_{T_0}\otimes\vartheta_T$ and $\delta_{T_1} \otimes \chi$ are
  seen as objects from $\MM(J;X_s')$. We say that the functions
  $(\vartheta,\chi) \in L^{s'}(J;W^{1,q'}) \times X_s'$ are a
  \emph{weak solution} of~\eqref{e-abstractadjoint}
  or~\eqref{e-adjoint}, if
  \begin{equation}
    \begin{split}
      \int_J \langle\partial_t \xi,\vartheta\rangle_{W^{1,q'}} \, \dd
      t & = -\int_J\int_\O \langle
      (\eta(\theta)\kappa\nabla\vartheta)\nabla\xi\, \dd x \dd t -
      \int_J\int_\Gamma\alpha\vartheta\xi \, \dd \omega \dd t \\ &
      \quad - \int_J \int_\O
      \left[(\eta'(\theta)\kappa\nabla\vartheta)\cdot\nabla \theta -
        (\sigma'(\theta)\vartheta\varepsilon\nabla
        \varphi)\cdot\nabla\varphi\right]\xi \, \dd x \dd t \\ & \quad
      - \int_J \int_\O
      (\sigma'(\theta)\varepsilon\nabla\varphi)\cdot\nabla\left[\JJ(\sigma(\theta))^*(-2\nabla\cdot\sigma(\theta)\vartheta\varepsilon\nabla\varphi
        + g)\right]\xi \, \dd x \dd t \\ & \quad + \langle
      \vartheta_T,\xi(T_1)\rangle_{X_r} - \langle
      \chi,\xi(T_0)\rangle_{X_r} + \langle f,\xi\rangle_\X
    \end{split}\label{e-weakadjoint}
  \end{equation}
  is true for all $\xi \in
  \X$. Equivalently,~\eqref{e-abstractadjoint} holds true in $\X'$.
\end{definition}

Note that the functionals $\delta_{T_0}\times\chi$ and
$\delta_{T_1}\otimes\vartheta_T$ are well-defined in $\X'$ due to $\X
\embeds C(\olJ;X_s)$. Of course, the inhomogeneities $f_1,f_2$ and
$g_1,g_2$ from ~\eqref{e-adjoint} are represented by $f = f_1 + f_2$
and $g=g_1+g_2$, respectively. Moreover, thanks to the symmetry of
$\varepsilon$, one easily sees that $\JJ(\sigma(\theta))^*$ is
formally selfadjoint, which is the basis of the following

\begin{rem}\label{r-adjointtwoequations}
  Similarly to Remark~\ref{r-linearizedtwoequations}, we introduce
  \[\psi(\vartheta) :=
  \JJ(\sigma(\theta))^*(-2\nabla\cdot\sigma(\theta)\vartheta\varepsilon\nabla\varphi
  + g),\] which allows to split~\eqref{e-abstractadjoint} back into
  two equations, namely
  \begin{align*}
    -\partial_t \vartheta + \partial_\theta \AA(\theta)\vartheta & =
    (\sigma'(\theta)\vartheta\varepsilon\nabla
    \varphi)\cdot\nabla\varphi -
    (\eta'(\theta)\kappa\nabla\vartheta)\cdot\nabla \theta -
    (\sigma'(\theta)\varepsilon\nabla\varphi)\cdot\nabla\psi \\ &
    \qquad + \delta_{T_1}\otimes\vartheta_T - \delta_{T_0}\otimes\chi
    + f, \\
    -\nabla\cdot \sigma(\theta)\varepsilon\nabla\psi & =
    -2\nabla\cdot\sigma(\theta)\vartheta\varepsilon\nabla\varphi + g,
  \end{align*}
  to be understood as in~\eqref{e-weakadjoint}. This is exactly a very
  weak abstract formulation of the formal adjoint
  system~\eqref{e-adjoint} with inhomogeneities $f = f_1+f_2$ and
  $g=g_1+g_2$ and terminal value $\vartheta_T$. Note that the first
  equation is supposed to hold in $\X'$, the second one in
  $L^{(2s)'}(J;W^{-1,q'}_\GD)$.
\end{rem}

We next show that the abstract adjoint~\eqref{e-abstractadjoint}
always admits a unique weak solution for $f \in \X'$ and $g \in
L^{(2s)'}(J;W^{-1,q'}_\GD)$. This will follow directly from
Theorem~\ref{t-open} using an adjoint-approach (see
e.g.~\cite[Ch.~7]{A05}). Since the inhomogeneity $f$
in~\eqref{e-abstractadjoint} will in fact contain the Lagrange
multiplier $\mu$ introduced in Definition~\ref{d-lagrangemultiplier},
we will not investigate the adjoint system more specifically under
additional regularity assumptions on $f$, since the Lagrange
  multipliers are in general only measures and thus limit said
  regularity in a crucial way anyhow. In particular, this lack
  of regularity is the very obstacle which permits time-derivatives
  for weak solutions to~\eqref{e-abstractadjoint},
  cf.~\cite[Prop.~6.1]{A05}.  Nevertheless, even in the absence of
  measure-valued Lagrange multipliers, the time regularity of
  the adjoint state is still limited by the differential operator
itself, since $(\eta'(\theta)\kappa\nabla\vartheta)\cdot\nabla\theta$
is only integrable in time (as opposed to $s'$-integrable).

\begin{theorem}
  \label{t-weakadjointexists} For every terminal value $\vartheta_T
  \in X_s' = (W^{1,q'},W^{-1,q'}_\emptyset)_{\frac{1}{s'},s'}$ and all
  imhomogeneities $f \in \X'$ and $g \in L^{(2s)'}(J;W^{-1,q'}_\GD)$,
  there exists a unique weak solution $(\vartheta,\chi) \in
  L^{s'}(J;W^{1,q'}) \times X_{s}'$ of~\eqref{e-abstractadjoint} in
  the sense of Definition~\ref{d-abstractadjointsys}.
\end{theorem}

\begin{proof} The equality $X_s' =
  (W^{1,q'},W^{-1,q'}_\emptyset)_{\frac{1}{s'},s'}$ follows from the
  usual duality properties of interpolation functors, see~\cite[Ch.\
  1.11.2 and 1.3.3]{T78}.  Recall the operator
  \[\BB \colon \X \times
  L^{2s}(J;W^{-1,q}_{\GD}) \to L^s(J;W^{-1,q}_\emptyset) \times
  (W^{1,q},W^{-1,q}_\emptyset)_{\frac{1}{s},s},\] from
  Theorem~\ref{t-open} with $r = s > \bar r(q, \varsigma) \geq
  r^*(q)$.
  The partial derivative w.r.t.\ $\theta$ of $\BB$ was given by
  \[\partial_\theta\BB(\theta,u)\xi = \left(\partial_t \xi +
    \AA(\theta)\xi - \nabla \cdot \eta'(\theta)\xi\kappa\nabla\theta
    - \partial_\theta \Psi_u(\theta)\xi,\xi(T_0)\right)\] with
  \[\partial_\theta \Psi_u(\theta)\xi =
  -2(\sigma(\theta)\varepsilon\nabla\varphi)\cdot\nabla\left[\JJ(\sigma(\theta)(-\nabla\cdot\sigma'(\theta)\xi\varepsilon\nabla\varphi)\right]+(\sigma'(\theta)\xi\varepsilon\nabla\varphi)\cdot\nabla\varphi,\]
  cf.~\eqref{e-Psiderivative}, and $\varphi =
  \JJ(\sigma(\theta))u$. Now, let $(\vartheta,\chi)$ be from
  $L^{s'}(J;W^{1,q'}) \times X_s'$. 
  We easily find
  \begin{equation}\left\langle - \nabla \cdot
      \eta'(\theta)\xi\kappa\nabla\theta,\vartheta\right\rangle_{W^{1,q'}}
    = \int_\O (\eta'(\theta)\xi\kappa\nabla\vartheta)\cdot\nabla\theta \, \dd
    x = \left\langle(\eta'(\theta)\kappa\nabla\vartheta)\cdot\nabla\theta,\xi\right\rangle_{W^{1,q}}\label{e-ellipadjoint}\end{equation}
  and 
  \begin{equation}\left\langle(\sigma'(\theta)\xi\varepsilon\nabla\varphi)\cdot\nabla\varphi,\vartheta\right\rangle_{W^{1,q'}}
    = \left\langle
      (\sigma'(\theta)\vartheta\varepsilon\nabla\varphi)\cdot\nabla\varphi,\xi\right\rangle_{W^{1,q}}.\label{e-righthandadjoint}\end{equation}
  Let us turn to the complicated term in
  $\partial_\theta\Psi_u(\theta)$. Analogously
  to~\eqref{e-ellipadjoint}, we find
  \begin{multline}\left\langle
      2(\sigma(\theta)\varepsilon\nabla\varphi)\cdot\nabla\left[\JJ(\sigma(\theta)(-\nabla\cdot\sigma'(\theta)\xi\varepsilon\nabla\varphi)\right],\vartheta\right\rangle_{W^{1,q'}}
    \\ = \left\langle
      \JJ(\sigma(\theta))(-\nabla\cdot\sigma'(\theta)\xi\varepsilon\nabla\varphi),-2\nabla\cdot\sigma(\theta)\vartheta\varepsilon\nabla\varphi\right\rangle_{W^{-1,q'}_\GD}
    \\ = \left\langle
      -\nabla\cdot\sigma'(\theta)\xi\varepsilon\nabla\varphi,\JJ(\sigma(\theta))^*\left(-2\nabla\cdot\sigma(\theta)\vartheta\varepsilon\nabla\varphi\right)\right\rangle_{W^{1,q'}_\GD}
    \\ = \left\langle
      \xi,\left(\sigma'(\theta)\varepsilon\nabla\varphi\right)\nabla\left[\JJ(\sigma(\theta))^*
        \left(-2\nabla\cdot\sigma(\theta)\vartheta\varepsilon\nabla\varphi\right)\right]\right\rangle_{W^{-1,q'}_\emptyset}.\label{e-ekeladjoint}
  \end{multline}
  Symmetry of $\kappa$ implies that $\AA(\theta)$ is formally
  self-adjoint, i.e., $\AA(\theta)^*$ maps $W^{1,q'}$ into
  $W^{-1,q'}_\emptyset$, but is still given as in
  Definitions~\ref{d-quasilinear} and~\ref{d-opera0},
  respectively. Using this and
  equations~\eqref{e-ellipadjoint},~\eqref{e-righthandadjoint}
  and~\eqref{e-ekeladjoint}, we obtain
  \begin{align*}\left\langle\partial_\theta\BB(\theta,u)^*(\vartheta,\chi),\xi\right\rangle_\X
    & =
    \left\langle(\vartheta,\chi),\partial_\theta\BB(\theta,u)\xi\right\rangle_{L^s(J;W^{-1,q}_\emptyset)
      \times X_s} \\
    & = \int_J \left\langle \partial_t
      \xi,\vartheta\right\rangle_{W^{1,q'}} \, \dd t + \int_J
    \left\langle \AA^*(\theta)\vartheta, \xi\right\rangle_{W^{1,q}} \,
    \dd t \\ & \qquad + \int_J
    \left\langle(\eta'(\theta)\kappa\nabla\vartheta)\cdot\nabla\theta,\xi\right\rangle_{W^{1,q}}
    \dd t \\ & \qquad - \int_J \left\langle
      (\sigma'(\theta)\vartheta\varepsilon\nabla\varphi)\cdot\nabla\varphi,\xi\right\rangle_{W^{1,q}}
    \, \dd t + \left\langle \chi,\xi(T_0)\right\rangle_{X_s} \\ &
    \qquad + \int_J\! \left\langle
      \left(\sigma'(\theta)\varepsilon\nabla\varphi\right)\nabla\left[\JJ(\sigma(\theta))^*\!\left(-2\nabla\!\cdot\!\sigma(\theta)\vartheta\varepsilon\nabla\varphi\right)\right],\xi\right\rangle_{W^{1,q}}
    \, \dd t
  \end{align*}
  for all $\xi \in \X$. Moreover, in the proof of
  Theorem~\ref{t-open}, $\partial_\theta\BB(\theta,u)$ was found to be
  a topological isomorphism between the spaces $\X$ and
  $L^s(J;W^{-1,q}_\emptyset) \times
  X_s$ 
  and consequently $\partial_\theta\BB(\theta,u)^*$ is also a
  topological isomorphism between $L^{s'}(J;W^{1,q'}) \times X_s'$ and
  $\X'$. In particular, for every $\mathfrak{f} \in \X'$ there exists
  a unique $p = p_{\mathfrak{f}} \in L^{s'}(J;W^{1,q'}) \times X_s'$
  such that $\partial_\theta\BB(\theta,u)^*p = \mathfrak{f}$. Hence,
  setting \begin{equation}\bar{\mathfrak{f}} = f +
    \delta_{T_1}\otimes\vartheta_T -
    \left(\sigma'(\theta)\varepsilon\nabla\varphi\right)\nabla\left[\JJ(\sigma(\theta))^*g\right],\label{e-adjointrhs}\end{equation}
  the pair $(\bar{\vartheta},\bar{\chi}) := p_{\bar{\mathfrak{f}}}$
  satisfies~\eqref{e-weakadjoint} by the above form of
  $\partial_\theta\BB(\theta,u)^*$, and is exactly the searched-for
  unique solution as in
  Definition~\ref{d-abstractadjointsys}.\hfill\end{proof}

As hinted above, we immediately obtain the following characterization
of $\s'(u)^*$ for given $u \in \U_g$:

\begin{corollary}\label{c-adjointsolop}
  Let $(\vartheta,\chi)$ be the solution of~\eqref{e-weakadjoint} in
  the sense of Definition~\ref{d-abstractadjointsys} with
  inhomogeneites $f$ and $g$ and terminal value $\vartheta_T$. The
  adjoint linearized solution operator $\s_\EE'(u)^*$ then assigns to
  $f,g$ and $\vartheta_T$ in the form $\mathfrak f \in \X'$ as
  in~\eqref{e-adjointrhs} the functional $\EE^* \psi \in \U'$, where
  $\psi(\vartheta) \in L^{(2s)'}(J;W^{1,q'}_\GD)$ is given
  by \[\psi(\vartheta) =
  \JJ(\sigma(\theta_u))^*(-\nabla\cdot\sigma(\theta_u)\vartheta\varepsilon\nabla\varphi_u),\]
  similarly to Remark~\ref{r-adjointtwoequations}.
\end{corollary}

\begin{proof}
  In Theorem~\ref{t-open}, we found $\s'(u) =
  -[\partial_\theta\BB(\s(u),u)]^{-1}\partial_u\BB(\s'(u),u)$. Hence,
  with $\s_\EE'(u) = \s'(u) \circ \EE$, we obtain
  \[\s_\EE'(u)^*\mathfrak{f}
  =-\EE^*\partial_u\BB(\s_\EE(u),u)^*\partial_\theta\BB(\s_\EE(u),u)^{-*}\mathfrak{f}.\]
  In view of Theorem~\ref{t-weakadjointexists} and its proof,
  $\partial_\theta \BB(\s_\EE(u),u)^{-*}\mathfrak{f}$ is exactly the
  unique solution $(\vartheta,\chi)$ of~\eqref{e-weakadjoint} in the
  sense of Definition~\ref{d-abstractadjointsys} with inhomogeneites
  $f,g$ and terminal value $\vartheta_T$. Moreover, a repetition of
  the first lines of~\eqref{e-ekeladjoint} shows that \[-\partial_u
  \BB(\s_\EE(u),u)^*(\vartheta,\chi) =
  \JJ(\sigma(\theta_u))^*(-\nabla\cdot\sigma(\theta_u)\vartheta\varepsilon\nabla\varphi_u)
  = \psi(\vartheta).\] An application of $\EE^* :
  L^{(2s)'}(J;W^{1,q'}_\GD) \embeds \U'$ yields the claim.\hfill
\end{proof}

Having $\s_\EE'(u)^*$ at hand, we now proceed to establish the actual
necessary optimality conditions by manipulating the variational
inequality in the KKT conditions~\eqref{e-lagrangemult3}.

For a concise ``strong'' formulation in the following theorem, we
decompose measures $\mu \in \MM(\olQ)$ by restriction into $\mu =
\mu_{(T_0,T_1)} + \mu_{\{T_0\}\times\{T_1\}}$, with $\mu_{(T_0,T_1)}
\in \MM((T_0,T_1) \times \olom)$ and $\mu_{\{T_0\} \times\{T_1\}} \in
\MM((\{T_0\}\times\{T_1\})\times\olom)$. Both measures may in turn be
further decomposed into $\mu_{(T_0,T_1)} = \mu_\O + \mu_\Gamma$, where
$\mu_\O \in \MM((T_0,T_1) \times \O)$ and $\mu_\Gamma \in
\MM((T_0,T_1) \times \Gamma)$, and $\mu_{\{T_0\}\times\{T_1\}} =
\delta_{T_0} \otimes \mu_{T_0} + \delta_{T_1} \otimes \mu_T$ with
$\mu_{T_0},\mu_T \in \MM(\olom)$.
\begin{theorem}[First Order Necessary
  Conditions]\label{t-neccondfinal}
  Let $\bar{u} \in \U_g^\ad$ be a locally optimal control such that
  the linearized Slater condition~\eqref{e-linearizedslater} is
  satisfied. Let $\theta_{\bar{u}} = \s_\EE(\bar{u})$ be the state
  associated with $\bar{u}$ and set $\varphi_{\bar{u}} :=
  \varphi_{\bar{u}}(\theta_{\bar{u}})$. Then there exists a Lagrangian
  multiplier $\bar \mu \in \MM(\olQ)$ in the sense of
  Definition~\ref{d-lagrangemultiplier} and adjoint states $\vartheta
  \in L^{s'}(J;W^{1,q'})$ and $\psi \in L^{(2s)'}(J;W^{1,q'}_\GD)$,
  such that the formal system
  \begin{align*}
    -\partial_t\vartheta -
    \dive(\eta(\theta_{\bar{u}})\kappa\nabla\vartheta) & =
    (\sigma'(\theta_{\bar{u}})\vartheta\varepsilon\nabla\varphi_{\bar{u}})\cdot\nabla\varphi_{\bar{u}}
    -
    (\sigma'(\theta_{\bar{u}})\varepsilon\nabla\varphi_{\bar{u}})\cdot\nabla\psi  \\
    & \qquad - (\eta'(\theta)\kappa\nabla\vartheta)\cdot\nabla\theta +
    \|\nabla\theta_{\bar{u}}\|_{L^{s}(J;L^q)}^{s-q}\Delta_q\theta_{\bar{u}}
    + \bar\mu_\O && \emph{in } Q, \\
    \nu \cdot \eta(\theta_{\bar{u}})\kappa\nabla\vartheta +
    \alpha\vartheta &
    = \bar\mu_\Gamma && \emph{on } \Si, \\
    \vartheta(T_1) & = \chi_E(\theta_{\bar{u}}(T_1)-\theta_d) +
    \bar\mu_{T_1}&& \emph{in } \O, \\[1em]
    -\dive(\sigma(\theta_{\bar{u}})\varepsilon\nabla\psi) & =
    -2\dive(\sigma(\theta_{\bar{u}})\vartheta\varepsilon\nabla\varphi_{\bar{u}})
    && \emph{in }
    Q, \\
    \nu \cdot \sigma(\theta_{\bar{u}})\varepsilon\nabla\psi & = 2\nu
    \cdot
    \sigma(\theta_{\bar{u}})\vartheta\varepsilon\nabla\varphi_{\bar{u}}
    && \emph{on }
    \Si_N, \\
    \psi & = 0 && \emph{on } \Si_D.
  \end{align*}
  is satisfied in the sense of Definition~\ref{d-abstractadjointsys}
  and Remark~\ref{r-adjointtwoequations}. Moreover, $\bar{u}$ is the
  solution of the variational inequality
  \begin{equation}\label{eq:vi}
    \begin{aligned}
      \int_{\Si_N}\partial_t\bar u\, \partial_t (u-\bar u) +
      \frac{p}{2} |\bar u|^{p-2}(u-\bar u)
      + \frac{1}{\beta}(\tau_\GN\psi)(u-\bar u) \, \dd\omega\,\dd t \geq 0 \quad\qquad &\\[-1ex]
      \text{for all } u \in \U^\ad = \{u \in \U \colon 0 \leq u
        \leq u_{\max} \text{ a.e.\ in } \Sigma_N\}. &
    \end{aligned}
  \end{equation}
\end{theorem}

Note that the Lagrange multiplier $\bar\mu$ is not active on the set
$\{T_0\} \times \olom$ due to
Assumption~\ref{assu-p}~\ref{assu-p-thetamax} and the positivity and
complementary conditions~~\eqref{e-lagrangemult1}
and~\eqref{e-lagrangemult2}. Hence, $\bar\mu_{T_0}$ is zero and does
not contribute to the system of equations in
Theorem~\ref{t-neccondfinal}.  Note moreover that the
  variational inequality in~\eqref{eq:vi} is just a (semilinear)
  variational inequality of
  obstacle-type in time.

\begin{proof}
  Let $\bar{u}$ be a locally optimal control such that the linearized
  Slater condition~\eqref{e-linearizedslater} is
  satisfied. Theorem~\ref{t-lagrangeexistence} then yields the
  existence of a Lagrangian multiplier $\bar{\mu} \in \MM(\olQ)$ such
  that~\eqref{e-lagrangemult1}-\eqref{e-lagrangemult3} hold true. We
  show that these lead to the assertions.

  First consider the linear continuous functional
  \[\langle\chi_E(\theta_{\bar{u}}(T_1)-\theta_d),\Theta\rangle_{L^2(\Omega)}
  := \int_E (\theta_{\bar{u}}(T_1)-\theta_d)\Theta \, \dd x.\] Due to
  the choice of $s$, we have $X_s \embeds C(\olom) \embeds L^2(\O)$,
  such that the functional is also an element of $X_s'$ and
  $\delta_{T_1}\otimes\chi_E(\theta_{\bar{u}}(T_1)-\theta_d) \in
  \X'$. Moreover, we set $\|\nabla\theta_{\bar{u}}\|_{L^q}^{s-q}
  \Delta_q\theta_{\bar{u}}$ as a functional on $\X \embeds
  L^{s}(J;W^{1,q})$ via
  \[\left\langle \|\nabla\theta_{\bar{u}}\|_{L^q}^{s-q}
    \Delta_q\theta_{\bar{u}},\xi\right\rangle_{\X} := \int_J
  \|\nabla\theta_{\bar{u}}(t)\|_{L^q}^{s-q} \langle
  \Delta_q\theta_{\bar{u}}(t),\xi(t)\rangle_{W^{1,q}}\, \dd t.\] The
  inclusion $\X \embeds C(\olQ)$ also implies $\bar\mu \in \MM(\olQ)
  \embeds \X'$.  Hence, inserting $\theta'_{\bar{u}} = \s_\EE'(u)h$
  in~\eqref{e-lagrangederivative}, we immediately obtain
  \begin{align*}\partial_{\bar{u}} \fL(u,\mu)h & = \langle
    \s_\EE'(u)^*\!\left[\delta_{T_1}^*\chi_E(\theta_{\bar{u}}(T_1)-\theta_d)
      + \gamma\|\nabla\theta_{\bar{u}}\|_{L^q}^{s-q} \Delta_q
      \theta_{\bar{u}} + \mu\right],h\rangle_{\U} \\ & \qquad + \beta
    \int_{\Sigma_N} \partial_tu\partial_t h + \frac{p}{2}|u|^{p-2}uh
    \, \dd\omega\,\dd t\end{align*} for $h \in \U$. Let us introduce
  $(\vartheta,\chi)$ as the unique solution
  of~\eqref{e-abstractadjoint} (cf.\
  Theorem~\ref{t-weakadjointexists}) with data $\vartheta_T =
  \chi_E(\theta_{\bar{u}}(T_1)-\theta_d)+\bar\mu_{T_1}$, $g = 0$ and
  $f = \gamma\|\nabla\theta_{\bar{u}}\|_{L^q}^{s-q} \Delta_q
  \theta_{\bar{u}} + \bar\mu_{(T_0,T_1)}$, which is then also the
  solution of the formal system~\eqref{e-adjoint} with the stated
  inhomogeneities $f$ and $g$ and terminal value $\vartheta_T$. Here,
  $\psi$ is obtained by $\psi(\vartheta) =
  \JJ(\sigma(\theta_{\bar{u}}))^*(-\nabla\cdot\sigma(\theta_{\bar{u}})\vartheta\varepsilon\nabla\varphi_{\bar{u}})$,
  cf.~Remark~\ref{r-adjointtwoequations}.
  Corollary~\ref{c-adjointsolop} now shows that
  \begin{equation}\partial_{\bar{u}} \fL(\bar u,\bar \mu)h =
    \langle \EE^*\psi,h\rangle_{\U} + \beta
    \int_{\Sigma_N} \partial_t\bar u\partial_t h + \frac{p}{2}|\bar u|^{p-2}\bar
    u h \, \dd
    \omega\,\dd t\label{e-lagrangederivativefinal}\end{equation} for
  $h \in \U$. It is convenient to write $\EE$ as $\EE = \tau_\GN^*
  \circ \mathfrak E$ with $\mathfrak E \colon \U \embeds L^{2s}(J;L^{\mathfrak
    p}(\GN))$ and $\tau_\GN^* \colon L^{2s}(J;L^{\mathfrak
    p}(\GN)) \to L^{2s}(J;W^{-1,q}_\GD)$ with $\mathfrak p > \frac23 q$, 
  see Proposition~\ref{p-compactcontrolspace} and Remark~\ref{r-einbett/ident}. Then
  we have
  \begin{equation}\langle \EE^*\psi,h\rangle_{\U} = \langle \tau_\GN
    \psi, \mathfrak Eh
    \rangle_{L^{(2s)'}(J;L^{\mathfrak p'}(\GN)),L^{2s}(J;L^{\mathfrak p}(\GN))} 
    = \int_{\Sigma_N} (\tau_\GN\psi) h\,
    \dd\omega\,\dd t,\label{e-adjidentity}\end{equation} again $h \in \U$. 
  Inserting~\eqref{e-adjidentity}
  and~\eqref{e-lagrangederivativefinal} into~\eqref{e-lagrangemult3},
  we obtain the stated variational inequality.  \hfill\end{proof}

\begin{rem}
  If the optimal control $\bar u$ in the previous theorem is an
  \emph{interior} point of $\U^\ad$, or if $\U^\ad$ is not present at
  all, then one may transform the variational inequality~\eqref{eq:vi} to the
  ordinary nonlinear differential equation of order
  two \[\partial_{tt}\bar u = \frac{1}{\beta} \tau_\GN \psi +
  \frac{p}{2}|\bar u|^{p-2}\bar u\] in the space $L^{p'}(\GN)$ as a
  boundary value problem with $\partial_t \bar u(T_0) = \partial_t
  \bar u(T_1) = 0$. In particular, $\partial_{tt} \bar u \in
  L^{(2s)'}(J;L^{p'}(\GN))$ in this case.
\end{rem}
    
\section{Application and numerical example}\label{sec:numerics}

As already outlined in~\cite{h/m/r/r} and the introduction, a typical example of an
application for a problem in the form~\eqref{P} is the optimal
heating of a conducting material such as steel by means of an
electric current. The aim of such procedures is to heat up a
workpiece by electric current and to cool it down
rapidly with water nozzles in order to harden it. In case of steel,
this treatment indeed produces a hard martensitic outer
layer, see for instance
\cite[Ch.~9.18]{Cal07} for a phase diagram and~\cite[Chapters
10.5/10.7 about~Martensite]{Cal07}, and is thus used for instance for rack-and-pinion
actuators, to be found e.g.\ in steering mechanisms.  The part of the
workpiece to be heated up corresponds to the
design area $E$ in the objective functional in~\eqref{P}. In order to
avoid thermal stresses in the material, it is
crucial to produce a homogeneous temperature distribution in the
design area, which is reflected by the first term of the objective
functional. The gradient term in the objective functional further
enforces minimal thermal stresses. Moreover, the temperatures necesssary for the
hardening process as described above are rather close to the melting
point of the material, thus the state constraints are used to prevent the
temperature exceeding the melting temperature $\theta_{\max}$. The
control constraints in~\eqref{P} represent a maximum electrical current
which can be induced in the workpiece.
\begin{figure}[ht]
  \centering
  \begin{subfigure}{0.95\textwidth}
    \centering
    \includegraphics[width=\textwidth]{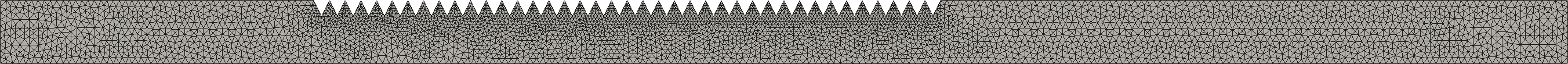}
    \caption{$\Omega$ with underlying mesh from the side
      ($x_1x_2$-plane).}
  \label{fig:meshgrid}
  \end{subfigure}
  \begin{subfigure}{0.95\textwidth}
    \centering
    \includegraphics[width=\textwidth]{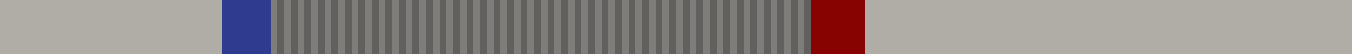}
    \caption{$\Omega$ from above ($x_1x_3$-plane) with $\GN$ (left)
      and $\GD$ (right) emphasized.}
  \label{fig:meshomega}
  \end{subfigure}
  \caption{The computational domain $\Omega$ used in the numerical
  example.}
  \label{fig:mesh}
\end{figure}

In the following we exhibit numerical examples for the optimal control of the three-dimensional thermistor
problem in the form~\eqref{P}, underlining in particular the importance of the
state-constraints. The considered computational domain $\Omega$ is a (simplified) three-dimensional gear-rack
as seen in Figure~\ref{fig:mesh}, where the design area $E$
consists of the sawteeth. The mesh consists of about 80000 nodes,
inducing 400000 cells with cell diameteres ranging from
$8.8\cdot10^{-4}$ to $7.6\cdot10^{-3}$.

The heat-equation we use in the computations is as follows:
\[\varrho C_p \partial_t \theta - \dive(\eta(\theta)\kappa\nabla\theta) = (\sigma(\theta)\nabla\varphi)\cdot\nabla\varphi.\]
It deviates from~\eqref{e-strongparabolic} by the factor $\varrho C_p$, the
so-called the volumetric heat capacity, where $\varrho$ is the density
of the material and $C_p$ is its specific heat capacity. However, since we assume $\varrho C_p$
to be constant, it certainly has no influence on the theory presented
above. In~\cite[Remarks~6.13/15]{hal/reh} and~\cite{hie/reh} it is
laid out how to modify the analysis if one wants to incorporate a
volumetric heat capacity depending on the temperature $\theta$. For a realistic modeling of the
process, we use the data gathered in~\cite{CEA93}, i.e., the workpiece
$\Omega$ is supposedly made of non-ferromagnetic stainless steel
(\#1.4301). The constants used can be found
in Table~\ref{tab:parameters} and the conductivity functions are given by
\[\sigma(\theta) := \frac{1}{a_\sigma + b_\sigma\theta + c_\sigma\theta^2 +
  d_\sigma\theta^3} \quad \text{for} \quad \theta \in [0,10000]~\text{K},\]
with the constants $a_\sigma =
4.9659\cdot 10^{-7}$, $b_\sigma = 8.4121 \cdot 10^{-10}$, $c_\sigma =
-3.7246\cdot10^{-13}$ and $d_\sigma = 6.1960\cdot10^{-17}$ for the
electrical conductivity (in $\Omega^{-1}\text{m}^{-1}$), and \[\eta(\theta) := 100 (a_\eta + b_\eta
\theta) \quad \text{for} \quad \theta \in [0,10000]~\text{K}\] with $a_\eta =
0.11215$ and $b_\eta = 1.4087\cdot10^{-4}$ for the
thermal conductivity (in $\text{Wm}^{-1}\text{K}^{-1}$). Both functions are extended outside of
$[0,10000]$ in a smooth and bounded way, such that Assumptions~\ref{a-assu1} and~\ref{a-sigmanemytskii} are satisfied. Note that $\varepsilon$ and $\kappa$ are
each chosen as the identity matrix, as we do not account for
heterogeneous materials in this numerical example. To counter-act on
the different scales inherent in the problem, cf.\ the value for
$u_{\max}$ and $\theta_0$ in Table~\ref{tab:parameters}, the model was
nondimensionalized for the implementation.
\begin{table}[ht]
  \centering
  \renewcommand{\arraystretch}{1.5}
  \begin{tabular}[t]{@{\extracolsep{\fill}}*{8}{c}}
   $\varrho$ & $C_p$ & $\alpha$ & $\theta_0$ &
   $\theta_l$ & $\theta_d$ & $\theta_{\max}$ & $u_{\max}$ \\ \hline\hline
   7900 $\frac{\text{kg}}{\text{m}^3}$ & 455 $\frac{\text{J}}{\text{kg\,K}}$    
& 20 $\frac{\text{W}}{\text{m}^3\,\text{K}}$ & 290 K & 290 K & 1500 K
& 1700 K & $10\cdot 10^7 \frac{\text{A}}{\text{m}^2}$ \\ \hline
  \end{tabular}
  \caption{Material parameters used in the numerical tests}
  \label{tab:parameters}
\end{table}

The optimization problem~\eqref{P} is solved by means
of a Nonlinear Conjugate-Gradients Method in the form as described in~\cite{daiyuan}, modified to a projected method to account for the
admissible set $\UU_{\text{ad}}$. The method needed up to 150 iterations to
meet the stopping criterion, which required the relative change in the
objective functional to be smaller than $10^{-5}$. The state
constraints in~\eqref{P} are incorporated by a quadratic penalty
approach---so-called Moreau-Yosida regularization---, cf.~\cite{HK06} and
the references therein, where the penalty-parameter was increased up to a
maximum of $10^{10}$, stopping earlier if the violation of the state
constraints was smaller than $10^{-2}$~K. This resulted in a violation
of $9.54\cdot10^{-2}$~K, which is about 0.0056$\%$ of the upper bound
of 1700~K. In each step of the
optimization algorithm, the nonlinear state equations~\eqref{e-strongparabolic}-\eqref{e-strongellipticdirichlet} and
the adjoint equations~\eqref{e-adjoint} have to be solved. We use an
Implicit Euler Scheme for the time-discretization of these equations, whereas the spatial
discretization is done via piecewise continuous linear finite
elements. The nonlinear system of equations arising in each time-step is
solved via Newton's method. Here, we do a semi-implicit pre-step to
obtain a suitable initial guess for the discrete $\varphi$ for
Newton's method. For the control, we also choose piecewise continuous
linear functions in space where the values in the first and last
timestep were pre-set to $0$. In the calculation of the gradient of the
reduced objective functional $j$, the gradient representation with respect to the
$L^2(J;L^2(\GN))$ scalar product of the derivative of $u \mapsto
\frac12(\partial_t u)^2$ is needed, which one formally computes as
$\partial^2_{tt} u$. We used the
second order central difference quotient
$\frac{u_{k+1}-2u_k-u_{k-1}}{\Delta t^2}$ to approximate
$(\partial^2_{tt} u)(t_k)$ at time step $k$ with the appropriate
modifications for the first and last time step, respectively. All
 computations were performed within the FEnICS framework~\cite{FEnICS}.

For the experiment duration, we set $T_1-T_0 = 2.0$~s with
timesteps $\Delta t = 0.02$~s and $T_0 = 0.0$~s, while we use $\gamma
= 10^{-8}$ and $\beta =
10^{-5}$ -- this small value for $\beta$ is only possible due to the
nondimensionalization performed. In the following, we elaborate on two
settings: one in which we enforce the state constraint $\theta \leq
\theta_{\max}$ and one in which we do not. 
\begin{figure}[ht]
  \centering
  \begin{subfigure}{0.425\textwidth}
    \centering
    \includegraphics[width=\textwidth]{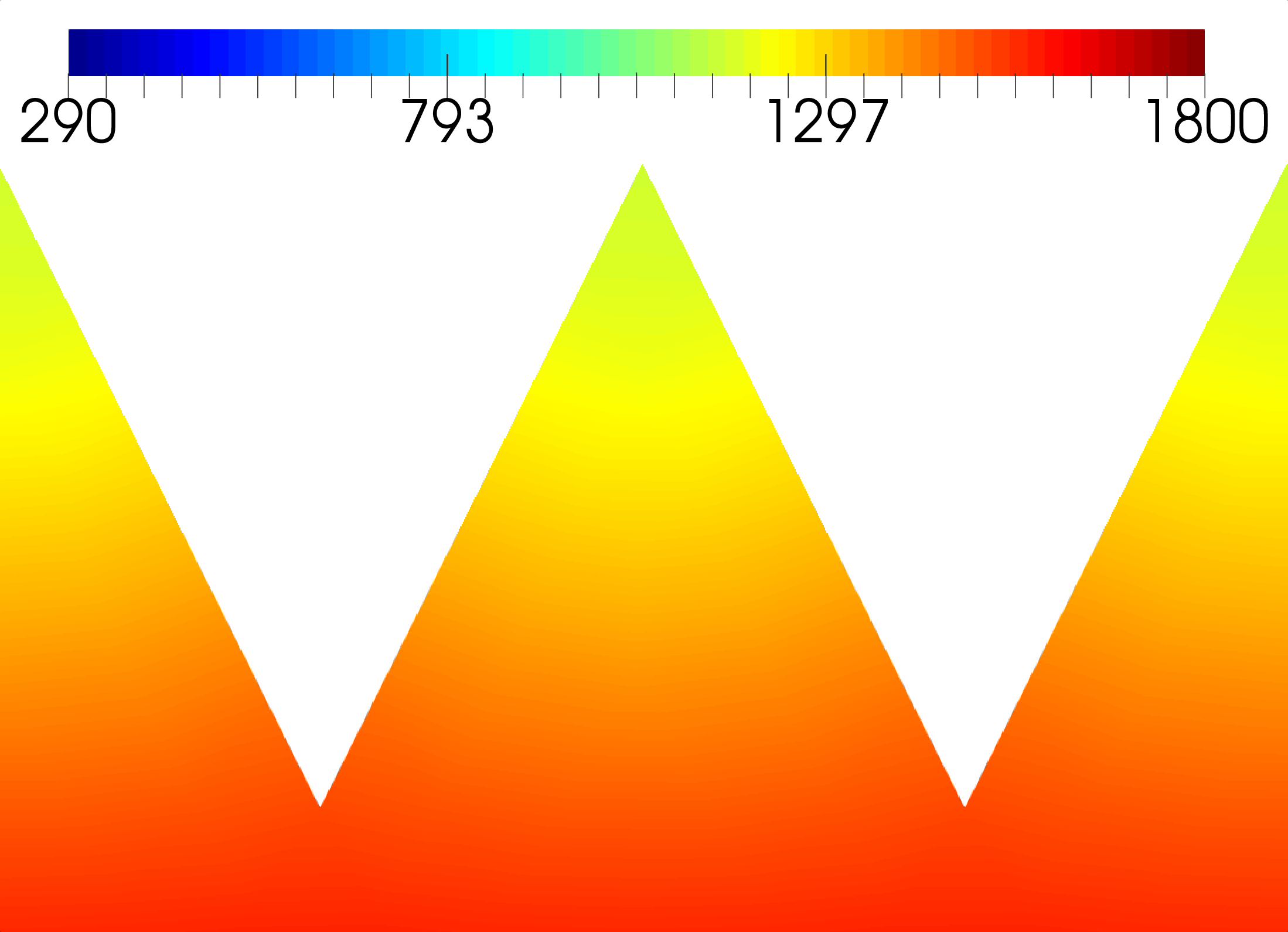}
    \caption{Free optimization.}
    \label{fig:objectivefree}
  \end{subfigure} \qquad
  \begin{subfigure}{0.425\textwidth}
    \centering
    \includegraphics[width=\textwidth]{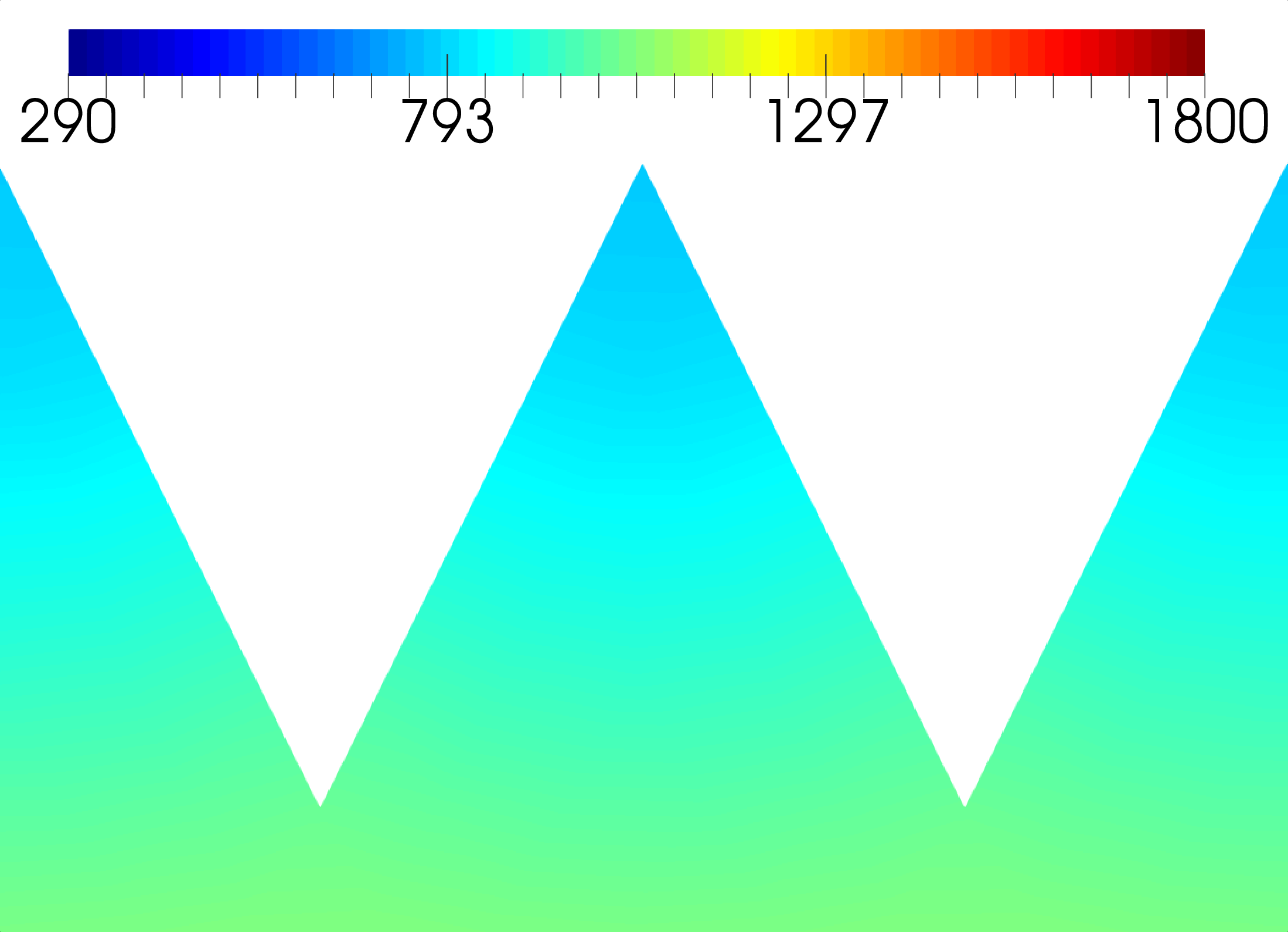}
    \caption{State constrained optimization.}
    \label{fig:objectivesc}
  \end{subfigure}
  \caption{Detail of the sawteeth in $E$ at end time $t = 2.0$~s with
    distribution of the temperature $\theta$ in~K.}
  \label{fig:objective}
\end{figure}

\begin{figure}[ht]
  \centering
  \begin{subfigure}{0.425\textwidth}
    \centering
    \includegraphics[width=\textwidth]{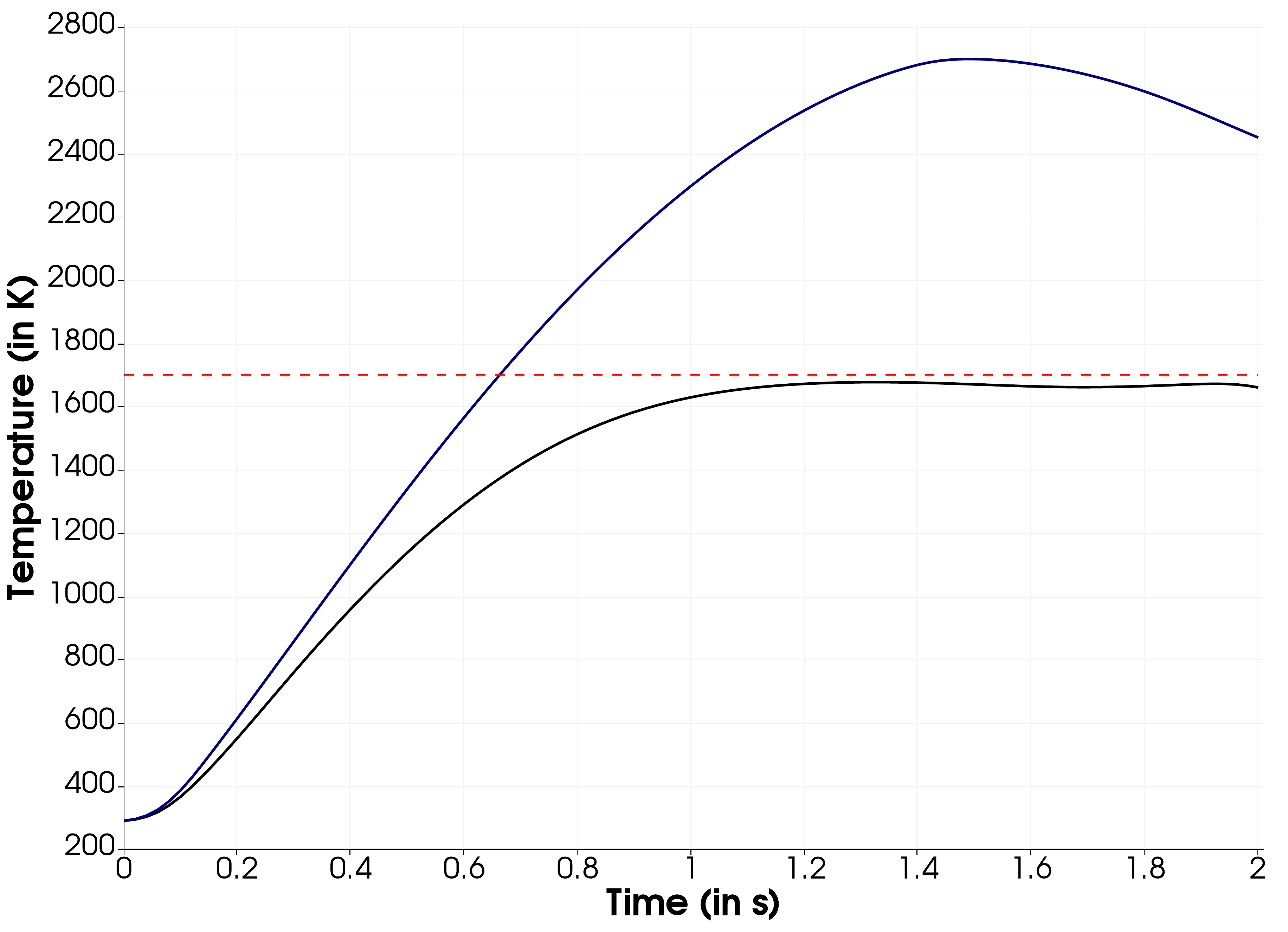}
    \caption{Time plot of the temperature in a point close to
      $\GN$.}
    \label{fig:violatedboundsplot}
  \end{subfigure} \qquad
  \begin{subfigure}{0.425\textwidth}
    \centering
    \includegraphics[width=\textwidth]{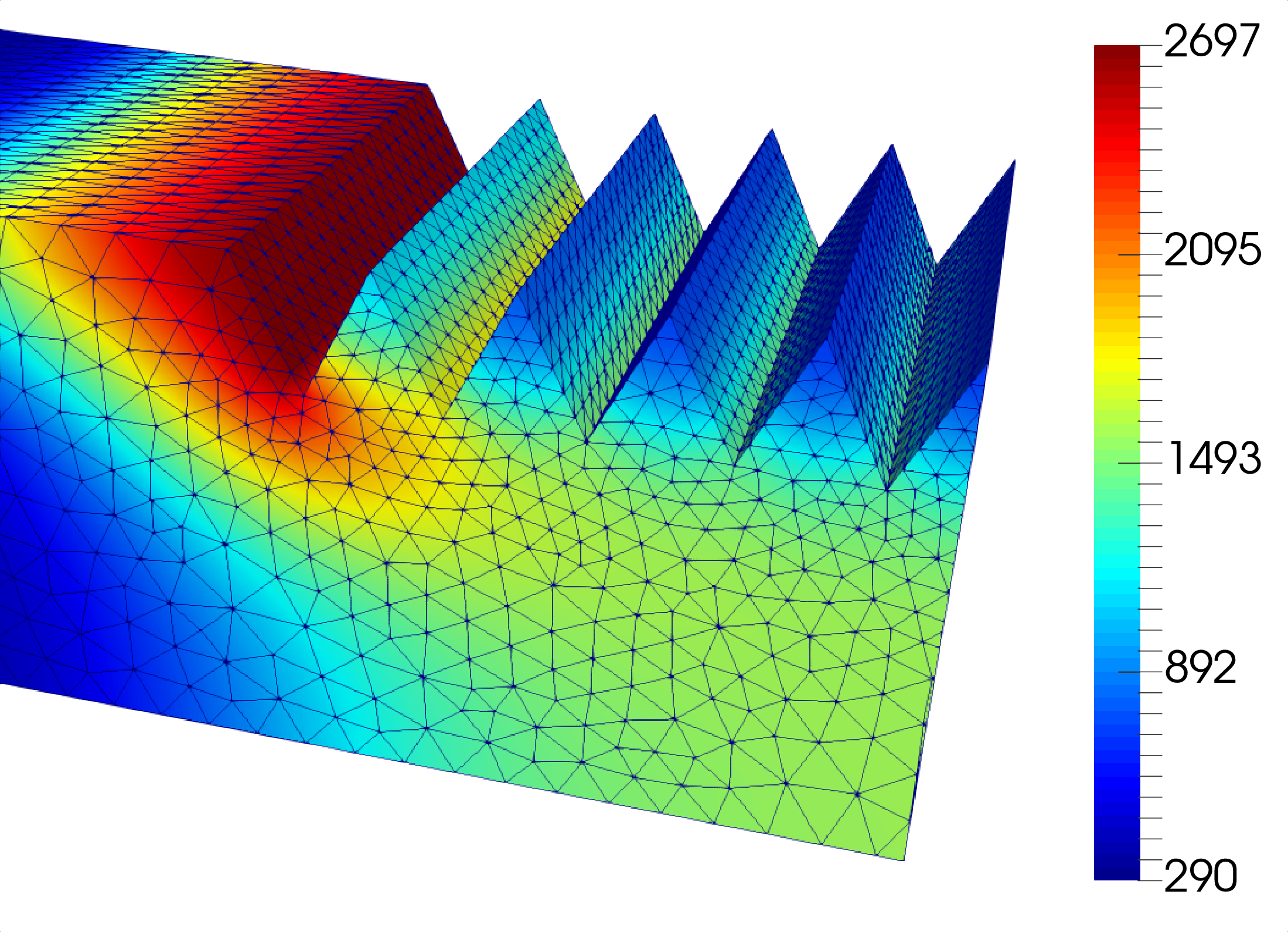}
    \caption{Temperature $\theta$ in K at the critical area near $\GN$ at time $t
      = 1.4$ s.}
    \label{fig:violatedboundscritical}
  \end{subfigure}
  \caption{Influence and necessity of state constraints.}
  \label{fig:violatedbounds}
\end{figure}

Figure~\ref{fig:objective} shows the temperature
distribution at end time $T_1 = 2.0$~s in $E$ in both cases. The desired
temperature distribution close to uniformly 1500~K has been nearly
achieved in the free optimization, see Figure~\ref{fig:objectivefree},
at the price of very high temperature values around $\GD$ and
$\GN$ already early in the heating process. We come back to this below, cf.\ also
Figure~\ref{fig:violatedbounds}. For the state-constrained
optimization, we achieve a much worse result (note the same scales in
both Figure~\ref{fig:objectivefree} and~\ref{fig:objectivesc}), which again corresponds
to the rapid evolution to high temperatures at the critical areas,
since these crucially limit the maximal amount of energy induced into the
workpiece if one wants to prevent the temperature rising higher than
the given bounds $\theta_{\max}$. This can also be seen in
the development of the optimal controls in both cases over time, see below.

\begin{figure}[htbp]
  \centering
  \includegraphics
  [width=0.95\textwidth]{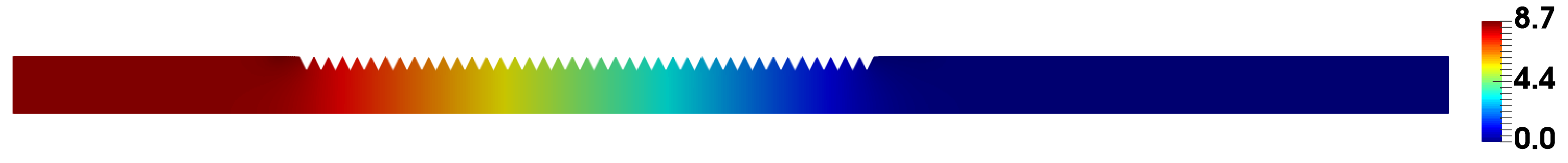}
  \caption{The potential $\varphi$ (in V) associated with the optimal
    solution at time $t = 1.0$~s, view from the side
    ($x_1x_2$-plane).}  \label{fig:pot}
\end{figure}

\begin{figure}[htbp]
  \centering
    \includegraphics
    [width=0.95\textwidth]{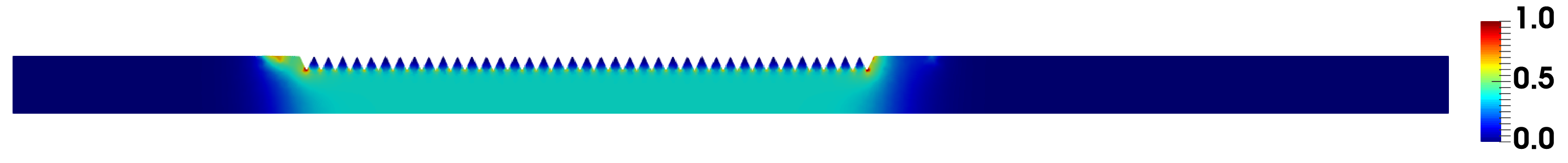}
    \caption{Magnitude of the gradient $\nabla \varphi$ (in V/m) associated with the optimal
    solution at time $t = 1.0$~s, view from the side
    ($x_1x_2$-plane).}
    \label{fig:potgrad}
\end{figure}
The potential $\varphi$ and its gradient $\nabla \varphi$ associated with the optimal control to the
state-constrained optimization problem, at time $t = 1.0$~s are
depicted in Figures~\ref{fig:pot} and~\ref{fig:potgrad}. Here, $\nabla\varphi$ is to be
understood as the
projection of the potentially discontinuous gradient of $\varphi$ to the space of continuous
linear finite elements. The potential $\varphi$ decreases from $\GN$ to the grounding with prescribed value $\varphi \equiv 0$ at
$\GD$, cf.\
Figure~\ref{fig:meshomega}, thus inducing a current flow and acting as
a heat source between $\GD$ and $\GN$, since the corresponding term in
the heat equation
$\sigma(\theta)\varepsilon\nabla\varphi\cdot\nabla\varphi$ is proportional to
$|\nabla\varphi|^2$ due to the coercivity and boundedness of
$\varepsilon$ and the bounds on $\sigma$. This is confirmed by the
magnitude of $\nabla \varphi$ as seen in Figure~\ref{fig:potgrad}. In
particular one observes that $\nabla \varphi$ is very small or 0 in
$E$, which means that the current flows only through the area between $\GD$ and
$\GN$ and right \emph{below} $E$, heating only this part of the workpiece.

\begin{figure}[ht]
  \centering
  \begin{subfigure}{0.425\textwidth}
    \centering
    \includegraphics[width=\textwidth]{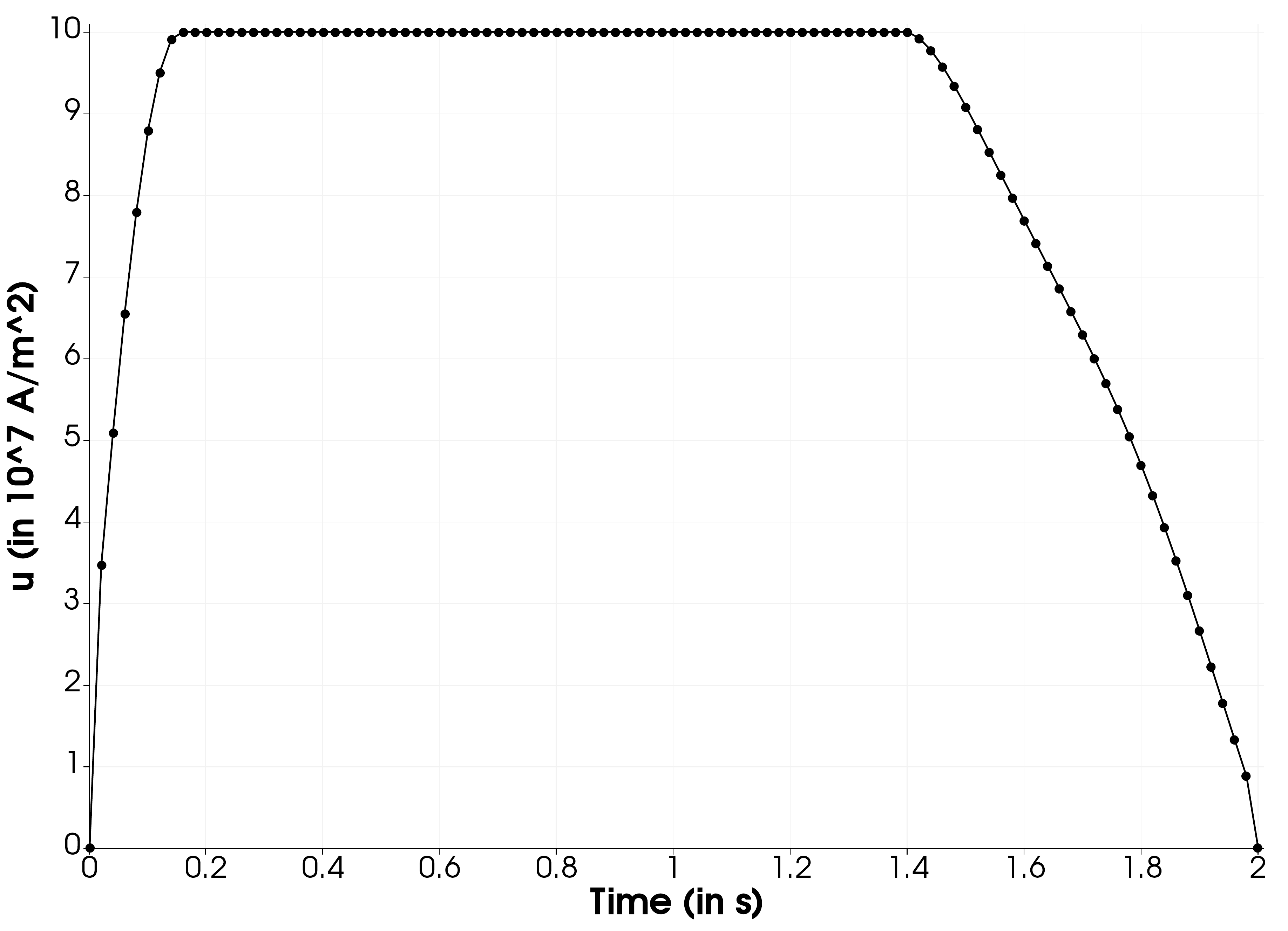}
    \caption{For the unconstrained problem.}
  \end{subfigure} \qquad
  \begin{subfigure}{0.425\textwidth}
    \centering
    \includegraphics[width=\textwidth]{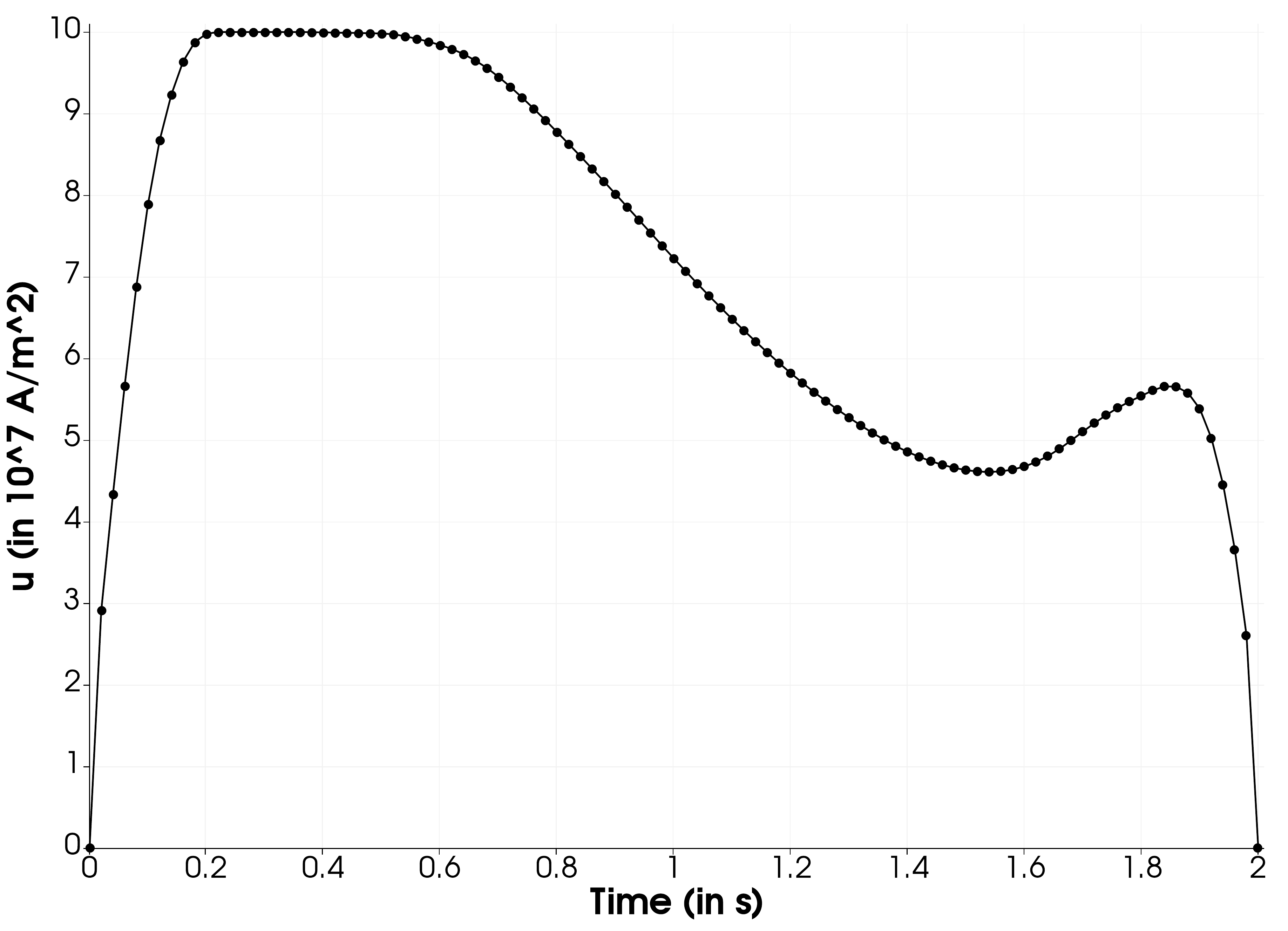}
    \caption{For the state-constrained problem.}
  \end{subfigure}
  \caption{Time plot of the optimal controls, taken at an arbitrary
    but fixed grid point in $\GN$.}
  \label{fig:optimalcontrols}
\end{figure}
The optimal controls are shown in Figure~\ref{fig:optimalcontrols}, taken at an
arbitrary but fixed grid point in $\GN$. The high
values in the control at the beginning of the process seem to be the
result of the inability to heat up the tooth rack in the design-area $E$
directly as explained above, which makes heating of
the teeth reliant on diffusion. This in turn requires the needed total
energy to be inserted into the system as fast as possible, resulting
in high control values, which also agrees with the requirement to obtain a
\emph{uniform} temperature distribution in the tooth rack. These considerations also underline the
necessity of control bounds in this example. In decreasing the control values after the
inital period, the opimization procedure in the free optimization is avoiding to
``over-shoot'', i.e., to produce a higher temperature than desired.
In the case of state-constrained optimization, the
presence of the state constraints forces an earlier decrease in control values in
order to not violate the upper bound $\theta_{\max}$, which is then
compensated by a slightly higher level of values towards the end of the
simulation. This, however, is clearly not enough to make up for the
earlier decrease as seen in Figure~\ref{fig:objective}.

Figure~\ref{fig:violatedbounds} illustrates why state constraints are
a necessary addition to an appropriate model of the industrial
steel heating process. Figure~\ref{fig:violatedboundsplot} shows the
temperature evolution in a point in one of the two critical regions,
which are the points near $\GD$ and $\GN$, see also
Figure~\ref{fig:violatedboundscritical} and the magnitude of
$\nabla\varphi$ at this region in Figure~\ref{fig:potgrad}. In this case, the point lies in $E$ close to $\GN$, but we emphasize that the state constraints
hold in the whole $\Omega$ and are not limited to $E$. The upper line
in Figure~\ref{fig:violatedboundsplot} corresponds to the temperature
associated to the optimal solution of the unconstrained optimization,
while the lower belongs to the state-constrained optimal solution,
with the upper bound $\theta_{\max} = 1700$~K marked by the dashed
line. In the free optimization case, the temperature exceeds the
bounds already at about one third of the
simulation time and continuous to rise to almost 1000~K above
$\theta_{\max}$. On the
other hand, the temperature obtained from the state-constrained case stays below the threshold, as required. Note here that
the evaluated point is chosen as one of those where the temperature
rises highest overall, compare the temperature distribution as seen in
Figure~\ref{fig:violatedboundscritical} and the maximal temperature
achieved in the free optimization case in
Figure~\ref{fig:violatedboundsplot}.

Concluding from the results presented above, it becomes apparent that
the prescribed time of 2.0~s is too short to heat up the workpiece in
the given geometry enough to reach the required temperature for Austenite to form in the
workpiece (cf.~\cite[Ch.~9.18]{Cal07}) in $E$, if melting is to be
prevented.

\appendix

\section{A ``minimum principle''}

\begin{proposition}\label{p-lowerboundsappendix} For every solution
  $(\theta,\varphi)$ in the sense of Theorem~\ref{t-locexist} with
  maximal existence interval $J_{\max}$, it is true that $\theta(x,t)
  \geq \min(\essinf_\Si \theta_l,\min_{\olom} \theta_0)$
  for all $(x,t) \in \olom \times [T_0,T_\bullet]$, where $T_\bullet
  \in J_{\max}$.
\end{proposition}

\begin{proof}
  We set $m_{\inf} :=  \min(\essinf_\Si \theta_l,\min_{\olom} \theta_0)$ and $\zeta(t) = \theta(t) - m_{\inf}$ and decompose $\zeta(t)$
  into its positive and negative part, that is, $\zeta(t) =
  \zeta^+(t)-\zeta^-(t)$ with both $\zeta^+(t)$ and $\zeta^-(t)$ being
  positive functions. By~\cite[Ch.~IV,~\S7, Prop.~6/Rem.~12]{daut/lionsII} we then have that
  $\zeta^-(t)$ is still an element of $W^{1,q}$ for almost every $t
  \in (T_0,T_\bullet)$. In particular, we may test~\eqref{e-paraex} against
  $-\zeta^-(t)$, insert $\theta = \zeta + m_{\inf}$ and use that
  $m_{\inf}$ is constant:
  \begin{multline*}-\int_\O \partial_t\zeta(t) \zeta^-(t) \, \dd x -
    \int_\O (\eta(\theta(t))\kappa\nabla\zeta(t))\cdot\nabla\zeta^-(t)
    \, \dd x - \int_\Gamma \alpha\zeta(t)\zeta^-(t) \, \dd x \\ =
    -\int_\Gamma \alpha(\theta_l(t)-m_{\inf})\zeta^-(t) - \int_\O
    \zeta^-(t)(\sigma(\theta(t))\varepsilon\nabla\varphi(t))\cdot\nabla\varphi(t)
    \, \dd x.\end{multline*} Observe that the support of products of
  $\zeta(t)$ and $\zeta^-(t)$ is exactly the support of $\zeta^-(t)$,
  and $\zeta(t) = -\zeta^-(t)$ there. We thus obtain
  \begin{multline}\frac{1}{2}\, \partial_t\!\left\|\zeta^-(t)\right\|_{L^2}^2
    + \int_\O
    (\eta(\theta(t))\kappa\nabla\zeta^-(t))\cdot\nabla\zeta^-(t) \,
    \dd x + \int_\Gamma \alpha\zeta^-(t)^2 \, \dd x \\ = -\int_\Gamma
    \alpha(\theta_l(t)-m_{\inf})\zeta^-(t) - \int_\O
    \zeta^-(t)(\sigma(\theta(t))\varepsilon\nabla\varphi(t))\cdot\nabla\varphi(t)
    \, \dd x.\label{e-lowerbounds}\end{multline} Let us show that
  $\partial_t\|\zeta^-(t)\|_{L^2}^2 \leq 0$. By
  Assumption~\ref{a-assu1},
  $(\eta(\theta(t))\kappa\nabla\zeta^-(t))\cdot\nabla\zeta^-(t) \geq
  \underline{\eta}\underline{\kappa}\|\nabla\zeta^-(t)\|_{\R^3}^2$ and
  $-(\sigma(\theta(t))\varepsilon\nabla\varphi(t))\cdot\nabla\varphi(t)
  \leq -
  \underline{\sigma}\underline{\varepsilon}\|\nabla\varphi(t)\|_{\R^3}^2$.
  This means that both integrals on the left-hand side
  in~\eqref{e-lowerbounds} are positive (since $\alpha \geq 0$), while
  the second term on the right-hand side is negative. The constant $m_{\inf}$ is
  constructed exactly such that $\theta_l(t) - m_{\inf}$ is greater or
  equal than zero almost everywhere, such that $-\alpha
  (\theta_l(t)-m_{\inf})\zeta^-(t) \leq 0$. Hence,
  from~\eqref{e-lowerbounds} it follows that $\partial_t
  \|\zeta^-(t)\|_{L^2}^2 \leq 0$. But, due to the construction of
  $\zeta$, we have $\zeta(T_0) \geq 0$, which means that $\zeta^-(T_0)
  \equiv 0$ and thus $\zeta^-(t) \equiv 0$ for all $t \in (T_0,T_\bullet)$.
  \hfill\end{proof}




\end{document}